\documentclass[11pt,a4paper]{article}
\usepackage{amsfonts}

\usepackage{indentfirst}
\usepackage{cite}
\usepackage{amsmath,amsthm}
\allowdisplaybreaks[2]
\usepackage{amssymb}
\usepackage{mathrsfs}

\usepackage[hidelinks,bookmarksnumbered,hypertexnames=false]{hyperref}
\hypersetup{pdftitle={Global strong solutions for the two-dimensional full compressible Navier--Stokes equations with vacuum at infinity}}
\usepackage[
  hmargin=2.2cm,
  top=2.4cm,
  bottom=2.6cm,
  footskip=0.9cm
]{geometry}

\newtheorem{theorem}{Theorem}[section]
\newtheorem{remark}{Remark}[section]
\newtheorem{definition}{Definition}[section]
\newtheorem{lemma}[theorem]{Lemma}

\newtheorem{proposition}[theorem]{Proposition}

\makeatletter
\@addtoreset{equation}{section}
\makeatother

\renewcommand{\div}{ {\rm div }  }
\newcommand{\na}{\nabla }

\newcommand{\OM}{\Omega}

\newcommand{\la}{\label}

\newcommand{\ba}{\begin{aligned}}
\newcommand{\ea}{\end{aligned}}
\newcommand{\be}{\begin{equation}}
\newcommand{\ee}{\end{equation}}

\def\rr{\mathbb{R}^2}

\title{Global Strong Solutions to the Cauchy problem of the Two-Dimensional Full Compressible Navier--Stokes Equations with Degenerate Heat Conductivity and Large Initial Data}

\author{Qinghao L{\small EI}$^{a,b}$, Xin S{\small{ONG}}$^{c}$ \thanks{Email addresses: leiqinghao22@mails.ucas.ac.cn (Q. H. Lei), songxin2357@email.ncu.edu.cn (X. Song)} \\
{\normalsize a. School of Mathematical Sciences,}\\
{\normalsize University of Chinese Academy of Sciences, Beijing 100049, P. R. China;} \\
{\normalsize b. Institute of Applied Mathematics,} \\
{\normalsize Academy of Mathematics and Systems Science,} \\
{\normalsize Chinese Academy of Sciences, Beijing 100190, P. R. China;} \\
{\normalsize c. School of Mathematics and Computer Sciences,}\\
{\normalsize  Nanchang University, Nanchang 330031, P. R. China.}
}

\date{}
\begin{document}
\maketitle

\begin{abstract}
The Cauchy problem for the two-dimensional full compressible Navier--Stokes equations with vacuum far-field density is considered.
Assuming that the shear viscosity $\mu$ is a positive constant, the bulk viscosity $\lambda(\rho)\ge0$ depends on the density, and the heat conductivity takes the power-law form $\kappa(\theta) = \kappa_0\theta^\alpha$ with $\kappa_0, \alpha > 0$, we establish the global existence of strong solutions for a broad class of density-dependent bulk viscosities.
In particular, our results cover the sublinear growth law $\lambda(\rho) = \rho/\log\log(e^e + \rho)$ for $\alpha \ge 1$, as well as all superlinear power laws for every $\alpha > 0$.
These results hold without any smallness restrictions on the initial data, and the initial density is allowed to vanish.
\end{abstract}

\section{Introduction and main results}\label{sec:intro}

We consider a compressible, viscous, heat-conducting polytropic fluid in $\mathbb{R}^2$. Its motion is governed by the two-dimensional full compressible Navier--Stokes system:
\begin{equation}\label{intro:conservative}
\begin{cases}
 \rho_t + \operatorname{div}(\rho u) = 0, \\
 (\rho u)_t + \operatorname{div}(\rho u \otimes u) + \nabla P = \operatorname{div}\mathbb{S}, \\
 (\rho E)_t + \operatorname{div}(\rho E u + P u) - \operatorname{div}\bigl(\kappa(\theta)\nabla\theta\bigr) = \operatorname{div}(\mathbb{S}u),
\end{cases}
\end{equation}
where the specific total energy $E$, pressure $P$, and viscous stress tensor $\mathbb{S}$ are given by
\begin{equation}\nonumber
 E = c_v\theta + \frac{1}{2}|u|^2,
 \qquad
 P = R\rho\theta,
 \qquad
 \mathbb{S} = 2\mu D(u) + \lambda(\rho)(\operatorname{div} u)\,\mathbb{I}_2,
\end{equation}
with $D(u) = \frac{1}{2}\bigl(\nabla u + (\nabla u)^\top\bigr)$ denoting the deformation tensor and $\mathbb{I}_2$ the $2 \times 2$ identity matrix. Here $t \ge 0$ and $x \in \mathbb{R}^2$ are the time and spatial variables, respectively, while $\rho \ge 0$, $u = (u^1, u^2)$, $P$, and $\theta > 0$ represent the fluid density, velocity, pressure, and absolute temperature.

The physical parameters $\mu, c_v, R, \kappa_0$ are positive constants; without loss of generality, we normalize $R = 1$. For a fixed exponent $\alpha > 0$, we take the heat conductivity and its primitive to be
\begin{equation}\nonumber
 \kappa(z) = \kappa_0 z^\alpha,
 \qquad
 K(z) = \int_0^z \kappa(s)\,ds = \frac{\kappa_0}{1 + \alpha} z^{1 + \alpha}
 \qquad (z \ge 0).
\end{equation}
Throughout the paper, the density-dependent bulk viscosity $\lambda$ is assumed to satisfy, for some constant $C_\lambda > 0$,
\begin{equation}\label{main:lambda-derivative}
 \lambda \in C^1([0,\infty)) \cap C^2((0,\infty)),
 \qquad
 \lambda \ge 0,
 \qquad
 r|\lambda'(r)| \le C_\lambda\bigl(2\mu + \lambda(r)\bigr)
 \quad (r > 0).
\end{equation}
For strong solutions, system \eqref{intro:conservative} can be rewritten as
\begin{equation}\label{main:pde}
\begin{cases}
 \rho_t + \operatorname{div}(\rho u) = 0, \\
 \rho(u_t + u \cdot \nabla u) = \mu\Delta u + \nabla\bigl((\mu + \lambda(\rho))\operatorname{div} u\bigr) - \nabla P, \\
 c_v\rho(\theta_t + u \cdot \nabla\theta) - \Delta K(\theta) = 2\mu|D(u)|^2 + \lambda(\rho)(\operatorname{div} u)^2 - P\operatorname{div} u.
\end{cases}
\end{equation}
We supplement \eqref{main:pde} with the initial data and the vacuum far-field condition
\begin{equation}\label{intro:cauchy}
 (\rho, \rho u, \rho\theta)|_{t=0} = (\rho_0, \rho_0 u_0, \rho_0\theta_0),
 \qquad
 \lim_{|x|\to\infty} \rho(t,x) = 0 \quad (t \ge 0).
\end{equation}

The existence of solutions to the full compressible Navier--Stokes equations have been studied extensively.
In one dimension, Kazhikhov and Shelukhin~\cite{KazhikhovShelukhin1977} obtained the global existence and uniqueness of strong solutions to the initial-boundary value problem with constant transport coefficients and no vacuum.
Kazhikhov~\cite{Kazhikhov1982} later extended this result to the Cauchy problem for large smooth perturbations of a positive constant state.
With density-dependent viscosity and temperature-dependent heat conductivity, Duan, Guo, and Zhu~\cite{DuanGuoZhu2017} proved global well-posedness under stress-free, thermally insulated boundary conditions, assuming the initial density has a positive lower bound.
In higher dimensions, Serrin~\cite{Serrin1959} proved the uniqueness of classical solutions away from vacuum, and Nash~\cite{Nash1962} obtained local existence for the three-dimensional Cauchy problem with strictly positive initial density.
When vacuum is present, Cho and Kim~\cite{ChoKim2006} established local strong solutions under compatibility conditions on the initial momentum and thermal energy equations.
For initial data close in $H^3$ to a constant state with positive density and temperature, Matsumura and Nishida~\cite{MatsumuraNishida1980} proved the global existence and uniqueness of classical solutions in $\mathbb{R}^3$, while Hoff~\cite{Hoff1997} obtained global weak solutions with discontinuous initial data in two and three dimensions.
Regarding global weak solutions with large initial data, Feireisl~\cite{Feireisl2004} proved the existence of variational solutions for pressure laws that exclude the ideal polytropic gas, whereas Bresch and Desjardins~\cite{BreschDesjardins2007} established weak stability by invoking a structural relation coupling the density-dependent shear and bulk viscosities.
In the presence of both vacuum and large oscillations, Huang and Li~\cite{HuangLi2018} obtained global classical solutions in $\mathbb{R}^3$ for smooth data with small initial energy and strictly positive far-field density and temperature; subsequent extensions to bounded and exterior domains were carried out by Li, L\"u, and Wang~\cite{LLW} and Li, Li, and L\"u~\cite{LiLiLu2025}, respectively.
Under other smallness assumptions, Wen and Zhu~\cite{WenZhu2017} established the global existence of strong solutions for sufficiently small initial mass, while Li~\cite{Li2020} proved global well-posedness under a scaling-invariant smallness condition provided that $2\mu>\lambda$; this viscosity restriction was later removed by Wen~\cite{Wen2025}.
Moreover, in the case of density-dependent viscosity, Shen, Xu, and Zhang~\cite{ShenXuZhang2025} proved the global existence of strong solutions with initial vacuum to the Dirichlet initial-boundary value problem in three-dimensional bounded domains, provided that the initial energy is sufficiently small.

In the isentropic case with $\lambda(\rho)=\rho^\beta$, the system (\ref{main:pde}) corresponds to the model introduced by Vaigant-Kazhikhov \cite{VK}, for which the authors proved global existence and uniqueness of strong solutions in rectangular domains for large, non-vacuum initial data whenever $\beta>3$.
This result was later extended to the periodic torus with initial vacuum permitted by Jiu, Wang, and Xin~\cite{JWX1}.
By combining commutator estimates with blow-up criteria, Huang and Li~\cite{HL2,HL3} (see also Jiu, Wang, and Xin~\cite{JWX2}) weakened the exponent requirement from $\beta>3$ to $\beta>\frac{4}{3}$ in both periodic domains and the whole space, where the density is allowed to vanish.
Subsequently, Fan, Li, and Li~\cite{FLL} extended the global existence of strong solutions for $\beta>\frac{4}{3}$ to general bounded simply connected domains under Navier-slip boundary conditions, while Fan, Li, and Wang~\cite{FLW} established uniform-in-time upper bounds on the density along with exponential decay rates in both periodic and bounded simply connected domains.
Turning to non-isentropic flows, Fan and Huang~\cite{FanHuang} recently investigated the two-dimensional compressible Navier--Stokes equations with transported entropy and $\lambda(\rho)=\rho^\beta$ ($\beta>\frac{4}{3}$), proving global existence and large-time behavior of strong solutions for arbitrarily large initial data containing vacuum on both the torus and bounded simply connected domains; in the bounded-domain setting, they leveraged the Navier-slip boundary condition to establish a new commutator estimate.
More recently, Fan, Jiang, Li, and L\"u~\cite{FJLL1,FJLL2} considered a general class of density-dependent bulk viscosities, covering all superlinear power laws and near-linear laws with logarithmic gains, and established the global existence and large-time behavior of strong solutions in periodic domains, the whole space, and smooth bounded domains with Navier-slip boundary conditions.
In the linear case $\lambda(\rho)=\rho$, global existence and uniqueness of strong solutions for large initial data were established by Li, Liang, and Ye~\cite{LLY} in $\mathbb{R}^2$ with vacuum at infinity, and by Fan, Huang, and Lei~\cite{FHL} on $\mathbb{T}^2$.

For the full compressible Navier--Stokes equations \eqref{main:pde} with degenerate heat conductivity, Lei and Song~\cite{LeiSong} obtained the global existence and large-time behavior of strong solutions for large initial data in periodic and bounded domains, as well as in the whole space with non-vacuum far-field density.
The purpose of this paper is to establish the global existence of strong solutions to the Cauchy problem in $\mathbb{R}^2$ with vacuum far-field density for arbitrarily large initial data.

To formulate the structural conditions on $\lambda$, we set
\begin{equation}\nonumber
 h(r) = 2\mu + \lambda(r),
 \qquad
 L(r) = \log(e + r),
 \qquad
 J(r) = \log\log(e^e + r),
\end{equation}
and define the viscosity primitive
\begin{equation}\nonumber
 \Lambda(r) = \int_1^r \frac{h(s)}{s}\,ds,
 \qquad r > 0.
\end{equation}

We now specify the admissible classes of bulk viscosities and heat conductivities.

\begin{definition}[Admissible bulk-viscosity and heat-conductivity laws]\label{main:classes}
Assume \eqref{main:lambda-derivative} together with one of the following four conditions, where each pointwise lower bound on $h$ or $\mu + \lambda$ is required for all $r \ge r_*$ with some fixed $r_* \ge 2$:

\smallskip\noindent\textup{(A) Pointwise lower bound ($\alpha > 1$).}
There exist constants $c, C > 0$, an exponent $0 \le \beta_D < 1/2$, and a nondecreasing function $D : [2,\infty) \to [1,\infty)$ such that
\begin{equation}\nonumber
 h(r) \ge \frac{c r}{D(r)},
 \qquad
 D(r) \le C L(r)^{\beta_D},
 \qquad
 \int_2^\infty \frac{dr}{(e + r)L(r)D(r)} = \infty.
\end{equation}

\smallskip\noindent\textup{(B) Integral criterion ($\alpha > 1$).}
With $D_h(r) = 1 + \sup_{0 \le z \le r} \frac{z}{\mu + \lambda(z)}$, the bulk viscosity satisfies
\begin{equation}\nonumber
 \lim_{r\to\infty} \frac{\log D_h(r)}{\log L(r)} = 0,
 \qquad
 \int_2^\infty \frac{\Lambda(r)}{r^2 L(r)}\,dr = \infty.
\end{equation}

\smallskip\noindent\textup{(C) Logarithmic viscosity gain.}
There exist constants $c > 0$ and $0 \le d < 1$ such that
\begin{equation}\nonumber
 \mu + \lambda(r) \ge c r L(r)^d,
 \qquad
 \alpha \ge \frac{2(1 - d)}{2 + d}.
\end{equation}

\smallskip\noindent\textup{(D) Linear heat conductivity ($\alpha = 1$).}
There exist constants $c > 0$, $p \ge 0$, and a function $\Psi \in C^1([1,\infty))$ with $\Psi \ge 1$ such that
\begin{equation}\nonumber
 \mu + \lambda(r) \ge \frac{c r}{\Psi(J(r))},
\end{equation}
together with
\begin{equation}\nonumber
 0 \le \frac{x\Psi'(x)}{\Psi(x)} \le p \quad (x \ge 1),
 \qquad
 \int_2^\infty \frac{dx}{\Psi(x)\log(e + x)} = \infty.
\end{equation}
\end{definition}

\begin{remark}[Representative coefficient laws]
\normalfont
Let $\lambda_0 \ge 0$ and $\kappa(\theta) = \kappa_0\theta^\alpha$ with $\kappa_0 > 0$. Definition~\ref{main:classes} covers the following representative growth laws and ranges of $\alpha$:
\begin{center}
\renewcommand{\arraystretch}{1.35}
\begin{tabular}{@{}l@{\hspace{1.5em}}c@{\hspace{1.5em}}c@{}}
\hline
Bulk viscosity $\lambda(r)$ & Admissible range of $\alpha$ & Case \\
\hline\noalign{\vskip 4pt}
$\lambda_0 + \dfrac{r}{[\log\log(e^e+r)]^\eta}$, $0 < \eta \le 1$ & $\alpha \ge 1$ & \textup{(A)}, \textup{(D)} \\[7pt]
$\lambda_0 + r$ & $\alpha \ge 1$ & \textup{(A)}, \textup{(D)} \\[3pt]
$\lambda_0 + r[\log(e+r)]^d$, $0 < d < 1$ & $\alpha \ge \dfrac{2(1-d)}{2+d}$ & \textup{(C)} \\[6pt]
$\lambda_0 + r[\log(e+r)]^d$, $d \ge 1$ & $\alpha > 0$ & \textup{(C)} \\[2pt]
$\lambda_0 + r\log(e+r)[\log\log(e^e+r)]^\delta$, $\delta \in \mathbb{R}$ & $\alpha > 0$ & \textup{(C)} \\[2pt]
$\lambda_0 + r^\beta[\log(e+r)]^s$, $\beta > 1$, $s \in \mathbb{R}$ & $\alpha > 0$ & \textup{(C)} \\
\hline
\end{tabular}
\end{center}
The first row follows from taking $D(r) = J(r)^\eta$ in \textup{(A)} and $\Psi(x) = x^\eta$ in \textup{(D)}; in particular, the borderline law $\lambda_0 + r/\log\log(e^e + r)$ is admissible when $\alpha = 1$. In the last three rows, for any prescribed $\alpha > 0$, one may choose $d' \in (0,1)$ sufficiently close to $1$ so that $2(1 - d')/(2 + d') \le \alpha$ and apply \textup{(C)}.

Further examples involving iterated logarithms are obtained by setting $L_1(r) = \log(e + r)$ and $L_{k+1}(r) = 1 + \log L_k(r)$ for $k \ge 1$, with the convention that an empty product equals $1$. For any integer $m \ge 2$, Definition~\ref{main:classes} admits
\begin{equation}\nonumber
\begin{aligned}
 \textup{(A)}\ (\alpha > 1):\quad & D(r) = \prod_{k=2}^{m} L_k(r)
  \qquad \Bigl(L_{m+1}'(r) = \frac{1}{(e + r)L_1(r)D(r)}\Bigr), \\
 \textup{(D)}\ (\alpha = 1):\quad & \Psi(x) = x \prod_{k=2}^{m} L_k(x), \\
 & \lambda(r) = \lambda_0 + \frac{r}{L_2(r)^\eta \prod_{k=3}^{m} L_k(r)^{s_k}}
  \qquad (0 \le \eta < 1,\ s_3, \ldots, s_m \in \mathbb{R}).
\end{aligned}
\end{equation}
\end{remark}

Before stating our main result, we fix the notation and conventions used throughout the paper. For a positive integer $s$ and $1 \le p \le \infty$, we denote the standard Lebesgue and Sobolev spaces on $\mathbb{R}^2$ by
\begin{equation}\nonumber
 L^p = L^p(\mathbb{R}^2),
 \qquad
 W^{s,p} = W^{s,p}(\mathbb{R}^2),
 \qquad
 H^s = W^{s,2},
 \qquad
 \|f\|_p = \|f\|_{L^p},
\end{equation}
and abbreviate spatial integration by
\begin{equation}\nonumber
 \int f\,dx = \int_{\mathbb{R}^2} f(x)\,dx.
\end{equation}
The material derivative is defined by
\begin{equation}\nonumber
 D_t f = \dot{f} = f_t + u \cdot \nabla f,
\end{equation}
while the effective viscous flux $G$ and vorticity $\omega$ are given by
\begin{equation}\nonumber
 G = h(\rho)\operatorname{div} u - P,
 \qquad
 \omega = \partial_1 u^2 - \partial_2 u^1.
\end{equation}
To measure the spatial decay of the density, we fix $\eta_0 > 0$ and introduce the weight
\begin{equation}\nonumber
 \bar{x} = (e + |x|^2)^{1/2} \log^{1 + \eta_0}(e + |x|^2).
\end{equation}

We now give the definition of strong solutions to \eqref{main:pde}.
\begin{definition}
A triple $(\rho, u, \theta)$ is called a strong solution to \eqref{main:pde} on $\mathbb{R}^2 \times (0,T)$ if all derivatives appearing in \eqref{main:pde} are regular distributions in $\mathbb{R}^2 \times (0,T)$ and the equations hold almost everywhere in $\mathbb{R}^2 \times (0,T)$.
\end{definition}

Our main result is stated as follows.

\begin{theorem}\label{main:global}
Assume one of the conditions \textup{(A)}--\textup{(D)} in Definition~\ref{main:classes}.
Let the initial data $(\rho_0, u_0, \theta_0)$ satisfy, for some $a > 1$, $q > 2$, and positive constants $\vartheta_*$ and $\theta_\infty$,
\begin{equation}\label{main:data}
\begin{cases}
\rho_0 \ge 0,\quad \bar{x}^a\rho_0 \in L^1 \cap H^1 \cap W^{1,q}(\mathbb{R}^2),\quad m_0 \triangleq \| \rho_0 \|_{L^1(\rr)}>0,\quad \nabla u_0 \in H^1(\mathbb{R}^2), \\
\sqrt{\rho_0} u_0, \sqrt{\rho_0} \theta_0 \in L^2(\mathbb{R}^2),\quad
\na \theta_0 \in L^2(\mathbb{R}^2),\quad \na K(\theta_0) \in L^2(\mathbb{R}^2), \quad \theta_0 \ge \vartheta_* > 0,
\end{cases}
\end{equation}
together with the momentum compatibility condition
\begin{equation}\label{main:compat}
 -\mu \Delta u_0
 - \nabla\big((\mu + \lambda(\rho_0))\operatorname{div} u_0\big)
 + \nabla(\rho_0 \theta_0)
 = \sqrt{\rho_0}\, g_1,
\end{equation}
for some $g_1 \in L^2(\mathbb{R}^2; \mathbb{R}^2)$.
Then the Cauchy problem \eqref{main:pde} admits a unique global strong solution $(\rho, u, \theta)$ in $\mathbb{R}^2 \times (0,\infty)$ satisfying, for any $0 < \tau < T < \infty$,
\begin{equation}\label{main:finite}
\begin{cases}
 \bar{x}^a \rho \in C([0,T]; L^1 \cap H^1 \cap W^{1,q}),\quad \rho_t \in L^\infty(0,T; L^q), \\
 u \in C([0,T]; H^1_{\rm loc}) \cap L^\infty(0,T; L^\infty),\quad \nabla u \in C([0,T]; L^2) \cap L^2(0,T; L^\infty), \\
 \nabla^2 u \in L^2(0,T; L^2 \cap L^q) \cap L^\infty(\tau,T; L^2), \\
 \theta \in L^\infty(\rr \times (\tau,T)), \quad \theta^{-1} \in L^\infty(\rr \times (0,T)), \\
 \nabla K(\theta) \in L^\infty(0,T; L^2) \cap L^2(0,T; H^1),\quad \nabla^2 K(\theta) \in L^\infty(\tau,T; L^2), \\
 \sqrt{\rho}\,\dot{u} \in L^\infty(0,T; L^2),\quad \nabla\dot{u} \in L^2(0,T; L^2),\quad u_t \in L^2(0,T; H^1_{\rm loc}), \\
 \sqrt{\rho\kappa(\theta)}\,\dot{\theta} \in L^2(0,T; L^2) \cap L^\infty(\tau,T; L^2),\quad \nabla D_t K(\theta) \in L^2(\tau,T; L^2).
\end{cases}
\end{equation}
In particular, $\lim_{|x|\to\infty}\rho(t,x) = 0$ for each $t \in [0,T]$, the initial velocity is attained locally in $H^1$, and the initial temperature is attained in the mass-weighted sense
\begin{equation}\label{main:trace}
 \int_{\mathbb{R}^2} \rho(t,x)|\theta(t,x) - \theta_0(x)|^2\,dx \le C_T t,
 \qquad 0 < t \le T.
\end{equation}
for a positive constant $C_T$ depending on $T$, the initial data, and the fixed coefficients. Moreover, for every $t \ge 0$, the mass, momentum, and total energy are conserved:
\begin{equation}\nonumber
 \int_{\mathbb{R}^2} \rho(t,x)\,dx = m_0,
 \qquad
 \int_{\mathbb{R}^2} (\rho u)(t,x)\,dx = \int_{\mathbb{R}^2} \rho_0 u_0\,dx,
 \qquad
 \int_{\mathbb{R}^2} (\rho E)(t,x)\,dx = E_0,
\end{equation}
where
\begin{equation}\nonumber
 E_0 = \int_{\mathbb{R}^2} \left(\frac{1}{2}\rho_0|u_0|^2 + c_v\rho_0\theta_0\right)dx.
\end{equation}
\end{theorem}

\begin{remark}
\normalfont
By \eqref{main:data} and the weighted entropy--Poincar\'e inequality in Lemma~\ref{spatial:entropy-poincare}, one has $\sqrt{\rho_0}\,u_0 \in L^2(\mathbb{R}^2)$, so both the initial momentum $\int_{\mathbb{R}^2}\rho_0 u_0\,dx$ and the initial total energy $E_0$ are finite.
Owing to the Galilean invariance of \eqref{main:pde}, the transformation
\begin{equation}\nonumber
 U = \frac{1}{m_0} \int_{\mathbb{R}^2} \rho_0 u_0\,dx,
 \qquad
 (\widetilde{\rho}, \widetilde{u}, \widetilde{\theta})(t,x)
 = (\rho, u - U, \theta)(t, x + tU),
\end{equation}
reduces the proof without loss of generality to the zero-momentum frame
\begin{equation}\nonumber
 \int_{\mathbb{R}^2} \rho_0 u_0\,dx = 0.
\end{equation}
\end{remark}

\section{Preliminaries}\label{sec:prelim}

In this section, we recall some known facts and elementary inequalities that will be used frequently later.

We begin with the local existence of strong solutions, which can be established by arguments similar to those in \cite{LXZthermal,LiangShuai,WX}.

\begin{lemma}\label{pre:local}
Under the assumptions on the coefficients and initial data in Theorem~\ref{main:global}, there exists a small time $T>0$ such that the problem \eqref{main:pde}--\eqref{intro:cauchy}, admits a unique strong solution $(\rho,u,\theta)$ on $\rr \times(0,T]$ satisfying \eqref{main:finite}.
\end{lemma}

Next, we record elementary growth and comparison estimates relating the bulk viscosity $\lambda$ to its primitive $\Lambda$ under \eqref{main:lambda-derivative} and Definition~\ref{main:classes}.

\begin{lemma}\label{pre:growth}
For classes \textup{(A)}, \textup{(B)}, \textup{(C)}, and \textup{(D)} in Definition~\ref{main:classes}, set $D$ as in \textup{(A)}, $D = D_h$, $D \equiv 1$, and $D = \Psi \circ J$, respectively, extending $D$ by a constant on $[0,2)$ when necessary so that $D : [0,\infty) \to [1,\infty)$ is nondecreasing. Then there exist constants $b_h \ge 1$, $0 \le \beta_D < 1/2$, and $c, C > 0$, depending only on the fixed coefficient laws, such that
\begin{equation}\label{vac:struct}
\begin{aligned}
 h(r) &\le C(1+r)^{b_h},
 \qquad \frac{r}{\mu + \lambda(r)} \le C D(r),
 \qquad D(r) \le C L(r)^{\beta_D}
 && (r \ge 0), \\
 \Lambda(r) &\ge c h(r),
 \qquad \Lambda(r) \ge \frac{c r}{D(r)}
 && (r \ge 2).
\end{aligned}
\end{equation}
\end{lemma}

\begin{proof}
Since $h(r) = 2\mu + \lambda(r) \ge 2\mu > 0$, assumption \eqref{main:lambda-derivative} implies
\begin{equation}\nonumber
 \left|\frac{d}{dr}\log h(r)\right|
 = \frac{|\lambda'(r)|}{h(r)}
 \le \frac{C_\lambda}{r}
 \qquad (r > 0).
\end{equation}
For any $r \ge 1$, integrating this bound over $[1,r]$ yields
\begin{equation}\nonumber
 \log\frac{h(r)}{h(1)}
 = \int_1^r \frac{h'(s)}{h(s)}\,ds
 \le C_\lambda \int_1^r \frac{ds}{s}
 = C_\lambda \log r,
\end{equation}
and hence $h(r) \le h(1)\,r^{C_\lambda}$. Taking $b_h = \max\{1, C_\lambda\}$ and noting that $h(r) \le \sup_{0 \le s \le 1} h(s)$ on $[0,1]$, we obtain
\begin{equation}\nonumber
 h(r) \le C(1+r)^{b_h} \qquad (r \ge 0).
\end{equation}

We next verify the bounds on $D$ and $r/(\mu + \lambda(r))$. In cases \textup{(A)} and \textup{(C)}, the bound $D(r) \le C L(r)^{\beta_D}$ follows directly from Definition~\ref{main:classes} (with $\beta_D = 0$ in \textup{(C)}). In case \textup{(B)}, the limit $\lim_{r\to\infty}\frac{\log D_h(r)}{\log L(r)} = 0$ implies $D_h(r) \le C_\varepsilon L(r)^\varepsilon$ for any fixed $\varepsilon \in (0, 1/2)$. In case \textup{(D)}, integrating $\Psi'(s)/\Psi(s) \le p/s$ over $[1,x]$ gives
\begin{equation}\nonumber
 \Psi(x) \le \Psi(1)\,x^p \qquad (x \ge 1),
\end{equation}
so that $D(r) = \Psi(J(r)) \le C J(r)^p \le C L(r)^{\beta_D}$ for any fixed $\beta_D \in (0, 1/2)$. Moreover, since $h(r) \le 2(\mu + \lambda(r))$, combining the large-density lower bounds in Definition~\ref{main:classes} with the elementary estimate $\frac{r}{\mu + \lambda(r)} \le \frac{r}{\mu}$ on bounded intervals and $D \ge 1$ yields
\begin{equation}\nonumber
 \frac{r}{\mu + \lambda(r)} \le C D(r) \qquad (r \ge 0).
\end{equation}

Finally, for any $r \ge 2$ and $s \in [r/2, r]$, integrating $h'(z)/h(z) \le C_\lambda / z$ over $[s,r]$ gives
\begin{equation}\nonumber
 \log\frac{h(r)}{h(s)}
 \le C_\lambda \log\frac{r}{s}
 \le C_\lambda \log 2,
 \qquad \text{and hence} \qquad
 h(s) \ge 2^{-C_\lambda} h(r).
\end{equation}
Since $[r/2, r] \subset [1,r]$ and $h > 0$, it follows that
\begin{equation}\nonumber
 \Lambda(r)
 = \int_1^r \frac{h(s)}{s}\,ds
 \ge \int_{r/2}^r \frac{h(s)}{s}\,ds
 \ge 2^{-C_\lambda} h(r) \int_{r/2}^r \frac{ds}{s}
 = (2^{-C_\lambda}\log 2)\,h(r).
\end{equation}
On the other hand, using $h(s) \ge \mu + \lambda(s) \ge c_0 s / D(s)$ and the monotonicity of $D$, we deduce that
\begin{equation}\nonumber
 \Lambda(r)
 \ge \int_{r/2}^r \frac{h(s)}{s}\,ds
 \ge c_0 \int_{r/2}^r \frac{ds}{D(s)}
 \ge \frac{c_0}{D(r)} \int_{r/2}^r ds
 = \frac{c_0 r}{2 D(r)}
 \qquad (r \ge 2).
\end{equation}
This completes the proof of Lemma~\ref{pre:growth}.
\end{proof}

We next recall the two-dimensional Gagliardo--Nirenberg inequality \cite{NI,TG}.
\begin{lemma}\la{gn1}
Assume that $\OM=\rr$.
For $2<p<\infty$, there exists a positive constant $C$ such that for any $f \in H^1(\rr)$,
\be\ba\la{gn11}
\| f \|_{L^p} \le Cp^{1/2}\| f \|^{2/p}_{L^2} \| \na f \|^{1-2/p}_{L^2}.
\ea\ee
\end{lemma}

Motivated by the weighted Poincar\'e inequality in
\cite{FJLL1}, we establish the following estimates under
the weaker assumption of $L\log L$ integrability of the density.

\begin{lemma}\label{spatial:entropy-poincare}
Let $\rho \ge 0$ have mass $m_0 = \int \rho\,dx \in (0,\infty)$, and assume that for some
$x_0\in\mathbb R^2$,
\[
 \mathcal A_\rho(x_0)=1+\int\rho\log(e+\rho)
                         +\int\rho\log(e+|x-x_0|) < \infty.
\]
For $v\in H^1_{\rm loc}$ with $\nabla v\in L^2$ one has
$v\in L^2(\rho\,dx)$ and
\begin{equation}\label{spatial:entropy-poincare-est}
 \int\rho|v-v_\rho|^2\le C\mathcal A_\rho(x_0)\|\nabla v\|_2^2,
 \qquad v_\rho=m_0^{-1}\int\rho v,
\end{equation}
where $C>0$ is independent of $\rho$, $v$ and $x_0$.

Moreover, for any $c > 2$, set
\begin{equation}\nonumber
 Q_0 = \left(-\frac{1}{2},\frac{1}{2}\right)^2,
 \qquad
 v_0 = \int_{x_0 + Q_0} v\,dx,
 \qquad
 w_{x_0}(x) = (1 + |x-x_0|^2)^{-c/2}.
\end{equation}
Then, for any $2\le p<\infty$, there exists a constant $C_c > 0$ depending only on $c$ such that
\begin{equation}\label{spatial:cube-sob}
 \|v-v_0\|_{L^p(w_{x_0}\,dx)}\le C_c\sqrt p\,\|\nabla v\|_2.
\end{equation}
\end{lemma}

\begin{proof}
By translation, we may assume without loss of generality that $x_0 = 0$. Set $E = \|\nabla v\|_2^2$. We may assume $E > 0$, since otherwise $v$ is constant and both \eqref{spatial:entropy-poincare-est} and \eqref{spatial:cube-sob} hold trivially.

We first derive a weighted Trudinger inequality. Decompose $\mathbb{R}^2$ into the dyadic regions
\begin{equation}\nonumber
 A_0 = B_2,
 \qquad
 A_j = B_{2^{j+1}} \setminus B_{2^j} \quad (j \ge 1),
\end{equation}
and denote the corresponding averages of $v$ by
\begin{equation}\nonumber
 a_j = \frac{1}{|A_j|} \int_{A_j} v\,dx \qquad (j \ge 0).
\end{equation}
Since $Q_0 \subset B_2$ and $|Q_0| = 1$, H\"older's inequality and the Poincar\'e inequality on $B_2$ yield
\begin{equation}\nonumber
 |v_0 - a_0|
 = \left|\int_{Q_0}(v - a_0)\,dx\right|
 \le \|v - a_0\|_{L^2(Q_0)}
 \le \|v - a_0\|_{L^2(B_2)}
 \le C\|\nabla v\|_{L^2(B_2)}.
\end{equation}
Next, for $j \ge 0$, introduce the overlapping annuli
\begin{equation}\nonumber
 D_0 = B_4,
 \qquad
 D_j = B_{2^{j+2}} \setminus \overline{B}_{2^j} \quad (j \ge 1),
 \qquad
 b_j = \frac{1}{|D_j|}\int_{D_j} v\,dx.
\end{equation}
Since $A_j, A_{j+1} \subset D_j$ up to null sets and $|A_j|, |A_{j+1}|, |D_j| \asymp 2^{2j}$, the Poincar\'e inequality and scaling give
\begin{equation}\nonumber
\begin{aligned}
 |a_{j+1} - a_j|
 &\le |a_{j+1} - b_j| + |a_j - b_j| \\
 &\le C 2^{-j}\|v - b_j\|_{L^2(D_j)}
  \le C\|\nabla v\|_{L^2(D_j)}.
\end{aligned}
\end{equation}
Because the family $\{D_k\}_{k \ge 0}$ has finite overlap, summing over $0 \le k \le j-1$ and applying the Cauchy--Schwarz inequality yield
\begin{equation}\label{spatial:annular-means}
 |a_j - v_0|
 \le |a_0 - v_0| + \sum_{k=0}^{j-1} |a_{k+1} - a_k|
 \le C\sqrt{j+1}\,E^{1/2}
 \qquad (j \ge 0).
\end{equation}

By the Trudinger inequality on the fixed domains $B_2$ and $B_2 \setminus \overline{B}_1$, together with scaling invariance of the planar Dirichlet integral, there exist constants $\varepsilon_0, C > 0$, independent of $j$, such that
\begin{equation}\nonumber
 \int_{A_j} \exp\left(\frac{\varepsilon_0|v - a_j|^2}{E}\right)dx
 \le C|A_j|
 \qquad (j \ge 0).
\end{equation}
Using $|v - v_0|^2 \le 2|v - a_j|^2 + 2|a_j - v_0|^2$, \eqref{spatial:annular-means}, $|A_j| \le C 2^{2j}$, and $w_0 \le C_c 2^{-cj}$ on $A_j$, we obtain for every $0 < 2\varepsilon \le \varepsilon_0$,
\begin{equation}\nonumber
\begin{aligned}
 \int_{A_j} w_0 \exp\left(\frac{\varepsilon|v - v_0|^2}{E}\right)dx
 &\le C_c 2^{-cj} e^{C\varepsilon(j+1)}
      \int_{A_j} \exp\left(\frac{2\varepsilon|v - a_j|^2}{E}\right)dx \\
 &\le C_c 2^{-(c-2)j} e^{C\varepsilon(j+1)}.
\end{aligned}
\end{equation}
Since $c > 2$, choosing $\varepsilon = \varepsilon_c \in (0, \varepsilon_0/2]$ small enough that $e^{C\varepsilon_c} < 2^{c-2}$ and summing over $j \ge 0$ lead to
\begin{equation}\label{spatial:weighted-trudinger}
 \int_{\mathbb{R}^2} w_0
 \exp\left(\frac{\varepsilon_c|v - v_0|^2}{\|\nabla v\|_2^2}\right)dx
 \le C_c.
\end{equation}

To prove \eqref{spatial:entropy-poincare-est}, fix $c = 4$ and apply Young's inequality $ab \le a\log(e+a) + e^b$ with $a = \rho/w_0$ and $b = \varepsilon_4|v - v_0|^2/E$ in the measure $w_0\,dx$. Since $w_0^{-1} = (1+|x|^2)^2 \le (e+|x|)^4$, we have
\begin{equation}\nonumber
 \log\left(e + \frac{\rho}{w_0}\right)
 \le \log(e + \rho) + 4\log(e + |x|).
\end{equation}
Combining this bound with \eqref{spatial:weighted-trudinger} yields
\begin{equation}\nonumber
\begin{aligned}
 \int \rho|v - v_0|^2\,dx
 &\le C E \left(1 + \int \rho\log\left(e + \frac{\rho}{w_0}\right)dx\right) \\
 &\le C\left(1 + \int \rho\log(e + \rho)\,dx
               + \int \rho\log(e + |x|)\,dx\right)\|\nabla v\|_2^2
  = C\mathcal{A}_\rho(0)\|\nabla v\|_2^2.
\end{aligned}
\end{equation}
Since $\rho$ has finite mass, this implies $v \in L^2(\rho\,dx)$. Minimizing the weighted $L^2$-distance over constants gives
\begin{equation}\nonumber
 \int \rho|v - v_\rho|^2\,dx \le \int \rho|v - v_0|^2\,dx,
\end{equation}
which proves \eqref{spatial:entropy-poincare-est}.

Finally, for any $c > 2$ and $p \ge 2$, using the elementary inequality $s^{p/2} \le \bigl(p/(2e)\bigr)^{p/2} e^s$ ($s \ge 0$) with $s = \varepsilon_c|v - v_0|^2/E$, we deduce from \eqref{spatial:weighted-trudinger} that
\begin{equation}\nonumber
 \int w_0|v - v_0|^p\,dx
 \le \left(\frac{p E}{2e\varepsilon_c}\right)^{p/2}
     \int w_0 \exp\left(\frac{\varepsilon_c|v - v_0|^2}{E}\right)dx
 \le C_c \left(\frac{p E}{2e\varepsilon_c}\right)^{p/2}.
\end{equation}
Taking $p$-th roots establishes \eqref{spatial:cube-sob} and completes the proof of Lemma~\ref{spatial:entropy-poincare}.
\end{proof}

\section{A priori estimates on the whole plane}\label{sec:apriori}

In this section, we assume that the initial data $(\rho_0,u_0,\theta_0)$ satisfy (\ref{main:data}).
Let $(\rho,u,\theta)$ be the strong solution to (\ref{main:pde})--(\ref{intro:cauchy}) on $\mathbb{R}^2 \times (0,T]$, provided by Lemma \ref{pre:local}.

Define
\[
Q = 2\mu |D(u)|^2+\lambda(\div u)^2,\quad
 W(x)=(1+|x|^2)^{1/2},\quad w=W^{-c},\quad
 \|f\|_{p,w}=\left(\int |f|^p w\,dx\right)^{1/p}.
\]

\subsection{Low-order estimates}\label{vac:low}

This subsection is devoted to deriving the upper bound for the density.

\begin{lemma}\label{vac:low-anchor}
There exists a positive constant $C$ such that
\begin{equation}\label{vac:entropy}
 \sup_{t\le T}\int\rho\bigl(W^{a_0}+L(\rho)+|u|^2+\theta\bigr)
 +\int_0^T\!\int\left(\frac Q\theta+
         \frac{\kappa(\theta)|\nabla\theta|^2}{\theta^2}\right)
 \le C.
\end{equation}
Moreover, there exist constants $N, L_0 > 1$ and $c > 0$, independent of $t$, such that the set
\begin{equation}\nonumber
 \Gamma_t = B_N \cap \{L_0^{-1} \le \rho \le L_0,\ \theta \le L_0\}
\end{equation}
satisfies $|\Gamma_t| \ge c$.
In addition, for any $f \in H^1_{\rm loc}(\mathbb{R}^2)$ with $\nabla f \in L^2(\mathbb{R}^2)$ and any $p \ge 2$,
\begin{equation}\label{vac:sob}
 \|f\|_{L^p(w\,dx)}
 \le C_{c}\sqrt p\,
 \bigl(\|\nabla f\|_2+\|f\|_{L^2(\Gamma_t)}\bigr).
\end{equation}
If further $z \in H^1_{\rm loc}(\mathbb{R}^2)$ satisfies $\nabla z \in L^2(\mathbb{R}^2)$, $z = 0$ a.e.\ on $\Gamma_t$, and $\|\nabla z\|_2 > 0$, then
\begin{equation}\label{vac:tm}
 \int w\exp\left(
 \frac{cz^2}{\|\nabla z\|_2^2}
 \right)dx\le C.
\end{equation}
In particular,
\begin{equation}\label{vac:tm-u}
 \int w\exp\left(
 \frac{c|u|^2}{(1+\|\nabla u\|_2)^2}
 \right)dx\le C.
\end{equation}
\end{lemma}
\begin{proof}
We first bound the spatial moment $X(t) = \int \rho W^{a_0}\,dx$. Using the continuity equation, $|\nabla W| \le 1$, energy conservation, and H\"older's inequality with $1 < a_0 < 2$, we obtain
\begin{equation}\nonumber
 \left|\frac{d}{dt}X(t)\right|
 \le a_0 \int \rho|u|W^{a_0-1}\,dx
 \le a_0(2E_0)^{1/2} \left(\int \rho W^{2a_0-2}\,dx\right)^{1/2}
 \le C X(t)^{1 - 1/a_0},
\end{equation}
which yields
\begin{equation}\label{vac:spatial-moment}
 \sup_{0 \le t \le T} \int \rho W^{a_0}\,dx \le C.
\end{equation}

Next, multiplying the temperature equation by $\theta^{-1}$ and using the continuity equation gives the entropy identity
\begin{equation}\label{vac:entropy-identity}
 \frac{d}{dt} \int \rho(c_v\log\theta - \log\rho)\,dx
 = \int \left(\frac{Q}{\theta}
        + \frac{\kappa(\theta)|\nabla\theta|^2}{\theta^2}\right)dx.
\end{equation}
Since $r\log(1/r) \le r W^{a_0} + e^{-1} e^{-W^{a_0}}$ for $0 < r \le 1$, it follows from \eqref{vac:spatial-moment} and energy conservation that
\begin{equation}\nonumber
 \sup_{0 \le t \le T} \int \rho(\log\rho)_-\,dx \le C,
 \qquad
 \int \rho(\log\theta)_+\,dx \le \int \rho\theta\,dx \le \frac{E_0}{c_v}.
\end{equation}
Integrating \eqref{vac:entropy-identity} over $(0,t)$ and separating the positive and negative parts of the logarithms, we deduce
\begin{equation}\nonumber
\begin{aligned}
 & c_v \int \rho(\log\theta)_-\,dx
   + \int \rho(\log\rho)_+\,dx
   + \int_0^t\!\int \left(\frac{Q}{\theta}
         + \frac{\kappa(\theta)|\nabla\theta|^2}{\theta^2}\right)dx\,ds \\
 & = c_v \int \rho(\log\theta)_+\,dx
          + \int \rho(\log\rho)_-\,dx
          - \int \rho_0(c_v\log\theta_0 - \log\rho_0)\,dx
   \le C,
\end{aligned}
\end{equation}
where the initial entropy is finite by \eqref{main:data}. Combining this bound with \eqref{vac:spatial-moment}, conservation of mass and energy, and $L(r) \le C + (\log r)_+$ establishes \eqref{vac:entropy}.

To bound $|\Gamma_t|$ from below, we use \eqref{vac:entropy} and Chebyshev's inequality to estimate
\begin{equation}\nonumber
\begin{aligned}
 \int_{\mathbb{R}^2 \setminus \Gamma_t} \rho\,dx
 &\le \int_{\mathbb{R}^2 \setminus B_N} \rho\,dx
      + \int_{B_N \cap \{\rho < L_0^{-1}\}} \rho\,dx
      + \int_{\{\rho > L_0\}} \rho\,dx
      + \int_{\{\theta > L_0\}} \rho\,dx \\
 &\le C N^{-a_0}
      + \frac{|B_N|}{L_0}
      + \frac{C}{\log L_0}
      + \frac{E_0}{c_v L_0}.
\end{aligned}
\end{equation}
Choosing $N > 1$ (with $Q_0 = (-1/2,1/2)^2 \subset B_N$) and then $L_0 > 1$ sufficiently large ensures that the right-hand side is at most $m_0/2$, whence
\begin{equation}\nonumber
 \frac{m_0}{2} \le \int_{\Gamma_t} \rho\,dx \le L_0 |\Gamma_t|,
 \qquad \text{and hence} \qquad
 |\Gamma_t| \ge \frac{m_0}{2L_0} > 0.
\end{equation}

We now prove \eqref{vac:sob}. Set $f_0 = \int_{Q_0} f\,dx$. Since $Q_0 \subset B_N$, the Poincar\'e inequality on $B_N$ yields $\|f - f_0\|_{L^2(B_N)} \le C_N \|\nabla f\|_{L^2(B_N)}$. As $\Gamma_t \subset B_N$, we have
\begin{equation}\nonumber
 |\Gamma_t|^{1/2}|f_0|
 \le \|f - f_0\|_{L^2(\Gamma_t)} + \|f\|_{L^2(\Gamma_t)}
 \le C_N \|\nabla f\|_2 + \|f\|_{L^2(\Gamma_t)},
\end{equation}
which together with $|\Gamma_t| \ge m_0/(2L_0)$ gives
\begin{equation}\label{vac:anchor-mean}
 |f_0| \le C \bigl(\|\nabla f\|_2 + \|f\|_{L^2(\Gamma_t)}\bigr).
\end{equation}
Applying \eqref{spatial:cube-sob} with $x_0 = 0$, using \eqref{vac:anchor-mean}, and noting that $\sup_{p \ge 2}\|1\|_{L^p(w\,dx)} < \infty$, we arrive at
\begin{equation}\nonumber
 \|f\|_{L^p(w\,dx)}
 \le \|f - f_0\|_{L^p(w\,dx)} + |f_0|\|1\|_{L^p(w\,dx)}
 \le C_{c}\sqrt{p}\,\bigl(\|\nabla f\|_2 + \|f\|_{L^2(\Gamma_t)}\bigr),
\end{equation}
which proves \eqref{vac:sob}.

Finally, combining \eqref{spatial:weighted-trudinger} with \eqref{vac:anchor-mean} and $|f|^2 \le 2|f - f_0|^2 + 2|f_0|^2$, we obtain for $\varepsilon = \varepsilon_c/2$,
\begin{equation}\nonumber
 \int w \exp\left(
   \frac{\varepsilon |f|^2}{(\|\nabla f\|_2 + \|f\|_{L^2(\Gamma_t)})^2}
 \right)dx \le C_{c},
\end{equation}
whenever the denominator is nonzero. Taking $f = z$ with $z = 0$ a.e.\ on $\Gamma_t$ gives \eqref{vac:tm}. Moreover, since $\rho \ge L_0^{-1}$ on $\Gamma_t$, we have
\begin{equation}\nonumber
 \|u\|_{L^2(\Gamma_t)}
 \le L_0^{1/2}\|\sqrt{\rho}\,u\|_2
 \le (2L_0 E_0)^{1/2}.
\end{equation}
Applying the preceding exponential estimate to $|u|$ and using $|\nabla|u|| \le |\nabla u|$ yields \eqref{vac:tm-u}, after decreasing $c$ if necessary. This completes the proof of Lemma~\ref{vac:low-anchor}.
\end{proof}

\begin{lemma}\label{vac:moment}
Set
\begin{equation}\nonumber
 U(t) = 1 + \int_0^t\!\int Q\,dx\,ds,
 \qquad
 N_p(t) = 2 + \sup_{0 \le s \le t} \int \rho(s)^p\,dx.
\end{equation}
Then, for every fixed $p > 2$, there exist positive constants $C_{p}$ and $c_p$ such that
\begin{equation}\label{vac:moment-bound}
 \log N_p(t) \le C_{p}\bigl(1 + U(t)\bigr)^{c_p},
 \qquad 0 \le t \le T.
\end{equation}
\end{lemma}

\begin{proof}
Define
\begin{equation}\nonumber
 \psi = (-\Delta)^{-1}\operatorname{div}(\rho u),
 \qquad
 R_i = \partial_i(-\Delta)^{-1/2}.
\end{equation}
Applying $(-\Delta)^{-1}\operatorname{div}$ to the momentum equation and combining the result with $D_t\Lambda(\rho) = -h(\rho)\operatorname{div} u$, we obtain
\begin{equation}\nonumber
 \psi_t = -\sum_{i,j=1}^2 R_i R_j(\rho u_i u_j) - h(\rho)\operatorname{div} u + P,
\end{equation}
and hence
\begin{equation}\nonumber
 D_t(\Lambda - \psi) + P = -F,
 \qquad
 F = \sum_{i,j=1}^2 [u_i, R_i R_j](\rho u_j).
\end{equation}

Choose $m \in (1/2 + \beta_D, 1)$, $\xi \in (m,1)$, and $k_* > p + 2$ sufficiently large such that
\begin{equation}\nonumber
 m + \frac{m + \beta_D}{k_*} < \xi < 1,
 \qquad
 k > m(k+1) \quad (k \in [k_*, k_*+1]).
\end{equation}
All constants below are uniform for $k \in [k_*, k_*+1]$. Fix $K_0 \ge 2$ large enough so that \eqref{vac:struct} holds for $\rho > K_0$, and put
\begin{equation}\nonumber
\begin{aligned}
 & N = 2 + \Lambda(K_0),
   \qquad \phi = (\Lambda - \psi - N)_+, \\
 & I_k = \int \rho\phi^k\,dx,
   \qquad M_k = 1 + I_k,
   \qquad J_k = 1 + \int_{\{\rho > K_0\}} \rho\Lambda^k\,dx,
\end{aligned}
\end{equation}
with $\rho\phi^k$ understood as zero on vacuum.

By the Hardy--Littlewood--Sobolev inequality, H\"older's inequality, and energy conservation,
\begin{equation}\nonumber
 \|\psi\|_{2k}
 \le C\|\rho u\|_{\frac{2k}{k+1}}
 \le C\|\sqrt{\rho}\,u\|_2\|\rho\|_k^{1/2}
 \le C\|\rho\|_k^{1/2}.
\end{equation}
Moreover, \eqref{vac:struct} and conservation of mass yield
\begin{equation}\nonumber
 \|\rho\|_k^k + \int_{\{\rho > K_0\}} \frac{\rho^{k+1}}{D(\rho)^k}\,dx
 \le C J_k.
\end{equation}
Since $\Lambda \le N + \phi + |\psi|$ on $\{\rho > K_0\}$, interpolating $\|\rho\|_2 \le \|\rho\|_1^{\frac{k-2}{2(k-1)}}\|\rho\|_k^{\frac{k}{2(k-1)}}$ gives
\begin{equation}\nonumber
\begin{aligned}
 J_k
 &\le C M_k + C \int \rho|\psi|^k\,dx
  \le C M_k + C\|\rho\|_2\|\psi\|_{2k}^k \\
 &\le C M_k + C\|\rho\|_k^{\frac{k}{2(k-1)} + \frac{k}{2}}
  \le C M_k + C J_k^{\frac{1}{2} + \frac{1}{2(k-1)}}.
\end{aligned}
\end{equation}
As $k > 2$, Young's inequality implies
\begin{equation}\label{vac:primmom}
 J_k + \|\rho\|_k^k + \|\psi\|_{2k}^{2k} \le C M_k.
\end{equation}

To estimate the commutator $F$, introduce the convex Orlicz function
\begin{equation}\nonumber
 B_k(s) = \int_0^s \frac{r^k}{L(r)^{m(k+1)}}\,dr,
 \qquad s \ge 0.
\end{equation}
Since $k > m(k+1)$, $B_k'$ is increasing and $B_k$ together with its Legendre conjugate $B_k^*$ satisfy
\begin{equation}\nonumber
 B_k(s) \asymp \frac{s^{k+1}}{L(s)^{m(k+1)}},
 \qquad
 B_k^*(z) \asymp z^{1 + 1/k} L(z)^{m + m/k}.
\end{equation}
For $T_u g = \sum_{i,j=1}^2 [u_i, R_i R_j] g_j$, the commutator estimate in \cite{CRW} and $\|u\|_{\rm BMO} \le C\|\nabla u\|_2$ yield
\begin{equation}\nonumber
 \|T_u g\|_s \le C_s \|\nabla u\|_2 \|g\|_s,
 \qquad 1 < s < \infty.
\end{equation}
Fix exponents $2 < q_- < k_* + 1 < k_* + 2 < q_+$. If $\|\nabla u\|_2 = 0$, then $T_u g = 0$; otherwise, splitting $g = g\mathbf{1}_{\{|g| > r\}} + g\mathbf{1}_{\{|g| \le r\}}$ and applying Chebyshev's inequality at $q_\pm$ give
\begin{equation}\nonumber
 \bigl|\{|T_u g| > 2r\|\nabla u\|_2\}\bigr|
 \le C r^{-q_-}\int_{\{|g| > r\}} |g|^{q_-}\,dx
   + C r^{-q_+}\int_{\{|g| \le r\}} |g|^{q_+}\,dx.
\end{equation}
Multiplying by $B_k'(r)$, integrating over $(0,\infty)$, and applying Fubini's theorem, we deduce
\begin{equation}\nonumber
\begin{aligned}
 \int B_k\left(\frac{|T_u g|}{2\|\nabla u\|_2}\right)dx
 &\le C \int |g|^{q_-}\left(\int_0^{|g|} r^{-q_-} B_k'(r)\,dr\right)dx
    + C \int |g|^{q_+}\left(\int_{|g|}^\infty r^{-q_+} B_k'(r)\,dr\right)dx \\
 &\le C \int B_k(|g|)\,dx,
\end{aligned}
\end{equation}
where we used $q_- < k + 1 < q_+$ to bound the inner integrals by $C|g|^{-q_\pm} B_k(|g|)$.

Set
\begin{equation}\nonumber
 V = 1 + \|\nabla u\|_2,
 \qquad v = \frac{|u|}{V},
 \qquad \gamma = \frac{c}{c + a_0},
\end{equation}
so that \eqref{vac:entropy} and \eqref{vac:tm} give
\begin{equation}\nonumber
 \int \rho W^{a_0}\,dx + \int w\exp(c v^2)\,dx \le C.
\end{equation}
For a large constant $A > 1$, on $\{v^2 \le A L(\rho)\}$ the monotonicity of $B_k$, $m > 1/2 + \beta_D$, and $D(\rho) \le C L(\rho)^{\beta_D}$ yield
\begin{equation}\nonumber
 B_k\left(\frac{\rho v}{\sqrt{A}}\right)
 \le B_k\bigl(\rho L(\rho)^{1/2}\bigr)
 \le C \rho^{k+1} L(\rho)^{-(m - 1/2)(k+1)}
 \le C\rho + C\mathbf{1}_{\{\rho > K_0\}} \frac{\rho^{k+1}}{D(\rho)^k}.
\end{equation}
On $\{v^2 > A L(\rho)\}$, choosing $A$ sufficiently large ensures $\rho^{k+1-\gamma} v^{k+1} \le C\exp\bigl((1-\gamma)c v^2\bigr)$, and hence by $a_0\gamma = c(1-\gamma)$,
\begin{equation}\nonumber
 B_k\left(\frac{\rho v}{\sqrt{A}}\right)
 \le C \rho^{k+1} v^{k+1}
 \le C (\rho W^{a_0})^\gamma \bigl(w\exp(c v^2)\bigr)^{1-\gamma}.
\end{equation}
Integrating over both regions via \eqref{vac:primmom} and H\"older's inequality, and applying the modular commutator bound to $g = \rho u / (\sqrt{A}\,V)$, we obtain
\begin{equation}\nonumber
 \int B_k\left(\frac{|F|}{C V^2}\right)dx
 \le C \int B_k\left(\frac{|\rho u|}{\sqrt{A}\,V}\right)dx
 \le C M_k.
\end{equation}

Next, we estimate $B_k^*(\rho\phi^{k-1})$. For a sufficiently large constant $C_1 > 1$, introduce
\begin{equation}\nonumber
 \mathcal{A} = \left\{\frac{\rho}{L(\rho)^{\beta_D}} \le C_1(1 + \phi)\right\}.
\end{equation}
On $\mathcal{A}$, we have $L(\rho) \le C L(1 + \phi)$ and $\rho \le C(1 + \phi)L(1 + \phi)^{\beta_D}$. In particular, $\rho$ is bounded on $\mathcal{A} \cap \{\phi < 1\}$, so $B_k^*(\rho\phi^{k-1}) \le C\rho$ there, whereas on $\mathcal{A} \cap \{\phi \ge 1\}$,
\begin{equation}\nonumber
\begin{aligned}
 B_k^*(\rho\phi^{k-1})
 &\le C (\rho\phi^{k-1})^{1 + 1/k} L(\rho\phi^{k-1})^{m + m/k} \\
 &\le C \rho\phi^k \left(\frac{\rho}{\phi}\right)^{1/k} L(\phi)^{m + m/k}
  \le C \rho\phi^k L(\phi)^{m + (m + \beta_D)/k}
  \le C \rho\phi^k L(\phi)^\xi.
\end{aligned}
\end{equation}
On $\mathcal{A}^c \cap \{\phi > 0\}$, choosing $C_1$ sufficiently large and using \eqref{vac:struct} give $\phi < \rho$ and
\begin{equation}\nonumber
 \psi = \Lambda - N - \phi
 \ge c\frac{\rho}{L(\rho)^{\beta_D}} - N - \phi
 \ge c_1\frac{\rho}{L(\rho)^{\beta_D}},
\end{equation}
so that
\begin{equation}\nonumber
 B_k^*(\rho\phi^{k-1})
 \le C \rho^{k+1} L(\rho)^{m + m/k}
 \le C|\psi|^{2k}.
\end{equation}
Combining these bounds with \eqref{vac:primmom} yields
\begin{equation}\nonumber
 \int B_k^*(\rho\phi^{k-1})\,dx
 \le C M_k + C \int_{\{\phi \ge 1\}} \rho\phi^k L(\phi)^\xi\,dx.
\end{equation}

Using the continuity equation, $P \ge 0$, and Young's inequality for $(B_k, B_k^*)$, we obtain
\begin{equation}\nonumber
\begin{aligned}
 \frac{d}{dt} I_k
 &= -k \int \rho\phi^{k-1}(P + F)\,dx
  \le k \int \rho\phi^{k-1}|F|\,dx \\
 &\le C V^2 \int \left[B_k\left(\frac{|F|}{C V^2}\right) + B_k^*(\rho\phi^{k-1})\right]dx
  \le C b(t)\left(M_k + \int_{\{\phi \ge 1\}} \rho\phi^k L(\phi)^\xi\,dx\right),
\end{aligned}
\end{equation}
where $b(t) = 1 + \|\nabla u(t)\|_2^2$. Since $0 < \xi < 1$, applying $L(r)^\xi \le C_\xi(1 + \delta^{-n}) + \delta\log r$ ($r \ge 1$, $0 < \delta \le 1$) with $n = \xi/(1-\xi)$ gives
\begin{equation}\nonumber
 \frac{d}{dt} I_k
 \le C(1 + \delta^{-n}) b(t) M_k
   + \delta b(t) \int_{\{\phi \ge 1\}} \rho\phi^k \log\phi\,dx.
\end{equation}

For a fixed $t \in [0,T]$, introduce the decreasing exponent and parameters
\begin{equation}\nonumber
 H = 1 + \int_0^t b(s)\,ds,
 \qquad
 k(s) = k_* + 1 - \frac{1}{H}\int_0^s b(\tau)\,d\tau \in [k_*, k_*+1],
 \qquad
 \delta = \frac{1}{2H}.
\end{equation}
Differentiating $M_{k(s)}$ in $s$ and noting that $\int_{\{0 < \phi < 1\}} \rho\phi^{k(s)}|\log\phi|\,dx \le C m_0$, we find
\begin{equation}\nonumber
\begin{aligned}
 \frac{d}{ds} M_{k(s)}
 &\le C(1 + \delta^{-n}) b(s) M_{k(s)}
    + \delta b(s)\int_{\{\phi \ge 1\}} \rho\phi^{k(s)}\log\phi\,dx
    - \frac{b(s)}{H}\int \rho\phi^{k(s)}\log\phi\,dx \\
 &\le C(1 + H^n) b(s) M_{k(s)}.
\end{aligned}
\end{equation}
Integrating over $(0,s)$ gives
\begin{equation}\label{ll21}
 \log M_{k(s)}(s) \le \log M_{k(0)}(0) + C H^{1+n},
 \qquad 0 \le s \le t,
\end{equation}
where $M_{k(0)}(0)$ is bounded in terms of the initial data since $\rho_0 \in L^1 \cap L^\infty$ and $\|\psi(0)\|_{2k(0)} \le C\|\sqrt{\rho_0}\,u_0\|_2\|\rho_0\|_{k(0)}^{1/2}$.

Finally, because $k(s) \ge k_* > p$, conservation of mass and \eqref{vac:primmom} yield
\begin{equation}\label{ll22}
 \int \rho(s)^p\,dx
 \le C\left(1 + \int \rho(s)^{k(s)}\,dx\right)
 \le C M_{k(s)}(s),
\end{equation}
whence $\log N_p(t) \le C_{p} H^{1+n}$. By Korn's inequality,
\begin{equation}\label{ll23}
 H = 1 + \int_0^t \bigl(1 + \|\nabla u\|_2^2\bigr)\,ds
 \le 1 + T + C \int_0^t\!\int Q\,dx\,ds
 \le C \bigl(1 + U(t)\bigr).
\end{equation}
Choosing $c_p = 1 + n$ in (\ref{ll21}) and using (\ref{ll22}), and (\ref{ll23}), we obtain \eqref{vac:moment-bound} and complete the proof of Lemma \ref{vac:moment}.
\end{proof}

\begin{lemma}\label{vac:heat-band}
Let $0 \le \nu < 1$ satisfy $a_\nu = \alpha + \nu - 1 > 0$, let $l \ge 2L_0 + 2$, and let $I \subset [0,T]$ be a closed subinterval. Set
\begin{equation}\nonumber
 f_l = \mathbf{1}_{\{\theta > l\}}\frac{\rho^2\theta^\nu}{\mu + \lambda(\rho)},
 \qquad
 F_l = 1 + \sup_{t \in I} \int f_l\log\left(e + \frac{f_l}{w}\right)dx,
 \qquad
 \mathcal{P}_l(I) = \int_I\!\int_{\{\theta > l\}} \frac{P^2}{\mu + \lambda(\rho)}\,dx\,dt.
\end{equation}
Then there exists a positive constant $C$ such that
\begin{equation}\label{vac:bandtest}
 \mathcal{P}_l(I)
 \le \frac{C(1 + a_\nu^{-1})F_l}{\kappa_0 l^{a_\nu}}
     \bigl(1 + l + \mathcal{P}_l(I)\bigr),
 \qquad
 \int_I\!\int Q\,dx\,dt
 \le C\bigl(1 + l + \mathcal{P}_l(I)\bigr).
\end{equation}
\end{lemma}

\begin{proof}
Consider the truncated function
\begin{equation}\nonumber
 z_l = \bigl(\theta^{1 - \nu/2} - (l/2)^{1 - \nu/2}\bigr)_+.
\end{equation}
Since $\theta \le L_0 < l/2$ on $\Gamma_t$, we have $z_l = 0$ on $\Gamma_t$, while $\theta^{2-\nu} \le C z_l^2$ on $\{\theta > l\}$. Applying Young's inequality $xy \le x\log(e+x) + e^y$ in the measure $w\,dx$ together with \eqref{vac:tm}, we obtain
\begin{equation}\nonumber
\begin{aligned}
 \int_{\{\theta > l\}} \frac{P^2}{\mu + \lambda(\rho)}\,dx
 &\le C \int f_l z_l^2\,dx \\
 &\le C\|\nabla z_l\|_2^2
      \left[\int f_l\log\left(e + \frac{f_l}{w}\right)dx
            + \int w\exp\left(\frac{c z_l^2}{\|\nabla z_l\|_2^2}\right)dx\right] \\
 &\le C F_l \|\nabla z_l\|_2^2,
\end{aligned}
\end{equation}
which holds trivially when $\|\nabla z_l\|_2 = 0$.
Using $2 - \nu - \alpha = 1 - a_\nu$, we deduce
\begin{equation}\nonumber
\begin{aligned}
 \|\nabla z_l\|_2^2
 &\le C \int_{\{\theta > l/2\}} \theta^{-\nu}|\nabla\theta|^2\,dx \\
 &\le \frac{C l^{1 - a_\nu}}{\kappa_0}
      \int_{\{l/2 < \theta \le 2l\}} \kappa(\theta)|\nabla\log\theta|^2\,dx
      + C \int_{\{\theta > 2l\}} \theta^{-\nu}|\nabla\theta|^2\,dx.
\end{aligned}
\end{equation}
Integrating over $I$ and applying \eqref{vac:entropy} yields
\begin{equation}\nonumber
 \mathcal{P}_l(I)
 \le C F_l \left(\frac{l^{1 - a_\nu}}{\kappa_0}
     + \int_I\!\int_{\{\theta > 2l\}} \theta^{-\nu}|\nabla\theta|^2\,dx\,dt\right).
\end{equation}

To estimate the last integral, introduce the Lipschitz cutoff
\begin{equation}\nonumber
 \chi_l(s) =
 \begin{cases}
  0, & s \le l, \\
  (s - l)/l, & l < s < 2l, \\
  1, & s \ge 2l.
 \end{cases}
\end{equation}
For any bounded Lipschitz function $b$, multiplying the temperature equation by $b(\theta)$ and using the continuity equation, we arrive at
\begin{equation}\nonumber
 c_v \left[\int \rho \int_0^\theta b(s)\,ds\,dx\right]_{\partial I}
 + \int_I\!\int \kappa(\theta) b'(\theta)|\nabla\theta|^2\,dx\,dt
 = \int_I\!\int b(\theta)(Q - P\operatorname{div} u)\,dx\,dt,
\end{equation}
where $[\cdot]_{\partial I}$ denotes the difference between the terminal and initial values on $I$.

Taking $b=\chi_l$ and using energy conservation, we obtain
\begin{equation}\nonumber
\begin{aligned}
 \int_I\!\int \chi_l(\theta) Q\,dx\,dt
 &\le C + \frac{1}{l}\int_I\!\int_{\{l < \theta < 2l\}} \kappa(\theta)|\nabla\theta|^2\,dx\,dt
        + \int_I\!\int \chi_l(\theta) P\operatorname{div} u\,dx\,dt \\
 &\le C + 4l \int_I\!\int \kappa(\theta)|\nabla\log\theta|^2\,dx\,dt
        + \frac{1}{2}\int_I\!\int \chi_l(\theta) Q\,dx\,dt
        + \frac{1}{2}\mathcal{P}_l(I),
\end{aligned}
\end{equation}
where we have used
\[
 Q\ge(\mu+\lambda)(\operatorname{div}u)^2,
 \qquad
 P|\operatorname{div}u|
 \le \frac12Q+\frac{P^2}{2(\mu+\lambda)}.
\]

Applying \eqref{vac:entropy} yields
\begin{equation}\nonumber
 \int_I\!\int_{\{\theta > 2l\}} Q\,dx\,dt
 \le C \bigl(1 + l + \mathcal{P}_l(I)\bigr).
\end{equation}
Together with
\begin{equation}\nonumber
 \int_I\!\int_{\{\theta \le 2l\}} Q\,dx\,dt
 \le 2l \int_I\!\int \frac{Q}{\theta}\,dx\,dt
 \le C l,
\end{equation}
this proves the second inequality in \eqref{vac:bandtest}.

Next, take $b(\theta) = -\chi_l(\theta)(l/\theta)^{a_\nu} \in [-1,0]$, whose derivative satisfies
\begin{equation}\nonumber
 b'(\theta) = -\chi_l'(\theta)\left(\frac{l}{\theta}\right)^{a_\nu}
              + a_\nu l^{a_\nu}\chi_l(\theta)\theta^{-a_\nu-1},
\end{equation}
so that $b'(\theta) \ge -1/l$ on $\{l < \theta < 2l\}$ and $\kappa(\theta)b'(\theta) = \kappa_0 a_\nu l^{a_\nu}\theta^{-\nu}$ on $\{\theta > 2l\}$. Moreover, completing the square yields
\begin{equation}\nonumber
 P\operatorname{div} u - Q
 \le P\operatorname{div} u - (\mu + \lambda)(\operatorname{div} u)^2
 \le \frac{P^2}{4(\mu + \lambda)}.
\end{equation}
Substituting this choice of $b$ into the multiplier identity and using \eqref{vac:entropy}, we deduce
\begin{equation}\nonumber
 \kappa_0 a_\nu l^{a_\nu} \int_I\!\int_{\{\theta > 2l\}} \theta^{-\nu}|\nabla\theta|^2\,dx\,dt
 \le C + \frac{1}{l}\int_I\!\int_{\{l < \theta < 2l\}} \kappa(\theta)|\nabla\theta|^2\,dx\,dt
       + \frac{1}{4}\mathcal{P}_l(I)
 \le C \bigl(1 + l + \mathcal{P}_l(I)\bigr).
\end{equation}
Inserting this bound into the preceding estimate for $\mathcal{P}_l(I)$ gives
\begin{equation}\nonumber
 \mathcal{P}_l(I)
 \le \frac{C F_l}{\kappa_0 l^{a_\nu}}
     \left[l + a_\nu^{-1}\bigl(1 + l + \mathcal{P}_l(I)\bigr)\right]
 \le \frac{C(1 + a_\nu^{-1})F_l}{\kappa_0 l^{a_\nu}}
     \bigl(1 + l + \mathcal{P}_l(I)\bigr),
\end{equation}
which establishes the first inequality in \eqref{vac:bandtest} and completes the proof of Lemma~\ref{vac:heat-band}.
\end{proof}

\begin{lemma}\label{vac:strict}
Under \textup{(A)} or \textup{(B)}, there exist constants $l \ge 2L_0 + 2$ and $C > 0$, independent of the density upper bound, such that
\begin{equation}\nonumber
 U(T) + \mathcal{P}_l([0,T]) \le C.
\end{equation}
\end{lemma}

\begin{proof}
Fix $p > 2$, and apply Lemma~\ref{vac:heat-band} with $\nu = 0$ and $I = [0,t]$ ($0 < t \le T$), so that $a_0 = \alpha - 1 > 0$ and
\begin{equation}\nonumber
 f_l = \mathbf{1}_{\{\theta > l\}}\frac{\rho^2}{\mu + \lambda(\rho)},
 \qquad
 F_l = 1 + \sup_{0 \le s \le t} \int f_l\log\left(e + \frac{f_l}{w}\right)dx.
\end{equation}
By \eqref{vac:struct}, we have $f_l \le C\rho D(\rho)$ and $\log(e + f_l) \le C L(\rho)$. Since $w = W^{-c} \le 1$, it follows that
\begin{equation}\nonumber
 f_l\log\left(e + \frac{f_l}{w}\right)
 \le f_l\bigl(\log(e + f_l) + c\log W\bigr)
 \le C\rho D(\rho)\bigl(L(\rho) + \log W\bigr).
\end{equation}
On $\{\log W \le L(\rho)\}$, the right-hand side is bounded by $C\rho D(\rho)L(\rho)$; on $\{\log W > L(\rho)\}$, \eqref{vac:struct} yields
\begin{equation}\nonumber
 D(\rho)\log W
 \le C L(\rho)^{\beta_D}\log W
 \le C(\log W)^{1 + \beta_D}
 \le C W^{a_0}.
\end{equation}
Using \eqref{vac:entropy}, we therefore obtain
\begin{equation}\nonumber
 F_l \le C\left(1 + \sup_{0 \le s \le t} \int \rho D(\rho)L(\rho)\,dx\right).
\end{equation}

Set $A = A(t) = \max\bigl\{2,\, N_p(t)^{2/(p-1)}\bigr\}$. On $\{\rho \le A\}$, monotonicity of $D$ and \eqref{vac:entropy} give
\begin{equation}\nonumber
 \int_{\{\rho \le A\}} \rho D(\rho)L(\rho)\,dx
 \le D(A)\int \rho L(\rho)\,dx
 \le C D(A).
\end{equation}
On $\{\rho > A\}$, since $D(r) \le C L(r)^{\beta_D}$ and $p > 2$, we have
\begin{equation}\nonumber
 \int_{\{\rho > A\}} \rho D(\rho)L(\rho)\,dx
 \le C \int_{\{\rho > A\}} \rho^{(p+1)/2}\,dx
 \le C A^{-(p-1)/2} N_p(t)
 \le C,
\end{equation}
and hence $F_l \le C\bigl(1 + D(A)\bigr)$.

Recall that $\alpha>1$ in both cases \textup{(A)} and \textup{(B)}.
We now choose
\[
 l=l(t):=2L_0+2+C_*\bigl(1+D(A)\bigr)^{1/(\alpha-1)},
\]
where $C_*$ is sufficiently large so that
\[
 \frac{C(1+(\alpha-1)^{-1})F_l}
      {\kappa_0l^{\alpha-1}}\le\frac12.
\]
It follows from \eqref{vac:bandtest} that
\[
 \mathcal P_l([0,t])
 \le\frac12\bigl(1+l+\mathcal P_l([0,t])\bigr),
\]
and hence
\[
 U(t)+\mathcal P_l([0,t])
 \le C(1+l)
 \le C\bigl(1+D(A)\bigr)^{1/(\alpha-1)}.
\]
By \eqref{vac:moment-bound} and the definition of $A$,
\[
 \log A(t)
 \le C_{p}\bigl(1+\log N_p(t)\bigr)
 \le C_1U(t)^{c_p}.
\]
Therefore,
\begin{equation}\label{vac:strict-bootstrap}
 U(t)\le C_2
 \left[1+D\left(\exp(C_1U(t)^{c_p})\right)\right]^{1/(\alpha-1)}.
\end{equation}
To close this estimate, we claim that for every $\gamma>0$,
there exist arbitrarily large $Y$ such that
\begin{equation}\label{vac:strict-profile}
 D(e^Y)\le Y^\gamma.
\end{equation}
Indeed, in case \textup{(A)}, the assumption
\[
 \int_2^\infty\frac{dr}{(e+r)L(r)D(r)}=\infty
\]
implies, by the change of variables $r=e^Y$, that
\[
 \int_2^\infty\frac{dY}{YD(e^Y)}=\infty.
\]
If \eqref{vac:strict-profile} were false for all sufficiently
large $Y$, then for some $Y_0>2$,
\[
 \int_{Y_0}^\infty\frac{dY}{YD(e^Y)}
 \le\int_{Y_0}^\infty\frac{dY}{Y^{1+\gamma}}<\infty,
\]
a contradiction. In case \textup{(B)}, \eqref{vac:strict-profile}
follows directly from $D=D_h$ and
\[
 \lim_{Y\to\infty}\frac{\log D_h(e^Y)}{\log Y}=0.
\]

Choose $\gamma>0$ such that
\[
 \frac{\gamma}{\alpha-1}<\frac1{c_p}.
\]
By \eqref{vac:strict-profile}, we can fix $Y$ sufficiently large
such that
\[
 D(e^Y)\le Y^\gamma,\qquad
 U_*:=\left(\frac{Y}{C_1}\right)^{1/c_p}
 >\max\left\{U(0),C_2(1+Y^\gamma)^{1/(\alpha-1)}\right\}.
\]
Suppose that $U(t)$ reaches $U_*$ at some time.
By continuity, there exists a first time $t_*\in(0,T]$ such that
$U(t_*)=U_*$. Applying \eqref{vac:strict-bootstrap} at $t_*$,
we obtain
\begin{align*}
 U_*
=U(t_*)
 &\le C_2
 \left[1+D\left(\exp(C_1U_*^{c_p})\right)\right]^{1/(\alpha-1)}\\
 &=C_2\bigl(1+D(e^Y)\bigr)^{1/(\alpha-1)}\\
 &\le C_2(1+Y^\gamma)^{1/(\alpha-1)}
 <U_*,
\end{align*}
which is a contradiction. Thus
\[
 \sup_{0\le t\le T}U(t)\le U_*\le C.
\]
Substituting this bound into \eqref{vac:moment-bound} gives
\[
 \sup_{0\le t\le T}A(t)\le C.
\]
We may therefore choose a fixed $l$, independent of the density
upper bound, such that the pressure term in
\eqref{vac:bandtest} can be absorbed on $I=[0,T]$.
Consequently,
\[
 U(T)+\mathcal P_l([0,T])\le C.
\]
This completes the proof of Lemma~\ref{vac:strict}.
\end{proof}

\begin{lemma}\label{vac:gain}
Under \textup{(C)}, there exist constants $l \ge 2L_0 + 2$ and $C > 0$, independent of the density upper bound, such that
\begin{equation}\nonumber
 U(T) + \mathcal{P}_l([0,T]) \le C.
\end{equation}
\end{lemma}

\begin{proof}
Set
\begin{equation}\nonumber
 s = \frac{1 - d}{2 + d},
 \qquad
 \nu = 1 - s,
 \qquad
 a_\nu = \alpha + \nu - 1 = \alpha - s \ge s > 0,
\end{equation}
and define, for a sufficiently large constant $r_* \ge 2$ above which the lower bound in \textup{(C)} holds,
\begin{equation}\nonumber
 \varphi(r) = \bigl(L(r)^{1 + d/2} - L(r_*)^{1 + d/2}\bigr)_+,
 \qquad
 \mathcal{M}(t) = \left(1 + \int \rho\varphi(\rho)^2\,dx\right)^{1/2}.
\end{equation}
By \textup{(C)}, we have
$$ r|r\varphi'(r)|^2 \le C r L(r)^d \mathbf{1}_{\{r > r_*\}} \le C\bigl(\mu + \lambda(r)\bigr).$$
It follows from the renormalized continuity equation that
\[
 \frac{d}{dt}\mathcal M^2(t)
 =-2\int\rho\varphi(\rho)\rho\varphi'(\rho)
                  \operatorname{div}u\,dx.
\]
Hence, by Cauchy--Schwarz and
$Q\ge(\mu+\lambda)(\operatorname{div}u)^2$,
\begin{align*}
 |\mathcal M'(t)|
\le \frac{1}{\mathcal M(t)}
       \left(\int\rho\varphi(\rho)^2\,dx\right)^{1/2}
       \left(\int\rho|\rho\varphi'(\rho)|^2
                    (\operatorname{div}u)^2\,dx\right)^{1/2}
\le C\left(\int Q\,dx\right)^{1/2}.
\end{align*}

Let $I=[t_1,t_2]\subset[0,T]$ and
$\mathcal M_I=\sup_{t\in I}\mathcal M(t)$.
We next estimate $F_l$ in Lemma~\ref{vac:heat-band}.
By the viscosity lower bound in \textup{(C)},
\[
 f_l=\mathbf1_{\{\theta>l\}}
          \frac{\rho^2\theta^\nu}{\mu+\lambda(\rho)}
 \le C\rho\theta^\nu L(\rho)^{-d}.
\]
Since $\theta^\nu\log(e+\theta)\le C(1+\theta)$, we have
\begin{equation}\nonumber
 f_l\log(e + f_l)
 \le C\rho\theta^\nu L(\rho)^{-d}\bigl(L(\rho) + \log(e + \theta)\bigr)
 \le C\rho\bigl(1 + \theta + \theta^\nu L(\rho)^{1 - d}\bigr).
\end{equation}
Since $(1 - d)/s = 2 + d$ and $L(r)^{2 + d} \le C\bigl(1 + \varphi(r)^2\bigr)$, H\"older's inequality and energy conservation give
\begin{equation}\nonumber
 \int f_l\log(e + f_l)\,dx
 \le C + C\left(\int \rho\theta\,dx\right)^\nu
          \left(\int \rho L(\rho)^{2 + d}\,dx\right)^s
 \le C\mathcal{M}(t)^{2s}.
\end{equation}
Similarly, by \eqref{vac:entropy},
\begin{equation}\nonumber
 \int f_l\log W\,dx
 \le C\left(\int \rho\theta\,dx\right)^\nu
      \left(\int \rho(\log W)^{1/s}\,dx\right)^s
 \le C\left(\int \rho W^{a_0}\,dx\right)^s
 \le C.
\end{equation}
Using $w = W^{-c} \le 1$ and $\mathcal{M}_I \ge 1$, we deduce
\begin{equation}\nonumber
 F_l
 \le 1 + \sup_{t \in I} \int f_l\bigl(\log(e + f_l) + c\log W\bigr)\,dx
 \le C\mathcal{M}_I^{2s}.
\end{equation}

Now take $l = C_*\mathcal{M}_I^2$ with $C_* \ge 2L_0 + 2$. Because $a_\nu \ge s$ and $\mathcal{M}_I \ge 1$, we have
\begin{equation}\nonumber
 \frac{C(1 + a_\nu^{-1})F_l}{\kappa_0 l^{a_\nu}}
 \le C C_*^{-a_\nu}\mathcal{M}_I^{2(s - a_\nu)}
 \le C C_*^{-a_\nu}
 \le \frac{1}{2}
\end{equation}
for $C_*$ sufficiently large. Applying \eqref{vac:bandtest} yields
\begin{equation}\nonumber
 \int_I\!\int Q\,dx\,dt + \mathcal{P}_l(I)
 \le C(1 + l)
 \le C\mathcal{M}_I^2.
\end{equation}
Combining this estimate with the differential bound on $\mathcal{M}$ gives, for every $t \in I$,
\begin{equation}\nonumber
 \mathcal{M}(t)
 \le \mathcal{M}(t_1) + C|I|^{1/2}\left(\int_I\!\int Q\,dx\,ds\right)^{1/2}
 \le \mathcal{M}(t_1) + C_0|I|^{1/2}\mathcal{M}_I,
\end{equation}
and hence $\mathcal{M}_I \le 2\mathcal{M}(t_1)$ whenever $|I| \le (2C_0)^{-2}$.

Subdividing $[0,T]$ into finitely many intervals of length at most $(2C_0)^{-2}$ and iterating this bound yields
$$\sup_{0 \le t \le T}\mathcal{M}(t) \le C.$$
Setting $l = C_*\sup_{0 \le t \le T}\mathcal{M}(t)^2$ in \eqref{vac:bandtest} on $I = [0,T]$ gives
$$U(T) + \mathcal{P}_l([0,T]) \le C,$$
which completes the proof of Lemma~\ref{vac:gain}.
\end{proof}

\begin{lemma}\label{vac:linear}
Under \textup{(D)}, fix $R_0 > r_*$ sufficiently large. Then there exist constants $l \ge 2L_0 + 2$ and $C > 0$, independent of the density upper bound, such that
\begin{equation}\nonumber
 \sup_{0 \le t \le T} \int_{\{\rho > 2R_0\}} \frac{\rho L(\rho)^2}{D(\rho)}\,dx
 + U(T) + \mathcal{P}_l([0,T]) \le C.
\end{equation}
\end{lemma}

\begin{proof}
Recall that $\alpha = 1$ and $D(r) = \Psi(J(r))$. Integrating $\Psi'(s)/\Psi(s) \le p/s$ and applying \textup{(D)} yield
\begin{equation}\nonumber
 \Psi(tx) \le t^p \Psi(x) \quad (t, x \ge 1),
 \qquad
 D(r) \le C J(r)^p,
 \qquad
 0 \le \frac{r D'(r)}{D(r)} \le \frac{C}{L(r)J(r)} \quad (r \ge r_*).
\end{equation}
Set
\begin{equation}\nonumber
 e_u(t) = \int \frac{Q}{\theta}\,dx,
 \qquad
 e_\theta(t) = \kappa_0 \int \frac{|\nabla\theta|^2}{\theta}\,dx,
 \qquad
 b(t) = \bigl[e_u(t)(1 + e_\theta(t))\bigr]^{1/2}.
\end{equation}
By \eqref{vac:entropy} and Young's inequality, we have
\begin{equation}\nonumber
 \int_0^T \bigl(e_u(t) + e_\theta(t) + b(t)\bigr)\,dt \le C.
\end{equation}
Throughout the proof, all constants are independent of the density upper bound and of the indices introduced below.

\smallskip\noindent\textit{Logarithmic moments of the density.}
Choose a smooth nondecreasing cutoff $\chi \in C^\infty([0,\infty))$ with $\chi = 0$ on $[0,R_0]$ and $\chi = 1$ on $[2R_0,\infty)$. For $m \ge 0$, set
\begin{equation}\nonumber
 \sigma_m = 2^{-m},
 \qquad
 w_m(r) = \frac{\chi(r)L(r)^{1 + \sigma_m}}{D(r)},
 \qquad
 H_m(t) = 1 + \int \rho w_m(\rho)\,dx.
\end{equation}
Clearly, $\sup_{m\ge0}H_m(0)\le C$.
We first prove that
\begin{equation}\label{crit:moment-rate}
 \left|\frac{d}{dt}\sqrt{H_0(t)}\right|
 \le C\left(\int Q\,dx\right)^{1/2},
 \qquad
 |H_m'(t)| \le C b(t)\sqrt{H_{m-1}(t)} \quad (m \ge 1).
\end{equation}
Indeed, from the renormalized continuity equation, we have
\begin{equation}\nonumber
 H_m'(t) = -\int \rho^2 w_m'(\rho)\operatorname{div} u\,dx,
 \qquad
 |r w_m'(r)| \le
 \begin{cases}
  C, & R_0 < r < 2R_0, \\
  C L(r)^{\sigma_m}/D(r), & r \ge 2R_0.
 \end{cases}
\end{equation}
For $m = 0$, H\"older's inequality, $Q \ge (\mu + \lambda)(\operatorname{div} u)^2$, and $\mu + \lambda(r) \ge c r/D(r)$ yield
\begin{equation}\nonumber
 |H_0'|
 \le \left(\int \frac{\rho^4|w_0'(\rho)|^2}{\mu + \lambda(\rho)}\,dx\right)^{1/2}
     \left(\int Q\,dx\right)^{1/2}
 \le C\sqrt{H_0}\left(\int Q\,dx\right)^{1/2}.
\end{equation}
For $m \ge 1$, the same estimate together with \eqref{vac:entropy} on $\{R_0 < \rho < 2R_0\}$ gives
\begin{equation}\nonumber
 |H_m'|
 \le \left(\int \frac{\rho^4\theta|w_m'(\rho)|^2}{\mu + \lambda(\rho)}\,dx\right)^{1/2} e_u^{1/2}
 \le C e_u^{1/2}\left(1 + \int_{\{\rho \ge 2R_0\}}
      \frac{\rho\theta L(\rho)^{2\sigma_m}}{D(\rho)}\,dx\right)^{1/2}.
\end{equation}
Applying \eqref{vac:tm} to $(\sqrt{\theta} - \sqrt{L_0})_+$ gives a constant $c_* > 0$ such that
\begin{equation}\nonumber
 Y = \frac{\theta}{1 + e_\theta},
 \qquad
 \int w e^{c_* Y}\,dx \le C.
\end{equation}
Set $q_m = \mathbf{1}_{\{\rho \ge 2R_0\}}\rho L(\rho)^{2\sigma_m}/D(\rho)$ and $X = \log W$. Since $2\sigma_m = \sigma_{m-1} \le 1$, splitting at $X = L(\rho)$ yields
\begin{equation}\nonumber
 q_m\log\left(e + \frac{q_m}{w}\right)
 \le C q_m\bigl(L(\rho) + X\bigr)
 \le C\rho w_{m-1}(\rho) + C\rho(1 + X^2).
\end{equation}
Applying Young's inequality and using \eqref{vac:entropy}, we obtain
\begin{equation}\nonumber
 \int q_m Y\,dx
 \le C\left(1 + \int q_m\log\left(e + \frac{q_m}{w}\right)dx\right)
 \le C H_{m-1},
\end{equation}
which gives \eqref{crit:moment-rate}.

\smallskip\noindent\textit{Estimates at successive doubling times.}
Fix $M_{\rm in} > \max\{2, H_0(0)\}$ and define, for each attained level $M_j = 2^j M_{\rm in}$ ($j \ge 1$),
\begin{equation}\nonumber
 t_j = \inf\{t \ge 0 : H_0(t) = M_j\},
 \qquad t_0 = 0,
 \qquad I_j = [t_{j-1}, t_j],
 \qquad B_j = \int_{I_j} b(t)\,dt,
 \qquad e_j = \int_{I_j} (e_u + e_\theta)\,dt.
\end{equation}
Select an integer $n_j \ge 1$ such that
\begin{equation}\nonumber
 \delta_j = 2^{-n_j} \in \left[\frac{1}{4\log(e+j)},\,\frac{1}{2\log(e+j)}\right],
 \qquad
 K_j = \sup_{0 \le t \le t_j} H_{n_j}(t),
\end{equation}
and introduce
\begin{equation}\nonumber
 a_j = \frac{\delta_j}{2(1 - \delta_j)},
 \qquad
 M^*(t) = \max\left\{M_{\rm in},\,\sup_{0 \le s \le t} H_0(s)\right\},
 \qquad
 \mathcal{V}_j(t) = 1 + \int_0^t b(s) M^*(s)^{a_j}\,ds.
\end{equation}
With $\gamma_k = 2(1 - \sigma_k)$, we claim that
\begin{equation}\label{crit:moment-induction}
 H_k(t) \le C M^*(t)^{\sigma_k - a_j\gamma_k}\mathcal{V}_j(t)^{\gamma_k},
 \qquad 0 \le k \le n_j,\quad 0 \le t \le t_j.
\end{equation}
For $k = 0$ this is immediate. Assuming \eqref{crit:moment-induction} holds at $k - 1$ with constant $C_{k-1}$, \eqref{crit:moment-rate} and $\gamma_k = 1 + \gamma_{k-1}/2$ yield
\begin{equation}\nonumber
 |H_k'(t)|
 \le C\sqrt{C_{k-1}}\,M^*(t)^{\sigma_k - a_j\gamma_k}
     \mathcal{V}_j(t)^{\gamma_k - 1}\mathcal{V}_j'(t).
\end{equation}
Since $\sigma_k - a_j\gamma_k = (\sigma_k - \delta_j)/(1 - \delta_j) \ge 0$ for $k \le n_j$ and $M^*$ is nondecreasing, integration gives
\begin{equation}\nonumber
 H_k(t)
 \le H_k(0) + \frac{C\sqrt{C_{k-1}}}{\gamma_k}
              M^*(t)^{\sigma_k - a_j\gamma_k}\mathcal{V}_j(t)^{\gamma_k}.
\end{equation}
The constants satisfy $C_k\le C(1+\sqrt{C_{k-1}})$ and are therefore
uniformly bounded. This proves \eqref{crit:moment-induction}.

At $k = n_j$, we have $\sigma_{n_j} - a_j\gamma_{n_j} = 0$ and $a_j\gamma_{n_j} = \delta_j$. Using $M^*(s) \le M_i = 2^{i-j} M_j$ on $I_i$ ($1 \le i \le j$), we obtain
\begin{equation}\nonumber
 \frac{\mathcal{V}_j(t_j)}{M_j^{a_j}}
 \le M_j^{-a_j} + \sum_{i=1}^j 2^{-a_j(j-i)} B_i
 \le M_j^{-a_j} + \sum_{i=1}^j e^{-c_1\delta_j(j-i)} B_i
 = h_j,
\end{equation}
and therefore \eqref{crit:moment-induction} implies
\begin{equation}\label{crit:moment-record}
 \frac{K_j}{M_j^{\delta_j}} \le C h_j^{2(1 - \delta_j)},
 \qquad
 h_j = M_j^{-a_j} + \sum_{i=1}^j e^{-c_1\delta_j(j-i)} B_i.
\end{equation}
Since $\delta_j \asymp 1/\log(e+j)$, splitting the sum over $j \ge i$ at $j = 2i + 4$ shows that
\begin{equation}\nonumber
 \sum_{j \ge 1} \delta_j M_j^{-a_j} \le C,
 \qquad
 \sup_{i \ge 1} \sum_{j \ge i} \delta_j e^{-c_1\delta_j(j-i)} \le C,
\end{equation}
which yields
\begin{equation}\label{crit:kernel-sum}
 \sum_{j=1}^N \delta_j h_j
 \le C\left(1 + \sum_{i=1}^N B_i\right)
 \le C.
\end{equation}

\smallskip\noindent\textit{Dissipation on temperature bands.}
Fix an attained level $M_j$ and write $\delta = \delta_j$, $K = K_j$. Choose $\ell_0 > 8\max\{R_0, L_0, 1\}$ sufficiently large and set
\begin{equation}\nonumber
 \ell = \max\left\{\ell_0,\, M_j^{1/2},\, K^{1/\delta}\Psi(e+j)^{(2+\delta)/\delta}\right\}.
\end{equation}
Interpolating $L^{1+\delta}/D$ between $L/D$ and $L^2/D$ with respect to $\chi\rho\,dx$ and using \eqref{vac:entropy} give $K \le C M_j^\delta$, whence polynomial growth of $\Psi$ implies
\begin{equation}\label{crit:cutoff-size}
 c(1 + j) \le \log(e + \ell) \le \log(e + \ell K) \le C(1 + j).
\end{equation}
For $k \ge -2$, define
\begin{equation}\nonumber
\begin{aligned}
 v_k &= 2^k \ell,
 &\qquad A_k &= \{v_k < \theta \le 2v_k\}, \\
 G_k &= \int_{I_j}\!\int_{A_k} |\nabla\theta|^2\,dx\,dt,
 &\qquad P_k &= \int_{I_j}\!\int_{A_k} \frac{P^2}{\mu + \lambda(\rho)}\,dx\,dt,
\end{aligned}
\end{equation}
and for $k \ge -1$, put $f_k = \mathbf{1}_{A_k}\rho^2/(\mu + \lambda(\rho))$ and $g_k = \sup_{t \in I_j}\|f_k/w\|_{L\log L(w\,dx)}$, where the Luxemburg norm is associated with $s\log(e+s)$. We claim that
\begin{equation}\label{low:band-pressure}
 P_k \le C g_k G_{k-1} \quad (k \ge -1),
 \qquad
 g_{-1} \le C,
 \qquad
 \sum_{k \ge 0} g_k \le \frac{C}{\delta}.
\end{equation}
Indeed, setting $z_k = (\min\{\theta, v_k\} - v_{k-1})_+$, we have $z_k = 0$ on $\Gamma_t$, $z_k = v_k/2 \ge \theta/4$ on $A_k$, and $\nabla z_k = \mathbf{1}_{A_{k-1}}\nabla\theta$. Orlicz--H\"older in $w\,dx$ and \eqref{vac:tm} yield
\begin{equation}\nonumber
 \int_{A_k} \frac{P^2}{\mu + \lambda(\rho)}\,dx
 \le C \int f_k z_k^2\,dx
 \le C g_k \|\nabla z_k\|_2^2
 = C g_k \int_{A_{k-1}} |\nabla\theta|^2\,dx,
\end{equation}
which integrates to $P_k \le C g_k G_{k-1}$.

To estimate $g_k$, set
\begin{equation}\nonumber
 r = \frac{1}{1 + \delta},
 \qquad s = 1 + r,
 \qquad u_k = v_k^{-1/3} + K^r v_k^{r-1}\Psi\bigl(\log(e + v_k K)\bigr)^s.
\end{equation}
By energy conservation, $\int_{A_k}\rho\,dx \le C/v_k$. On $A_k \cap \{\rho \le v_k\}$, monotonicity of $D$ gives
\begin{equation}\nonumber
 \int_{A_k \cap \{\rho \le v_k\}} \rho D(\rho)L(\rho)\,dx
 \le \frac{C D(v_k)L(v_k)}{v_k}
 \le C v_k^{-1/3}.
\end{equation}
On $A_k \cap \{\rho > v_k\}$, putting $Z = L(\rho)^{1+\delta}/D(\rho)$ and using $D(\rho) \le C L(\rho)^{1/2}$ yield
\begin{equation}\nonumber
 D(\rho)L(\rho) = Z^r D(\rho)^s \le C\Phi(Z),
 \qquad
 \Phi(z) = z^r \Psi\bigl(\log(e+z)\bigr)^s,
 \qquad
 \int_{\{\rho > v_k\}} \rho Z\,dx \le K.
\end{equation}
Since $\Phi$ is nondecreasing and
\begin{equation}\nonumber
 \frac{d}{dz}\log\frac{\Phi(z)}{z}
 \le \frac{1}{z}\left(r - 1 + \frac{sp}{\log(e+z)}\right) \le 0
 \qquad (z \ge v_k K),
\end{equation}
by \eqref{crit:cutoff-size} for $\ell_0$ sufficiently large, splitting at $Z = v_k K$ gives
\begin{equation}\nonumber
 \int_{A_k \cap \{\rho > v_k\}} \rho D(\rho)L(\rho)\,dx
 \le C\Phi(v_k K)\int_{A_k} \rho\,dx
   + C\frac{\Phi(v_k K)}{v_k K}\int_{\{\rho > v_k\}} \rho Z\,dx
 \le C K^r v_k^{r-1}\Psi\bigl(\log(e + v_k K)\bigr)^s.
\end{equation}
Combining both density regions yields $\int_{A_k}\rho D(\rho)L(\rho)\,dx \le C u_k$. Using $f_k \le C\mathbf{1}_{A_k}\rho D(\rho)$ and $u_k \ge v_k^{-1/3}$, we obtain
\begin{equation}\nonumber
 \int f_k\log\left(e + \frac{f_k}{u_k}\right)dx
 \le C \int_{A_k} \rho D(\rho)\bigl(L(\rho) + L(v_k)\bigr)\,dx
 \le C u_k.
\end{equation}
For the spatial weight, splitting at $\log W = L(\rho)$ and applying \eqref{vac:entropy} give
\begin{equation}\nonumber
 \int f_k\log W\,dx
 \le C\int_{A_k}\rho D(\rho)L(\rho)\,dx
   + C\left(\int_{A_k}\rho\,dx\right)^{1/2}\left(\int \rho W^{a_0}\,dx\right)^{1/2}
 \le C(u_k + v_k^{-1/2})
 \le C u_k.
\end{equation}
Hence $\int f_k\log\bigl(e + f_k/(u_k w)\bigr)\,dx \le C u_k$, and thus $g_k \le C u_k$. Because $K^r \ell^{r-1} \le \Psi(e+j)^{-s}$, \eqref{crit:cutoff-size} and $\Psi(tx) \le t^p\Psi(x)$ imply
\begin{equation}\nonumber
 u_{-1} \le C,
 \qquad
 u_k \le v_k^{-1/3} + C 2^{-k\delta/2}\left(1 + \frac{k}{e+j}\right)^{2p}
 \quad (k \ge 0).
\end{equation}
Using $(e+j)^{-1} \le C\delta$, summation over $k \ge 0$ gives $\sum_{k \ge 0} u_k \le C/\delta$, completing the proof of \eqref{low:band-pressure}.

Next, for $k \ge 0$, define the piecewise logarithmic hat function
\begin{equation}\nonumber
 b_k(\theta) = \left(1 - \left|\log_2\frac{\theta}{v_k}\right|\right)_+,
 \qquad
 \kappa(\theta)b_k'(\theta) = \frac{\kappa_0}{\log 2}\bigl(\mathbf{1}_{A_{k-1}} - \mathbf{1}_{A_k}\bigr).
\end{equation}
Testing the temperature equation by $-b_k(\theta)$ over $I_j$ and applying Young's inequality yield
\begin{equation}\nonumber
 \frac{\kappa_0}{\log 2}(G_k - G_{k-1})
 + \frac{1}{2}\int_{I_j}\!\int b_k(\theta)Q\,dx\,dt
 \le a_k + C(P_{k-1} + P_k),
\end{equation}
where
\begin{equation}\nonumber
 a_k = c_v\left|\left[\int \rho\int_0^\theta b_k(s)\,ds\,dx\right]_{\partial I_j}\right|,
 \qquad
 \sum_{k \ge 0} a_k \le 2E_0.
\end{equation}
Since $G_{-2} + G_{-1} \le C\ell e_j$, setting $\overline{G}_k = \max_{-2 \le i \le k} G_i$ and substituting \eqref{low:band-pressure} give
\begin{equation}\nonumber
 \overline{G}_k \le \bigl(1 + C(g_{k-1} + g_k)\bigr)\overline{G}_{k-1} + C a_k,
\end{equation}
which iterates via \eqref{low:band-pressure} to $\sup_{k \ge -1} G_k \le C e^{C/\delta}(1 + \ell e_j)$. Consequently,
\begin{equation}\nonumber
 \mathcal{P}_\ell(I_j)
 = \sum_{k \ge 0} P_k
 \le C\left(\sum_{k \ge 0} g_k\right)\sup_{k \ge -1} G_k
 \le C e^{C/\delta}(1 + \ell e_j).
\end{equation}
Testing the temperature equation by $\chi_\ell(\theta)$ on $I_j$ as in Lemma~\ref{vac:heat-band} then gives
\begin{equation}\label{crit:heat-window}
 \int_{I_j}\!\int Q\,dx\,dt
 \le C\bigl(1 + \ell e_j + \mathcal{P}_\ell(I_j)\bigr)
 \le e^{C/\delta}(1 + \ell e_j).
\end{equation}

\smallskip\noindent\textit{Boundedness of the logarithmic moment.}
Since $H_0(t_j) = M_j \ge 2 H_0(t_{j-1})$, integrating the first bound of \eqref{crit:moment-rate} over $I_j$ and applying \eqref{crit:heat-window} yield
\begin{equation}\nonumber
 c M_j
 \le \bigl(\sqrt{H_0(t_j)} - \sqrt{H_0(t_{j-1})}\bigr)^2
 \le C|I_j|\int_{I_j}\!\int Q\,dx\,dt
 \le C|I_j|e^{C/\delta_j}(1 + \ell e_j).
\end{equation}
Because $M_j = 2^j M_{\rm in}$ and $\delta_j^{-1} \le 4\log(e+j)$, we have $e^{C/\delta_j}(1 + \ell_0 + M_j^{1/2})/M_j \to 0$ as $j \to \infty$. Substituting the definition of $\ell$ therefore gives, for all sufficiently large attained $j$,
\begin{equation}\nonumber
 1 \le C e^{C/\delta_j}|I_j|e_j
       \left(\frac{K_j}{M_j^{\delta_j}}\right)^{1/\delta_j}
       \Psi(e+j)^{(2 + \delta_j)/\delta_j}.
\end{equation}
Raising both sides to the power $\delta_j/2$, applying \eqref{crit:moment-record}, multiplying by $\delta_j$, and using Young's inequality, we deduce
\begin{equation}\nonumber
 \delta_j \Psi(e+j)^{-1 - \delta_j/2}
 \le C\delta_j |I_j|^{\delta_j/2} e_j^{\delta_j/2} h_j^{1 - \delta_j}
 \le C\bigl(|I_j| + e_j + \delta_j h_j\bigr).
\end{equation}
Since $\Psi(e+j)^{\delta_j/2} \le C$, summing over $1 \le j \le N$ and using \eqref{crit:kernel-sum} give, for every attained level $M_N$,
\begin{equation}\nonumber
 \sum_{j=1}^N \frac{1}{\Psi(e+j)\log(e+j)}
 \le C\left(1 + \sum_{j=1}^N (|I_j| + e_j + \delta_j h_j)\right)
 \le C.
\end{equation}
By \textup{(D)}, the series $\sum_{j \ge 1}[\Psi(e+j)\log(e+j)]^{-1}$ diverges, so the number of attained doubling levels is bounded and
\begin{equation}\nonumber
 \sup_{0 \le t \le T} \int_{\{\rho > 2R_0\}} \frac{\rho L(\rho)^2}{D(\rho)}\,dx
 \le \sup_{0 \le t \le T} H_0(t) \le C.
\end{equation}

Finally, since $D(r)\le CJ(r)^p\le CL(r)^{1/4}$, we have
\[
 \sup_{0\le t\le T}\int\rho L(\rho)^{3/2}\,dx\le C.
\]
Taking $\nu=1/12$ in Lemma~\ref{vac:heat-band}, so that
$a_\nu=\nu>0$, we obtain
\[
 f_l\le C\rho L(\rho)^{1/4}\theta^\nu.
\]
Using Young's inequality and splitting the spatial term at
$\log W=L(\rho)$, we have
\[
 f_l\log(e+f_l/w)
 \le C\rho\left(
 1+\theta+W^{a_0}
 +L(\rho)^{\frac{5}{4(1-\nu)}}
 +L(\rho)^{\frac{1}{4(1-2\nu)}}
 \right).
\]
Since both logarithmic exponents are less than $3/2$,
it follows from \eqref{vac:entropy} that
$F_l\le C$ independently of $l$.
Choosing $l$ sufficiently large in \eqref{vac:bandtest}, we conclude
\[
 U(T)+\mathcal P_l([0,T])\le C.
\]

We also record for later use that $D(r)\le CJ(r)^p$ implies
$\log D(r)=o(\log L(r))$, and, by the substitution $x=J(r)$,
\[
 \int_{r_*}^\infty\frac{dr}{(e+r)L(r)D(r)}
 \asymp\int_{J(r_*)}^\infty\frac{dx}{\Psi(x)}
 =\infty
\]
by \textup{(D)}.
This completes the proof of Lemma~\ref{vac:linear}.
\end{proof}

\begin{lemma}\label{vac:thermal-tail}
For every $p > 1$ and $0 < s < a_0$, one has
\begin{equation}\label{vac:density-weight}
 \sup_{0 \le t \le T} \int \rho^p W^s\,dx
 + \int_0^T\!\int \left(Q + \frac{P^2}{h}\right)dx\,dt
 \le C_{p,s}.
\end{equation}
Set $m = 1 + \alpha$ and $H = c_v\int \rho\int_0^\theta K(v)\,dv\,dx$. Then, for $q \ge \max\{2,m\}$, $0 < \epsilon \le \min\{1/2, m/2\}$, and $0 < \zeta \le 1/2$,
\begin{equation}\label{vac:thermal-tails}
 \int_0^T \|\theta\|_{L^q(w\,dx)}^{m - \epsilon}\,dt
 \le C_{c}\frac{q}{\epsilon},
 \qquad
 \int_0^T (1 + H)^{1 - \zeta}\,dt
 \le C \zeta^{-2}.
\end{equation}
Moreover, if $\rho \le R$ on $\mathbb{R}^2 \times [0,T]$, then
\begin{equation}\label{vac:floor}
 \theta(x,t) \ge
 \left(\vartheta_*^{-1}
       + \frac{t}{4c_v}\sup_{0 \le r \le R}\frac{r}{\mu + \lambda(r)}
 \right)^{-1}.
\end{equation}
\end{lemma}

\begin{proof}
By Lemmas~\ref{vac:strict}--\ref{vac:linear}, there exists a fixed $l \ge 2L_0 + 2$ such that
$$U(T) + \mathcal{P}_l([0,T]) \le C.$$
This, together with Lemma~\ref{vac:moment} and conservation
of mass, yields for any $p>1$,
\[
 \sup_{0\le t\le T}\int\rho^p\,dx\le C_{p}.
\]
Applying H\"older's inequality and \eqref{vac:entropy}, we obtain
\begin{equation}\nonumber
 \int \rho^p W^s\,dx
 \le \left(\int \rho W^{a_0}\,dx\right)^{s/a_0}
     \left(\int \rho^{\frac{p - s/a_0}{1 - s/a_0}}\,dx\right)^{1 - s/a_0}
 \le C_{p,s}.
\end{equation}
Moreover, splitting at $\theta = l$ yields
\begin{equation}\nonumber
 \int_0^T\!\int \frac{P^2}{h}\,dx\,dt
 \le C_l \int_0^T\!\int \rho^2\,dx\,dt + \mathcal{P}_l([0,T])
 \le C,
\end{equation}
which establishes \eqref{vac:density-weight}.

Next, testing the temperature equation by $-\chi_l(\theta)(l/\theta)^\epsilon$ as in the proof of Lemma~\ref{vac:heat-band} gives
\begin{equation}\nonumber
 \kappa_0\epsilon l^\epsilon
 \int_0^T\!\int_{\{\theta > 2l\}} \theta^{\alpha - 1 - \epsilon}|\nabla\theta|^2\,dx\,dt
 \le C + C l \int_0^T\!\int \kappa(\theta)|\nabla\log\theta|^2\,dx\,dt
       + C\mathcal{P}_l([0,T])
 \le C.
\end{equation}
Putting $b = (m - \epsilon)/2$ and $z = \bigl(\theta^b - (2l)^b\bigr)_+$, we note that $2b - 2 = \alpha - 1 - \epsilon$ and hence
\begin{equation}\nonumber
 \int_0^T \|\nabla z\|_2^2\,dt
 = b^2 \int_0^T\!\int_{\{\theta > 2l\}} \theta^{\alpha - 1 - \epsilon}|\nabla\theta|^2\,dx\,dt
 \le \frac{C}{\epsilon}.
\end{equation}
Since $z = 0$ on $\Gamma_t$ and $q/b \ge 2$, applying \eqref{vac:sob} yields
\begin{equation}\nonumber
 \|\theta\|_{L^q(w\,dx)}^{m - \epsilon}
 = \|\theta^b\|_{L^{q/b}(w\,dx)}^2
 \le 2\|z\|_{L^{q/b}(w\,dx)}^2 + 2(2l)^{2b}\|1\|_{L^{q/b}(w\,dx)}^2
 \le C_{c}\bigl(1 + q\|\nabla z\|_2^2\bigr).
\end{equation}
Integrating over $(0,T)$ gives the first estimate in \eqref{vac:thermal-tails}.

To estimate $H$, we note that
\[
 H=\frac{c_v\kappa_0}{m(m+1)}\int\rho\theta^{m+1}\,dx.
\]
Choose $p=1+a_0/(2c)$ and let $p'=p/(p-1)$.
Since $c(p-1)=a_0/2$, \eqref{vac:density-weight} yields
\[
 \|\rho W^{c/p'}\|_p^p
 =\int\rho^pW^{a_0/2}\,dx\le C.
\]
For any $0 < b_1 < 1$, two applications of H\"older's inequality yield
\begin{equation}\nonumber
\begin{aligned}
 \int \rho\theta^{m+1}\,dx
 &\le \left(\int \rho\theta\,dx\right)^{1 - b_1}
      \left(\int \rho\theta^{(m + b_1)/b_1}\,dx\right)^{b_1} \\
 &\le \left(\int \rho\theta\,dx\right)^{1 - b_1}
      \|\rho W^{c/p'}\|_p^{b_1} \|\theta\|_{L^q(w\,dx)}^{m + b_1}
  \le C\|\theta\|_{L^q(w\,dx)}^{m + b_1},
 \qquad q = \frac{p'(m + b_1)}{b_1}.
\end{aligned}
\end{equation}
For $\zeta > 0$ sufficiently small, choosing
\begin{equation}\nonumber
 \sigma = 1 - \zeta,
 \qquad
 b_1 = \frac{m\zeta}{2(1 - \zeta)},
 \qquad
 \epsilon = \frac{m\zeta}{2}.
\end{equation}
Then
\[
 \sigma(m+b_1)=m-\epsilon,\qquad
 q=\frac{p'(2-\zeta)}{\zeta},
 \qquad \frac{q}{\epsilon}\le C\zeta^{-2}.
\]
Consequently, the first estimate in \eqref{vac:thermal-tails}
implies
\[
 \int_0^T(1+H)^{1-\zeta}\,dt
 \le C\left(
 1+\int_0^T\|\theta\|_{L^q(w\,dx)}^{m-\epsilon}\,dt\right)
 \le C\zeta^{-2}.
\]
The remaining values of $\zeta \in (0, 1/2]$ follow from the monotonicity of $(1 + H)^{1-\zeta}$ in $\zeta$, which proves \eqref{vac:thermal-tails}.

Finally, using $Q \ge (\mu + \lambda)(\operatorname{div} u)^2$ and $\rho \le R$, we have
\begin{equation}\nonumber
 Q - P\operatorname{div} u
 \ge -\frac{P^2}{4(\mu + \lambda(\rho))}
 = -\frac{\rho^2\theta^2}{4(\mu + \lambda(\rho))}
 \ge -c_v a_R\rho\theta^2,
 \qquad
 a_R = \frac{1}{4c_v}\sup_{0 \le r \le R}\frac{r}{\mu + \lambda(r)}.
\end{equation}
Set $y(t) = (\vartheta_*^{-1} + a_R t)^{-1}$, so that $y' = -a_R y^2$. Subtracting $c_v\rho y'$ from the temperature equation yields
\begin{equation}\nonumber
 c_v\rho D_t(\theta - y) - \Delta K(\theta)
 \ge c_v a_R\rho(y^2 - \theta^2).
\end{equation}
Multiplying by $-v$, where $v = (y - \theta)_+$, and using the continuity equation, we deduce
\begin{equation}\nonumber
 \frac{c_v}{2}\frac{d}{dt}\int \rho v^2\,dx
 + \int \kappa(\theta)|\nabla v|^2\,dx
 \le -c_v a_R\int \rho(\theta + y)v^2\,dx \le 0.
\end{equation}
Since $v(0) = 0$, integrating in time gives $\sqrt{\rho}\,v = 0$ and $\nabla v = 0$ almost everywhere. Because $\int \rho\,dx = m_0 > 0$, it follows that $v = 0$, which proves \eqref{vac:floor} and completes the proof of Lemma~\ref{vac:thermal-tail}.
\end{proof}

To obtain the density upper bound, we recall the effective viscous flux $G = h(\rho)\operatorname{div} u - P$, the vorticity $\omega = \partial_1 u^2 - \partial_2 u^1$, and $m = 1 + \alpha$, and use the convention $(\operatorname{div}\mathcal{T})_i = \sum_{j=1}^2 \partial_j \mathcal{T}_{ij}$ for a matrix-valued field $\mathcal{T}$. By Lemmas~\ref{vac:low-anchor} and~\ref{vac:thermal-tail},
\begin{equation}\label{vh:inputs}
\begin{aligned}
 &\sup_{0 \le t \le T} \int \rho\bigl(1 + W^{a_0} + \theta + |u|^2\bigr)\,dx
  + \int_0^T\!\int \left(Q + \frac{P^2}{h}\right)dx\,dt \le C, \\
 &\sup_{0 \le t \le T} \int \rho^p\,dx \le C_p \quad (1 < p < \infty),
  \qquad \int \rho\,dx = m_0 > 0,
\end{aligned}
\end{equation}
where here and below constants may depend on $T$. With $b_h$ and $D$ from Lemma~\ref{pre:growth}, we introduce
\begin{equation}\nonumber
\begin{gathered}
 A_1^2 = 1 + \int \left(\frac{G^2}{h} + \mu\omega^2\right)dx,
 \qquad A_2^2 = \int \rho|\dot{u}|^2\,dx,
 \qquad Z = \|\nabla K\|_2, \\
 B = \int \rho\kappa(\theta)|\dot{\theta}|^2\,dx,
 \qquad H = c_v \int \rho \int_0^\theta K(z)\,dz\,dx,
 \qquad D_u = \int \bigl(h(\operatorname{div} u)^2 + \mu\omega^2\bigr)\,dx.
\end{gathered}
\end{equation}
Fix $\ell > \max\{6, m, 1 + 2/\alpha\}$ and a sufficiently large exponent $p_\theta > 2\ell$, and set
\begin{equation}\nonumber
 b_\ell = \frac{\ell + 1}{\ell - 1} < m,
 \qquad \Theta = 1 + \|\theta\|_{p_\theta,w},
 \qquad g = A_1^2 + \Theta^{b_\ell},
 \qquad s = \frac{2}{2 + \alpha},
\end{equation}
where $\|\cdot\|_{p,w} = \|\cdot\|_{L^p(w\,dx)}$. Choose $\sigma \in \bigl((1 + s)/2, 1\bigr)$ sufficiently close to $1$ so that
\begin{equation}\nonumber
\begin{cases}
 \displaystyle \sigma > \frac{3}{2 + \alpha}, & \alpha > 1, \\[4pt]
 0 < r_c = 3 - (2 + \alpha)\sigma < 1, & 0 < \alpha \le 1,
\end{cases}
\end{equation}
and define
\begin{equation}\nonumber
 S = 1 + A_1^2 + H,
 \qquad S_0 = 1 + D_u + (1 + H)^\sigma,
 \qquad g_* = 1 + g + (1 + H)^{(2s - 1)_+}.
\end{equation}

\begin{lemma}\label{vh:pressure}
For all sufficiently small $\zeta > 0$, one has
\begin{equation}\label{vh:tails}
 \int_0^T (g_* + S_0 + F_\alpha)\,dt \le C,
 \qquad
 \int_0^T S^{1 - \zeta}\,dt \le C\zeta^{-2},
\end{equation}
where $F_\alpha = 0$ if $\alpha > 1$, and $F_\alpha = \|\theta\|_{p_c,w}^{r_c}$ for a sufficiently large fixed $p_c$ if $0 < \alpha \le 1$. Moreover,
\begin{align}
 &\int \frac{P^2}{h}\,dx \le C(1 + H)^s,
  \qquad
  \int \frac{P^3}{h^2}\,dx \le C(1 + H)^\sigma(1 + F_\alpha),
  \label{vh:p23} \\
 &\|\nabla u\|_2^2 + A_1^2 \le C S_0,
  \qquad S_0 \le C S,
  \qquad S \le C S_0^{1/\sigma},
  \qquad \|\nabla u\|_2^4 \le C\bigl(A_1^4 + (1 + H)^{2s}\bigr).
  \label{vh:compare}
\end{align}
\end{lemma}

\begin{proof}
Since $\theta\le L_0$ on $\Gamma_t$ and $\Gamma_t\subset B_N$,
we have
\[
 \|K\|_{L^2(\Gamma_t)}
 \le K(L_0)|B_N|^{1/2}\le C.
\]
Applying \eqref{vac:sob} to $K=\kappa_0\theta^m/m$, we obtain
\begin{equation}\label{vh:anchor}
 \|K\|_{p,w}\le C\sqrt p(1+Z)\quad(p\ge2),
 \qquad \Theta\le C(1+Z)^{1/m},
\end{equation}
where the second inequality follows by taking $p=p_\theta/m$.

We next estimate the pressure terms. Recall that
\[
 H=\frac{c_v\kappa_0}{m(m+1)}\int\rho\theta^{m+1}\,dx,
 \qquad \frac{\rho}{h(\rho)}\le CD(\rho).
\]
For $1 < i < m + 1 = 2 + \alpha$, H\"older's inequality, \eqref{vh:inputs}, and $D(r) \le C L(r)^{\beta_D}$ give
\begin{equation}\nonumber
 \int \frac{P^i}{h^{i-1}}\,dx
 \le C \int \rho D(\rho)^{i-1}\theta^i\,dx
 \le C\left(\int \rho\theta^{m+1}\,dx\right)^{\frac{i}{m+1}}
      \left(\int \rho D(\rho)^{\frac{(i-1)(m+1)}{m+1-i}}\,dx\right)^{1 - \frac{i}{m+1}}
 \le C H^{\frac{i}{m+1}}.
\end{equation}
Taking $i=2$ gives
\[
 \int\frac{P^2}{h}\,dx\le C(1+H)^s.
\]
If $\alpha>1$, taking $i=3$ and using $\sigma>3/(m+1)$ yields
\[
 \int\frac{P^3}{h^2}\,dx
 \le CH^{\frac{3}{m+1}}\le C(1+H)^\sigma.
\]
For $0<\alpha\le1$, we use $(m+1)\sigma+r_c=3$ to obtain
\begin{align*}
 \int\frac{P^3}{h^2}\,dx
\le C\int
       (\rho\theta^{m+1})^\sigma
       \rho^{1-\sigma}D(\rho)^2\theta^{r_c}\,dx
\le CH^\sigma
 \left(\int\rho D(\rho)^{\frac{2}{1-\sigma}}
              \theta^{\frac{r_c}{1-\sigma}}\,dx
 \right)^{1-\sigma}.
\end{align*}
Choose $k > 1$ sufficiently large so that $c/(k - 1) < a_0$ and $p_c = k r_c/(1 - \sigma) \ge \max\{2, m\}$, and put $k' = k/(k - 1)$. By H\"older's inequality, \eqref{vac:density-weight}, and the logarithmic growth of $D$,
\begin{equation}\nonumber
 \int \rho D(\rho)^{\frac{2}{1-\sigma}}\theta^{\frac{r_c}{1-\sigma}}\,dx
 \le \left\|\rho D(\rho)^{\frac{2}{1-\sigma}} W^{c/k}\right\|_{k'}
     \left(\int w\theta^{p_c}\,dx\right)^{1/k}
 \le C\|\theta\|_{p_c,w}^{\frac{r_c}{1-\sigma}},
\end{equation}
which gives
$$\int \frac{P^3}{h^2}\,dx \le C H^\sigma \|\theta\|_{p_c,w}^{r_c}$$
and thus proves \eqref{vh:p23}.

From $h(\rho)\operatorname{div} u = G + P$ and \eqref{vh:p23}, we obtain
\begin{equation}\nonumber
 D_u \le C\left(A_1^2 + \int \frac{P^2}{h}\,dx\right)
     \le C\bigl(A_1^2 + (1 + H)^s\bigr),
 \qquad
 A_1^2 \le C\left(1 + D_u + \int \frac{P^2}{h}\,dx\right)
       \le C\bigl(1 + D_u + (1 + H)^s\bigr).
\end{equation}
Moreover, the div--curl identity and $h\ge2\mu$ imply
\[
 \|\nabla u\|_2^2
 =\|\operatorname{div}u\|_2^2+\|\omega\|_2^2
 \le CD_u.
\]
Since $s<\sigma<1$, these estimates directly yield \eqref{vh:compare}.

Finally, integration by parts gives
$$2\int |D(u)|^2\,dx = \int\bigl(2(\operatorname{div} u)^2 + \omega^2\bigr)\,dx,$$
and hence $D_u = \int Q\,dx$. By \eqref{vh:inputs},
\begin{equation}\nonumber
 \int_0^T (D_u + A_1^2)\,dt
 \le C\left(1 + \int_0^T\!\int \left(Q + \frac{P^2}{h}\right)dx\,dt\right)
 \le C.
\end{equation}
Since $b_\ell < m$, $r_c < 1 < m$ when $F_\alpha \ne 0$, and $(2s - 1)_+ < \sigma < 1$, applying \eqref{vac:thermal-tails} gives
\begin{equation}\nonumber
 \int_0^T \bigl(\Theta^{b_\ell} + F_\alpha + (1 + H)^\sigma + (1 + H)^{(2s-1)_+}\bigr)\,dt \le C,
\end{equation}
which implies the first inequality in \eqref{vh:tails}.
For $0 < \zeta \le 1/2$, combining the above bounds with \eqref{vac:thermal-tails} yields
\begin{equation}\nonumber
 \int_0^T S^{1 - \zeta}\,dt
 \le C \int_0^T \bigl(1 + A_1^2 + (1 + H)^{1 - \zeta}\bigr)\,dt
 \le C(1 + \zeta^{-2})
 \le C\zeta^{-2}.
\end{equation}
This completes the proof of Lemma~\ref{vh:pressure}.
\end{proof}

\begin{lemma}\label{vh:products}
Set $X = e + A_1^2 + A_2^2 + Z^2$.
Then there exists a fixed constant $n > 0$ such that, for every $0 < z \le 1$,
\begin{align}
& \int K\bigl(Q + |P\operatorname{div} u|\bigr)\,dx
 \le z(A_2^2 + Z^2) + C z^{-n} A_1^4 \log(e + A_1^2),
 \label{vh:k1} \\
& \int K^2\bigl(|\nabla u|^2 + h|\operatorname{div} u|^2\bigr)\,dx
 \le C g X \log X.
 \label{vh:k2}
\end{align}
\end{lemma}

\begin{proof}
Write $\mathcal{R} = \nabla^2\Delta^{-1}$ and decompose
\begin{equation}\nonumber
 M = \bigl(\nabla u - \mathcal{R}(P/h),\, G/\sqrt{h}\bigr),
 \qquad
 N = \bigl(\mathcal{R}(P/h),\, P/\sqrt{h}\bigr).
\end{equation}
Since
\[
 h\operatorname{div}u=G+P,
 \qquad
 \|\nabla u\|_2^2
 =\|\operatorname{div}u\|_2^2+\|\omega\|_2^2,
\]
the div--curl decomposition gives
\begin{equation}\label{vh:MN}
 Q + |P\operatorname{div} u| + |\nabla u|^2 + h|\operatorname{div} u|^2
 \le C(|M|^2 + |N|^2),
 \qquad
 \|M\|_2 \le C A_1.
\end{equation}

We first derive weighted and $L^4$ bounds for $M$.
The momentum equation gives
\[
 f:=\rho\dot u=\nabla G+\mu\nabla^\perp\omega.
\]
Fix $0 < \delta \ll 1$ and put $p = 4/(4 - \delta)$. By Cauchy--Schwarz and \eqref{vac:density-weight},
\begin{equation}\label{vh:weighted-f}
 \|W^{2\delta} f\|_1 + \|W^{2\delta} f\|_p + \|f\|_{4/3} \le C A_2.
\end{equation}
Let $H_1$ be any component of the kernel of $\nabla\Delta^{-1}$ and set $I_1 = (-\Delta)^{-1/2}$. For $|x| > 2$, subtracting $H_1(x)\int f\,dy = 0$, splitting at $|y| = |x|/2$, and using $|H_1(x - y) - H_1(x)| \le C|y|^{2\delta}|x|^{-1-2\delta}$ on $\{|y| < |x|/2\}$, we obtain
\begin{equation}\nonumber
 |(H_1 * f)(x)|
 \le C A_2 |x|^{-1 - 2\delta}
   + C |x|^{-2\delta} I_1(W^{2\delta}|f|)(x).
\end{equation}
By the Hardy--Littlewood--Sobolev inequality, $I_1 : L^p(\mathbb{R}^2) \to L^{4/(2-\delta)}(\mathbb{R}^2)$. Since $W^{-\delta} \in L^{4/\delta}(\mathbb{R}^2)$, combining this pointwise bound on $\{|x| > 2\}$ with $\|\nabla\Delta^{-1} f\|_4 \le C\|f\|_{4/3}$ and \eqref{vh:weighted-f} gives
\begin{equation}\nonumber
 \|W^\delta \nabla\Delta^{-1} f\|_2 + \|\nabla\Delta^{-1} f\|_4 \le C A_2.
\end{equation}
The potential representations of $G$ and $\omega$ are unique, since any difference is harmonic and lies in $L^2 + L^4$, hence vanishes. As $W^{2\delta} \in A_2$, the weighted div--curl estimate yields
\begin{equation}\label{vh:M-high}
 \|W^\delta M\|_2 + \|M\|_4 \le C A_2.
\end{equation}

For $b$ sufficiently large, set $p_b = 2b/(b - 1)$, $\eta = c/(2\delta b)$, and $\xi = 2/b$. Interpolating between $L^2$, $L^2(W^{2\delta}\,dx)$, and $L^4$ and applying \eqref{vh:MN} and \eqref{vh:M-high}, we find a fixed constant $C_1 > 0$ such that
\begin{equation}\label{vh:M-interpolation}
 \|W^{c/(2b)} M\|_{p_b}^2
 \le \|M\|_2^{2(1 - \eta - \xi)} \|W^\delta M\|_2^{2\eta} \|M\|_4^{2\xi}
 \le C A_1^2 (1 + A_2)^{C_1/b}.
\end{equation}
By weighted H\"older's inequality and \eqref{vh:anchor}, we have
\begin{equation}\nonumber
 \int K|M|^2\,dx \le C\sqrt{b}(1 + Z)\|W^{c/(2b)} M\|_{p_b}^2,
 \qquad
 \int K^2|M|^2\,dx \le C b(1 + Z)^2\|W^{c/(2b)} M\|_{p_b}^2.
\end{equation}
Choosing $b = b_* + \log(e + A_1^2)$ with $b_* > 4C_1$ in the first inequality, applying Young's inequality, and choosing $b = b_* + \log(e + A_2^2)$ in the second inequality, we deduce from \eqref{vh:M-interpolation} that
\begin{equation}\label{vh:M-products}
\begin{aligned}
 \int K|M|^2\,dx &\le z(Z^2 + A_2^2) + C z^{-n} A_1^4 \log(e + A_1^2), \\
 \int K^2|M|^2\,dx &\le C A_1^2 X \log X.
\end{aligned}
\end{equation}

We now estimate the contributions of $N$. By \eqref{vac:density-weight}, energy conservation, and H\"older's inequality, decreasing $\delta$ and enlarging $p_\theta$ if necessary give
$$\|W^\delta P\|_\ell \le C\Theta, \quad \|P\|_1 \le C.$$
Interpolating these bounds at $p_b = 2b/(b - 1)$ with $t = \ell(1 - 1/p_b)/(\ell - 1)$ yields
$$\|W^{\delta t} P\|_{p_b} \le C\Theta^t.$$
Since $h \ge 2\mu$ and $W^{p_b\delta t} \in A_{p_b}$ for $\delta$ small, the weighted Calder\'on--Zygmund estimate transfers this bound to $N$. Fixing $b \ge \ell$ large enough that $c/(2b) \le \delta t$, we have
\begin{equation}\nonumber
 \|W^{c/(2b)} N\|_{p_b} \le C\Theta^t,
 \qquad
 2t \le \frac{\ell + 1}{\ell - 1} = b_\ell.
\end{equation}
Weighted H\"older's inequality and \eqref{vh:anchor} therefore imply
\begin{align}
 \int K|N|^2\,dx &\le C(1 + Z)\Theta^{b_\ell}, \label{vh:N1}
\end{align}
and
\begin{align}
 \int K^2|N|^2\,dx &\le C(1 + Z)^2\Theta^{b_\ell}. \label{vh:N2}
\end{align}
Since $\Theta \le C(1 + Z)^{1/m}$ and $b_\ell < m$, Young's inequality gives
\begin{equation}\nonumber
 \int K|N|^2\,dx \le z Z^2 + C z^{-n},
 \qquad
 \int K^2|N|^2\,dx \le C X \Theta^{b_\ell}.
\end{equation}
Combining \eqref{vh:MN}, \eqref{vh:M-products}, \eqref{vh:N1}, and \eqref{vh:N2} establishes \eqref{vh:k1} and \eqref{vh:k2}. This completes the proof of Lemma~\ref{vh:products}.
\end{proof}

Fix $r \ge 2$ such that $\rho \le r$ on the time interval under consideration, and set
\begin{equation}\nonumber
 R = 2 + r,
 \qquad D_R = 1 + D(r),
 \qquad L_R = R\log(e + R).
\end{equation}
By \eqref{vac:floor} and \eqref{vac:struct}, we have
\begin{equation}\label{vh:floor}
 \theta^{-1} \le (\inf\theta_0)^{-1} + C T D_R,
 \qquad
 \Xi = \max\left\{1,\, \sup \frac{\rho}{h\kappa}\right\}
 \le C D_R^{1 + \alpha}.
\end{equation}
Introduce the corrected acceleration and energy functionals
\begin{equation}\nonumber
\begin{gathered}
 v = 2\nabla\Delta^{-1}\det\nabla u,
 \qquad w_1 = \dot{u} + v,
 \qquad E_w = \int \rho|w_1|^2\,dx, \\
 D_w = \int \bigl(h|\operatorname{div} w_1|^2 + \mu|\operatorname{curl} w_1|^2\bigr)\,dx,
 \qquad E_\theta = \frac{1}{2} Z^2 - \int K Q\,dx.
\end{gathered}
\end{equation}

\begin{lemma}\label{vh:diff}
There exist fixed constants $n > 0$ and $c_1 > 0$ such that, for any $0 < z \le 1$ on the sublevel $\rho \le r$,
\begin{align}
 (A_1^2)' + A_2^2
 &\le z Z^2 + C z^{-n} L_R g A_1^2
      + C R\|\nabla u\|_2^4 + C A_1\sqrt{\Xi B},
 \label{vh:Aprime} \\
 H' + \frac{3}{4} Z^2
 &\le \frac{1}{4} A_2^2 + C A_1^4 \log(e + A_1^2),
 \label{vh:Hprime} \\
 E_w' + c_1 D_w
 &\le C \Xi B + C L_R g(X + E_w),
 \label{vh:Wprime} \\
 E_\theta' + c_1 B
 &\le z D_w + C z^{-n} D_R^n g X(\log X + L_R).
 \label{vh:Kprime}
\end{align}
\end{lemma}

\begin{proof}
From $\eqref{main:pde}_2$, we deduce that
\begin{equation}\label{vx:flux-system}
\Delta G = \operatorname{div}(\rho\dot{u}), \quad \mu\Delta\omega = \operatorname{curl}(\rho\dot{u}).
\end{equation}
By elliptic regularity, \eqref{vac:density-weight}, and $h(\rho) \le C R^{b_h}$, we obtain for any $2 < p < \infty$,
\[
\|G\|_p + \|\omega\|_p \le C_p A_2, \quad \|\nabla G\|_2 + \|\nabla \omega\|_2 \le C R^{1/2} A_2, \quad \|G\|_2 \le C R^{b_h/2} A_1.
\]
Fixing $b > 1$ with $c/(b - 1) < a_0$, applying \eqref{vac:sob} to $\dot{u}$, and using H\"older's inequality, we also obtain for $k \ge 8$,
\begin{equation}\nonumber
 \|\nabla G\|_k
 \le C k \|\rho\dot{u}\|_k
 \le C k \|\rho W^{c/(bk)}\|_{\frac{bk}{b - 1}} \|\dot{u}\|_{bk,w}
 \le C k^{3/2} R(A_2 + \|\nabla\dot{u}\|_2).
\end{equation}

Set $\gamma = [4\log(e + R)]^{-1}$ and $p = 2(1 + \gamma)/\gamma$. By H\"older's inequality, the Gagliardo--Nirenberg inequality, $R^{b_h\gamma/2} \le C$, and $p^{1+\gamma} \le C\log(e + R)$,
\begin{equation}\nonumber
 \left\|\frac{G^2}{h}\right\|_2
 \le C\left\|\frac{G}{\sqrt{h}}\right\|_2^{1 - \gamma} \|G\|_p^{1 + \gamma}
 \le C A_1^{1 - \gamma} p^{(1 + \gamma)/2} \|G\|_2^\gamma \|\nabla G\|_2
 \le C\sqrt{L_R}\, A_1 A_2.
\end{equation}
Similarly,
\[
 \|\omega\|_4^4
 \le C\|\omega\|_2^2\|\nabla\omega\|_2^2
 \le CR A_1^2A_2^2.
\]
As in Lemma~\ref{vh:products}, we have
$\|P\|_1\le C$ and $\|P\|_\ell\le C\Theta$. Thus
\[
 \|P\|_p\le C\Theta^{a(1-1/p)},
 \qquad 1\le p\le\ell,\qquad a=\frac{\ell}{\ell-1}.
\]
Since $a<m$, $3a-b_\ell<2m$ and
$\Theta\le C(1+Z)^{1/m}$ by \eqref{vh:anchor},
\begin{equation}\nonumber
 \|P\|_4^4 \le C\Theta^{3a} \le C\Theta^{b_\ell}(1 + Z^2),
 \qquad
 \|P\|_2 \|P\|_4^2 \le C\Theta^{2a} \le z Z^2 + C z^{-n}.
\end{equation}
Using $h\operatorname{div} u = G + P$ and the div--curl estimate, we deduce
\begin{equation}\label{vh:quartic}
 Q_4 = \|\nabla u\|_4^4 + \|\sqrt{h}\,\operatorname{div} u\|_4^4
       + \|P\operatorname{div} u\|_2^2 + \|P\|_4^4
     \le C L_R g X,
 \quad
 \|\nabla u\|_2^2 \le C(A_1^2 + \|P\|_2^2) \le C g.
\end{equation}
We now prove \eqref{vh:Aprime} and \eqref{vh:Hprime}. Testing $\eqref{main:pde}_2$ by $2\dot{u}$ gives
\[
 2A_2^2=-2\int G\operatorname{div}\dot u\,dx
         -2\mu\int\omega\operatorname{curl}\dot u\,dx.
\]
Using
\begin{equation}\nonumber
\begin{aligned}
 \operatorname{div}\dot{u} = D_t(\operatorname{div} u) + (\operatorname{div} u)^2 - 2\det\nabla u,
 \quad \operatorname{curl}\dot{u} = \dot{\omega} + \omega\operatorname{div} u, \quad \dot{P} = \rho\dot{\theta} - P\operatorname{div} u,
\end{aligned}
\end{equation}
and integrating by parts, we arrive at
\begin{equation}\nonumber
\begin{aligned}
 (A_1^2)' + 2A_2^2
 &= -\mu\int \omega^2\operatorname{div} u\,dx
    + 4\int G\det\nabla u\,dx
    - \int \frac{(h + \rho\lambda')G^2}{h^2}\operatorname{div} u\,dx \\
 &\quad - 2\int \frac{\rho\lambda' P G}{h^2}\operatorname{div} u\,dx
    - 2 \int \frac{\rho}{h} G\dot{\theta}\,dx.
\end{aligned}
\end{equation}
By H\"older's inequality, the Jacobian--BMO estimate \cite{CLMS}, and $|\rho\lambda'| \le C h$, the first four integrals are bounded by
\begin{equation}\nonumber
\begin{aligned}
 & C\bigl(\|\omega\|_4^2 + \|\nabla G\|_2\|\nabla u\|_2\bigr)\|\nabla u\|_2
   + C\|\nabla u\|_2\left(\left\|\frac{G^2}{h}\right\|_2 + \|P\|_4^2\right) \\
 &\quad \le C\sqrt{L_R}\,A_1 A_2\|\nabla u\|_2
          + C R^{1/2} A_2\|\nabla u\|_2^2
          + C(A_1 + \|P\|_2)\|P\|_4^2 \\
 &\quad \le A_2^2 + z Z^2 + C z^{-n} L_R g A_1^2 + C R\|\nabla u\|_2^4,
\end{aligned}
\end{equation}
where we used \eqref{vh:quartic} and
$$A_1\|P\|_4^2 \le C A_1\Theta^{b_\ell/2}(1 + Z) \le z Z^2 + C z^{-1} g A_1^2.$$
For the last integral, H\"older's inequality and \eqref{vh:floor} yield
\begin{equation}\nonumber
 \left|\int \frac{\rho}{h} G\dot{\theta}\,dx\right|
 \le \left(\int \frac{G^2}{h}\,dx\right)^{1/2}
     \left(\int \rho\kappa|\dot{\theta}|^2\,dx\right)^{1/2}
     \left\|\frac{\rho}{h\kappa}\right\|_\infty^{1/2}
 \le A_1\sqrt{\Xi B},
\end{equation}
which proves \eqref{vh:Aprime}.

Multiplying the heat equation by $K(\theta)$ and integrating,
we obtain
\[
 H'+Z^2=\int K(Q-P\operatorname{div}u)\,dx.
\]
Applying \eqref{vh:k1} with $z=1/4$ yields
\[
 H'+\frac34Z^2
 \le\frac14A_2^2+CA_1^4\log(e+A_1^2),
\]
which is \eqref{vh:Hprime}.

We next prove \eqref{vh:Wprime}. Write $\mathcal{J}(f, g) = \det(\nabla f, \nabla g)$. By the Hardy-space Jacobian estimate \cite{CLMS} and elliptic regularity, we have
\begin{equation}\nonumber
 \|v\|_2 \le C\|\nabla u\|_2^2,
 \qquad
 \|\nabla v\|_2 \le C\|\nabla u\|_4^2.
\end{equation}
Differentiating $\operatorname{div} v = 2\det\nabla u$ along the flow yields
\begin{equation}\nonumber
 (\partial_t + \operatorname{div}(u\,\cdot))\det\nabla u
 = \mathcal{J}(\dot{u}^1, u^2) + \mathcal{J}(u^1, \dot{u}^2),
\end{equation}
and hence
\begin{equation}\nonumber
 \dot{v}_i
 = 2\partial_i\Delta^{-1}\bigl(\mathcal{J}(\dot{u}^1, u^2) + \mathcal{J}(u^1, \dot{u}^2)\bigr)
   + \sum_{j=1}^2 [u_j, \partial_i\partial_j\Delta^{-1}]\operatorname{div} v.
\end{equation}
Applying \cite{CLMS} together with the commutator estimate \cite{CRW} gives
\begin{equation}\nonumber
 \|\dot{v}\|_2
 \le C\|\nabla u\|_2\bigl(\|\nabla\dot{u}\|_2 + \|\nabla v\|_2\bigr)
 \le C\|\nabla u\|_2\bigl(\|\nabla w_1\|_2 + \|\nabla u\|_4^2\bigr).
\end{equation}

By the definition of $w_1 = \dot{u} + v$, we have
\begin{equation}\nonumber
 \operatorname{div} w_1 = D_t(\operatorname{div} u) + (\operatorname{div} u)^2,
 \qquad
 \operatorname{curl} w_1 = \dot{\omega} + \omega\operatorname{div} u,
\end{equation}
from which it follows that
\begin{equation}\nonumber
 \dot{G} + G\operatorname{div} u
 = h\operatorname{div} w_1 - \rho\lambda'(\operatorname{div} u)^2 - \rho\dot{\theta}.
\end{equation}
Introduce the stress-like tensor
\begin{equation}\nonumber
 \mathcal{T} = G I + \mu\omega \begin{pmatrix} 0 & -1 \\ 1 & 0 \end{pmatrix},
 \qquad \rho\dot{u} = \operatorname{div}\mathcal{T}.
\end{equation}
Applying $D_t + \operatorname{div} u$ to $\rho\dot{u} = \operatorname{div}\mathcal{T}$ yields
\begin{equation}\nonumber
 \rho D_t\dot{u}
 = \operatorname{div}\bigl(\dot{\mathcal{T}} + \mathcal{T}\operatorname{div} u - \mathcal{T}(\nabla u)^\top\bigr).
\end{equation}
Using $D_t w_1 = D_t\dot{u} + \dot{v}$ and testing this identity by $w_1$, we obtain
\begin{equation}\nonumber
 \frac{1}{2} E_w' + D_w
 = \int \bigl(\rho\lambda'(\operatorname{div} u)^2 + \rho\dot{\theta}\bigr)\operatorname{div} w_1\,dx
   + \int \mathcal{T}(\nabla u)^\top : \nabla w_1\,dx
   + \int \rho\dot{v}\cdot w_1\,dx.
\end{equation}

We estimate the three integrals on the right-hand side in turn. By $|\rho\lambda'| \le C h$, \eqref{vh:floor}, and H\"older's inequality, the first integral satisfies
\begin{equation}\nonumber
 \left|\int \bigl(\rho\lambda'(\operatorname{div} u)^2 + \rho\dot{\theta}\bigr)\operatorname{div} w_1\,dx\right|
 \le C\bigl(Q_4^{1/2} + \sqrt{\Xi B}\bigr) D_w^{1/2}.
\end{equation}
For the second integral, observe that
\begin{equation}\nonumber
 (\nabla u)^\top : \nabla w_1
 = \operatorname{div} u\operatorname{div} w_1
   - \mathcal{J}(u^1, w_1^2) - \mathcal{J}(w_1^1, u^2),
\end{equation}
and, since $G = h\operatorname{div} u - P$,
\begin{equation}\nonumber
 \int \frac{G^2}{h}(\operatorname{div} u)^2\,dx
 \le C\int \left(h(\operatorname{div} u)^4 + \frac{P^2}{h}(\operatorname{div} u)^2\right)dx
 \le C Q_4.
\end{equation}
Using \cite{CLMS} for the Jacobian terms and H\"older's inequality for the remaining terms, we deduce
\begin{equation}\nonumber
 \left|\int \mathcal{T}(\nabla u)^\top : \nabla w_1\,dx\right|
 \le C\bigl(Q_4^{1/2} + R^{1/2} A_2\|\nabla u\|_2 + \|\omega\|_4\|\nabla u\|_4\bigr) D_w^{1/2}
 \le C\bigl(Q_4^{1/2} + R^{1/2} A_2\|\nabla u\|_2\bigr) D_w^{1/2}.
\end{equation}
For the third integral, the bound on $\|\dot{v}\|_2$ gives
\begin{equation}\nonumber
 \left|\int \rho\dot{v}\cdot w_1\,dx\right|
 \le R^{1/2} E_w^{1/2} \|\dot{v}\|_2
 \le C R^{1/2} E_w^{1/2} \|\nabla u\|_2 \bigl(D_w^{1/2} + Q_4^{1/2}\bigr).
\end{equation}
Applying Young's inequality and \eqref{vh:quartic} yields \eqref{vh:Wprime}.

Finally, to prove \eqref{vh:Kprime}, we test the temperature equation by $\dot{K} = \kappa(\theta)\dot{\theta}$:
\begin{equation}\nonumber
 c_v B = \int \Delta K\,\dot{K}\,dx + \int Q\dot{K}\,dx - \int P\operatorname{div} u\,\dot{K}\,dx.
\end{equation}
Integrating by parts and applying $\|\nabla K\|_4^2 \le C Z\|\Delta K\|_2$, we obtain
\begin{equation}\nonumber
\begin{aligned}
 \int \Delta K\,\dot{K}\,dx
 &= -\frac{1}{2}(Z^2)' - \int \nabla K \cdot (\nabla u)\nabla K\,dx
    + \frac{1}{2}\int |\nabla K|^2\operatorname{div} u\,dx \\
 &\le -\frac{1}{2}(Z^2)' + C\|\nabla u\|_2 Z \|\Delta K\|_2.
\end{aligned}
\end{equation}
From the temperature equation and \eqref{vh:floor}, we have
\begin{equation}\nonumber
 \|\Delta K\|_2^2
 \le C\left\|\frac{\rho}{\kappa}\right\|_\infty B + C Q_4
 \le C R D_R^\alpha B + C Q_4,
\end{equation}
which together with \eqref{vh:quartic} gives
\begin{equation}\nonumber
 \int \Delta K\,\dot{K}\,dx
 \le -\frac{1}{2}(Z^2)' + \frac{c_v}{4} B + C D_R^\alpha L_R g X.
\end{equation}

For the heat production term, integration by parts yields
\begin{equation}\nonumber
 \int Q\dot{K}\,dx - \frac{d}{dt}\int K Q\,dx
 = -\int K(\dot{Q} + Q\operatorname{div} u)\,dx.
\end{equation}
Using $\dot{\lambda} = -\rho\lambda'\operatorname{div} u$, $D_t\nabla u = \nabla\dot{u} - (\nabla u)^2$, and $|\rho\lambda'| \le C h$, we estimate
\begin{equation}\nonumber
 |\dot{Q} + Q\operatorname{div} u|
 \le C\bigl(|\nabla u||\nabla\dot{u}| + h|\operatorname{div} u||D_t(\operatorname{div} u)|
            + |\nabla u|^3 + h|\operatorname{div} u|^3\bigr).
\end{equation}
Noting that $\int K h|\operatorname{div} u|^3\,dx \le \|K\sqrt{h}\,\operatorname{div} u\|_2 \|\sqrt{h}\,\operatorname{div} u\|_4^2$, Young's inequality gives
\begin{equation}\nonumber
\begin{aligned}
 \int Q\dot{K}\,dx - \frac{d}{dt}\int K Q\,dx
 &\le z\left(\|\nabla\dot{u}\|_2^2 + \int h|D_t(\operatorname{div} u)|^2\,dx\right) \\
 &\quad + C z^{-1}\left(Q_4 + \int K^2\bigl(|\nabla u|^2 + h(\operatorname{div} u)^2\bigr)\,dx\right).
\end{aligned}
\end{equation}
Since $\nabla\dot{u} = \nabla w_1 - \nabla v$ and $D_t(\operatorname{div} u) = \operatorname{div} w_1 - (\operatorname{div} u)^2$, we have
\begin{equation}\nonumber
 \|\nabla\dot{u}\|_2^2 + \int h|D_t(\operatorname{div} u)|^2\,dx
 \le C(D_w + Q_4).
\end{equation}
Combining the last two bounds with \eqref{vh:k2} and \eqref{vh:quartic}, and rescaling $z$, we deduce
\begin{equation}\nonumber
 \int Q\dot{K}\,dx - \frac{d}{dt}\int K Q\,dx
 \le z D_w + C z^{-n} g X(\log X + L_R).
\end{equation}

For the pressure work term, using $\theta^2\kappa(\theta) \le C(1 + K^2)$, $\rho/h \le C D_R$, \eqref{vh:k2}, and $D_u \le C g$, we obtain
\begin{equation}\nonumber
 \left|\int P\operatorname{div} u\,\dot{K}\,dx\right|
 \le \frac{c_v}{4} B + C D_R \int (1 + K^2)h(\operatorname{div} u)^2\,dx
 \le \frac{c_v}{4} B + C D_R g X \log X.
\end{equation}
Collecting the estimates for the three terms establishes \eqref{vh:Kprime} and completes the proof of Lemma~\ref{vh:diff}.
\end{proof}

\begin{lemma}\label{vh:joint}
Set
\begin{equation}\nonumber
 \mathcal{H} = 1 + H,
 \qquad N_R = R + \log(e + A_1^2),
 \qquad F_R = A_1^4 N_R + R\mathcal{H}^{2s}.
\end{equation}
Then there exist $t_R = c_2 D_R^{-(1 + \alpha)}$ and $1 \le C_R, k_R \le M_R = C D_R^C$ such that, with
\begin{equation}\nonumber
\begin{gathered}
 Y_R = e + E_\theta + t_R E_w + C_R F_R,
 \quad V_R = k_R\log Y_R + \log S,\quad
 \mathcal{D}_R = \frac{D_w + B}{Y_R} + \frac{Y_R + A_2^2 + Z^2}{S},
\end{gathered}
\end{equation}
one has
\begin{equation}\label{vh:jointest}
 \sup_{0 \le t \le T} V_R(t) + \int_0^T \mathcal{D}_R\,dt
 \le Q_R = C L_R \exp(C D_R^C).
\end{equation}
\end{lemma}

\begin{proof}
Choose $c_2 > 0$ sufficiently small in $t_R = c_2 D_R^{-(1+\alpha)}$ so that $C t_R \Xi \le c_1/4$ by \eqref{vh:floor}. Adding \eqref{vh:Kprime} with $z = c_1 t_R/4$ to $t_R$ times \eqref{vh:Wprime} yields
\begin{equation}\nonumber
 (E_\theta + t_R E_w)' + c(B + t_R D_w)
 \le C D_R^C g X(\log X + L_R) + C t_R L_R g E_w.
\end{equation}

We next verify that $Y_R$ dominates the higher-order quantities. By \eqref{vh:k1}, for any $0 < \eta \le 1$,
\begin{equation}\nonumber
 \int K Q\,dx
 \le \eta(A_2^2 + Z^2) + C\eta^{-n} A_1^4\log(e + A_1^2).
\end{equation}
Moreover, the $L^2$ bound for $v$ in Lemma~\ref{vh:diff} and \eqref{vh:compare} give
\begin{equation}\nonumber
 \int \rho|v|^2\,dx
 \le C R\|\nabla u\|_2^4
 \le C R(A_1^4 + \mathcal{H}^{2s})
 \le C F_R,
\end{equation}
so $w_1 = \dot{u} + v$ implies
$$\frac{1}{2} A_2^2 - C F_R \le E_w \le 2 A_2^2 + C F_R.$$
Choosing $\eta \le \min\{1/4, t_R/4\}$ and then $1 \le C_R \le C D_R^C$ sufficiently large, we obtain $Y_R \ge e$ and
\begin{equation}\label{vh:Y-equivalence}
 Y_R \asymp e + t_R A_2^2 + Z^2 + C_R F_R.
\end{equation}
In particular,
\begin{equation}\label{vh:Y-control}
 X \le C D_R^C Y_R,
 \qquad t_R E_w + C_R F_R \le C Y_R,
 \qquad N_R \le C(R + \log Y_R).
\end{equation}
Dividing the differential inequality for $E_\theta + t_R E_w$ by $Y_R$ and applying \eqref{vh:Y-control} gives
\begin{equation}\label{vh:joint-preliminary}
 \frac{(E_\theta + t_R E_w)'}{Y_R} + c\frac{B + t_R D_w}{Y_R}
 \le C D_R^C g(L_R + \log Y_R).
\end{equation}

To include the derivative of $F_R$, write $a = A_1^2$. Since $(a^2 N_R)' = \bigl(2a N_R + a^2/(e + a)\bigr)a'$ and $2a N_R + a^2/(e + a) \le C a N_R$, \eqref{vh:Aprime} with $z = 1$ yields
\begin{equation}\nonumber
 \frac{C_R(a^2 N_R)'}{Y_R}
 \le \frac{C C_R a N_R}{Y_R}
     \bigl(Z^2 + L_R g a + R\|\nabla u\|_2^4 + A_1\sqrt{\Xi B}\bigr).
\end{equation}
By \eqref{vh:Y-equivalence}--\eqref{vh:Y-control}, we have $C_R a^2 N_R \le C Y_R$ and $Z^2 + R\|\nabla u\|_2^4 \le C Y_R$, while Young's inequality gives
\begin{equation}\nonumber
 \frac{C C_R A_1^3 N_R\sqrt{\Xi B}}{Y_R}
 \le \varepsilon\frac{B}{Y_R} + C_\varepsilon\frac{\Xi C_R^2 A_1^6 N_R^2}{Y_R}
 \le \varepsilon\frac{B}{Y_R} + C_\varepsilon\Xi C_R A_1^2 N_R
 \le \varepsilon\frac{B}{Y_R} + C_\varepsilon D_R^C g(R + \log Y_R).
\end{equation}
Hence,
\begin{equation}\label{vh:F-density}
 \frac{C_R(A_1^4 N_R)'}{Y_R}
 \le \varepsilon\frac{B}{Y_R} + C_\varepsilon D_R^C g(L_R + \log Y_R).
\end{equation}
Similarly, by \eqref{vh:Hprime} and \eqref{vh:Y-equivalence}--\eqref{vh:Y-control},
\begin{equation}\nonumber
 \frac{C_R R(\mathcal{H}^{2s})'}{Y_R}
 = \frac{2s C_R R \mathcal{H}^{2s-1}}{Y_R} H'
 \le \frac{C C_R R \mathcal{H}^{2s-1}}{Y_R}\bigl(A_2^2 + A_1^4\log(e + A_1^2)\bigr)
 \le C D_R^C R \mathcal{H}^{(2s-1)_+}.
\end{equation}
Adding the last two bounds to \eqref{vh:joint-preliminary} and choosing $\varepsilon > 0$ sufficiently small yields
\begin{equation}\label{vh:logY}
 (\log Y_R)' + c\frac{B + t_R D_w}{Y_R}
 \le C D_R^C g_*(L_R + \log Y_R).
\end{equation}

Next, adding \eqref{vh:Aprime} with $z = 1/4$ to \eqref{vh:Hprime} gives
\begin{equation}\nonumber
 S' + \frac{3}{4} A_2^2 + \frac{1}{2} Z^2
 \le C L_R g A_1^2 + C R\|\nabla u\|_2^4
     + C A_1^4\log(e + A_1^2) + C A_1\sqrt{\Xi B}.
\end{equation}
Since $A_1^2, \mathcal{H} \le S$, it follows from \eqref{vh:compare} that
\begin{equation}\nonumber
 \frac{R\|\nabla u\|_2^4}{S} \le C R\bigl(A_1^2 + \mathcal{H}^{(2s-1)_+}\bigr),
 \qquad
 \frac{A_1^4\log(e + A_1^2)}{S} \le C A_1^2\log Y_R,
\end{equation}
and therefore
\begin{equation}\label{vh:logS-preliminary}
 (\log S)' + c\frac{A_2^2 + Z^2}{S}
 \le C g_*(L_R + \log Y_R) + C\frac{A_1\sqrt{\Xi B}}{S}.
\end{equation}
By \eqref{vh:Y-equivalence}--\eqref{vh:Y-control} and $A_1^2, \mathcal{H} \le S$, we have
\begin{equation}\nonumber
 \frac{Y_R}{S}
 \le C D_R^C\left(\frac{A_2^2 + Z^2}{S} + 1 + A_1^2 N_R + R\mathcal{H}^{(2s-1)_+}\right)
 \le C D_R^C\left(\frac{A_2^2 + Z^2}{S} + g_*(R + \log Y_R)\right),
\end{equation}
while Young's inequality and \eqref{vh:floor} yield
\begin{equation}\nonumber
 C\frac{A_1\sqrt{\Xi B}}{S}
 \le \varepsilon\frac{Y_R}{S} + C_\varepsilon\frac{A_1^2\Xi B}{S Y_R}
 \le \varepsilon\frac{Y_R}{S} + C_\varepsilon D_R^{1+\alpha}\frac{B}{Y_R}.
\end{equation}
Inserting these two bounds into \eqref{vh:logS-preliminary} and taking $\varepsilon = c D_R^{-C}$ sufficiently small gives
\begin{equation}\label{vh:logS}
 (\log S)' + c D_R^{-C}\frac{Y_R + A_2^2 + Z^2}{S}
 \le C D_R^C\frac{B}{Y_R} + C D_R^C g_*(L_R + \log Y_R).
\end{equation}
Multiplying \eqref{vh:logY} by a sufficiently large $k_R$ and adding \eqref{vh:logS} and yields
\begin{equation}\label{vh:joint-differential}
 V_R' + M_R^{-1}\mathcal{D}_R \le M_R g_*(L_R + V_R),
 \qquad M_R = C D_R^C,
\end{equation}
where we used $t_R = c_2 D_R^{-(1+\alpha)}$ and $\log Y_R \le V_R$.

Finally, \eqref{vh:Y-equivalence} and the initial bounds imply
$$Y_R(0) \le C D_R^C(1 + R), \quad V_R(0) + L_R \le C D_R^C L_R.$$
Since $\int_0^T g_*\,dt \le C$ by \eqref{vh:tails}, Gr\"onwall's inequality applied to \eqref{vh:joint-differential} yields
\begin{equation}\nonumber
 V_R(t) + L_R
 \le \bigl(V_R(0) + L_R\bigr)\exp\left(M_R \int_0^t g_*\,ds\right)
 \le C L_R \exp(C D_R^C).
\end{equation}
Integrating \eqref{vh:joint-differential} over $[0,T]$ then gives
\begin{equation}\nonumber
 \int_0^T \mathcal{D}_R\,dt
 \le M_R V_R(0) + M_R^2 \int_0^T g_*(L_R + V_R)\,dt
 \le C L_R \exp(C D_R^C),
\end{equation}
which establishes \eqref{vh:jointest} and completes the proof of Lemma~\ref{vh:joint}.
\end{proof}

Take $r_0 > 2\|\rho_0\|_\infty$ sufficiently large and set
\begin{equation}\nonumber
 r_j = 2^j r_0,
 \qquad R_j = 2 + r_j,
 \qquad L_j = R_j\log(e + R_j).
\end{equation}
For each attained level $r_j$, let $t_j$ be the first time at which $\|\rho(t_j)\|_\infty = r_j$. Fix $\zeta_* > 0$ small enough so that $D(r_j) \le (\log r_j)^{\zeta_*}$ implies $Q_{R_j} = R_j^{1 + o(1)}$ as $j \to \infty$ in \eqref{vh:jointest}. For each such attained index $j$, let $j^-$ denote the preceding index with the same property, setting $t_{j^-} = 0$ if none exists. With $\eta = 1/8$, define
\begin{equation}\nonumber
\begin{gathered}
 I_j^0 = [t_{j^-}, t_j],
 \quad E_j = I_j^0 \cap \{S_0 > R_j^\eta\},\quad
 d_j = \frac{1}{L_j} \int_{E_j} \frac{A_2^2}{S_0}\,dt,
 \quad e_j = \int_{E_j} S_0\,dt.
\end{gathered}
\end{equation}

\begin{lemma}\label{vh:sum}
For every attained level $r_N$, one has
\begin{equation}\label{vh:sumest}
 \sum_{\substack{j \le N \\ D(r_j) \le (\log r_j)^{\zeta_*}}} (d_j + e_j) \le C,
\end{equation}
where $C$ is independent of $N$.
\end{lemma}

\begin{proof}
Testing the momentum equation by $\dot{u}$ and using $\operatorname{div}\dot{u} = \operatorname{div} w_1 - 2\det\nabla u$ and $\operatorname{curl}\dot{u} = \dot{\omega} + \omega\operatorname{div} u$, we obtain
\begin{equation}\nonumber
 \frac{1}{2} D_u' + A_2^2
 = \int P\operatorname{div} w_1\,dx
   + 2\int G\det\nabla u\,dx
   - \frac{1}{2}\int (h + \rho\lambda')(\operatorname{div} u)^3\,dx
   - \frac{\mu}{2}\int \omega^2\operatorname{div} u\,dx.
\end{equation}
By the Jacobian--BMO estimate \cite{CLMS}, Lemma~\ref{vh:diff}, and \eqref{vh:compare}, we have
\begin{align*}
 \left|\int G\det\nabla u\,dx\right|
 + \left|\int \omega^2\operatorname{div} u\,dx\right|
& \le C\bigl(\|\nabla G\|_2\|\nabla u\|_2^2 + \|\omega\|_2\|\nabla\omega\|_2\|\nabla u\|_2\bigr) \\
& \le C R^{1/2} A_2\|\nabla u\|_2^2
 \le \varepsilon A_2^2 + C_\varepsilon R S_0^2.
\end{align*}
For the cubic divergence term, $|\rho\lambda'| \le C h$ and $h\operatorname{div} u = G + P$ yield
\begin{equation}\nonumber
 \int |h + \rho\lambda'||\operatorname{div} u|^3\,dx
 \le C\int \frac{|G|^3 + P^3}{h^2}\,dx.
\end{equation}
Using the bound on $\|G^2/h\|_2$ from Lemma~\ref{vh:diff} together with \eqref{vh:p23} and \eqref{vh:compare}, we find
\begin{equation}\nonumber
 \int \frac{|G|^3}{h^2}\,dx
 \le C A_1\left\|\frac{G^2}{h}\right\|_2
 \le C\sqrt{L_R}\, A_1^2 A_2
 \le \varepsilon A_2^2 + C_\varepsilon L_R S_0^2,
 \quad
 \int \frac{P^3}{h^2}\,dx
 \le C S_0(1 + F_\alpha).
\end{equation}
Choosing $\varepsilon > 0$ sufficiently small gives
\begin{equation}\nonumber
 D_u' + A_2^2
 \le 2\left|\int P\operatorname{div} w_1\,dx\right|
     + C L_R S_0^2 + C S_0(1 + F_\alpha).
\end{equation}

Multiplying \eqref{vh:Hprime} by $\sigma(1 + H)^{\sigma - 1}$ and adding the result to the preceding inequality yields
\begin{equation}\nonumber
 S_0' + c_3 A_2^2
 \le 2\left|\int P\operatorname{div} w_1\,dx\right|
     + C L_R S_0^2 + C S_0(1 + F_\alpha) + C A_1^4\log(e + A_1^2).
\end{equation}
Dividing by $S_0$ and using $A_1^2 \le C S_0$ from \eqref{vh:compare}, we arrive at
\begin{equation}\nonumber
 (\log S_0)' + c_3\frac{A_2^2}{S_0}
 \le \frac{2}{S_0}\left|\int P\operatorname{div} w_1\,dx\right|
     + C b(t)(L_R + \log S_0),
\end{equation}
where $b(t) = 1 + S_0 + F_\alpha$ satisfies $\int_0^T b(t)\,dt \le C$ by \eqref{vh:tails}.

We next estimate the pressure term on $E_j$. By H\"older's inequality, \eqref{vh:p23}, and \eqref{vh:compare}, with $\delta = 1 - (1 + s)/(2\sigma) > 0$,
\begin{equation}\nonumber
 \frac{1}{S_0}\left|\int P\operatorname{div} w_1\,dx\right|
 \le \frac{C(1 + H)^{s/2} S^{1/2}}{S_0}
     \left(\frac{D_w}{Y_R}\frac{Y_R}{S}\right)^{1/2}
 \le C S_0^{\frac{1 + s}{2\sigma} - 1}
     \left(\frac{D_w}{Y_R} + \frac{Y_R}{S}\right)
 \le C S_0^{-\delta}\mathcal{D}_R.
\end{equation}
Since $\rho \le r_j$ on $[0, t_j]$ and $S_0 > R_j^\eta$ on $E_j$, applying \eqref{vh:jointest} gives
\begin{equation}\nonumber
 \frac{1}{L_j}\int_{E_j} \frac{\left|\int P\operatorname{div} w_1\,dx\right|}{S_0}\,dt
 \le \frac{C R_j^{-\eta\delta}}{L_j} \int_0^{t_j} \mathcal{D}_{R_j}\,dt
 \le C R_j^{-\eta\delta}\frac{Q_{R_j}}{L_j}
 \le C R_j^{-\eta\delta/2}.
\end{equation}

To sum across density levels, set
\begin{equation}\nonumber
 U_j = (\log S_0 - \eta\log R_j)_+,
 \qquad
 b_*(t) = \exp\left(-C\int_0^t b(s)\,ds\right).
\end{equation}
On $I_j^0$, we have $U_j' = \mathbf{1}_{E_j}(\log S_0)'$ almost everywhere, while on $E_j$ it holds that $\log S_0 = U_j + \eta\log R_j$ and $L_j + \eta\log R_j \le C L_j$. Choosing the constant in $b_*$ sufficiently large and using $0 < c \le b_* \le 1$, we obtain
\begin{equation}\nonumber
 \left(\frac{b_* U_j}{L_j}\right)'
 + \frac{c_3 b_*}{L_j}\mathbf{1}_{E_j}\frac{A_2^2}{S_0}
 \le C\mathbf{1}_{E_j} b
   + \frac{C}{L_j}\mathbf{1}_{E_j}\frac{\left|\int P\operatorname{div} w_1\,dx\right|}{S_0}.
\end{equation}
Integrating over $I_j^0$ and inserting the pressure bound yields
\begin{equation}\label{vh:tel}
 c d_j + \left[\frac{b_* U_j}{L_j}\right]_{t_{j^-}}^{t_j}
 \le C\int_{E_j} b(t)\,dt + C R_j^{-\eta\delta/2}.
\end{equation}

For any two consecutive indices $j < j^+$ satisfying $D(r_j) \le (\log r_j)^{\zeta_*}$, monotonicity of $R_j$ and $L_j$ implies
\begin{equation}\nonumber
 \frac{U_{j^+}(t_j)}{L_{j^+}} \le \frac{U_j(t_j)}{L_j},
\end{equation}
so the interior boundary contributions in the sum of \eqref{vh:tel} are nonnegative, while the initial boundary term is bounded by $(\log S_0(0))_+ \le C$. Summing \eqref{vh:tel} over all such $j \le N$ and using the disjointness of the interiors of $I_j^0$ together with \eqref{vh:tails}, we conclude
\begin{equation}\nonumber
 \sum_{\substack{j \le N \\ D(r_j) \le (\log r_j)^{\zeta_*}}} d_j
 \le C\left(1 + \int_0^T b(t)\,dt + \sum_{j \ge 0} R_j^{-\eta\delta/2}\right)
 \le C,
 \quad
 \sum_{\substack{j \le N \\ D(r_j) \le (\log r_j)^{\zeta_*}}} e_j
 \le \int_0^T S_0\,dt
 \le C,
\end{equation}
which proves \eqref{vh:sumest} and completes the proof of Lemma~\ref{vh:sum}.
\end{proof}

\begin{proposition}\label{vh:cap}
Under any of the conditions \textup{(A)}--\textup{(D)} in Definition~\ref{main:classes}, there exists a positive constant $C$ such that
\begin{equation}\label{v2:cap-input}
\sup_{0 \le t \le T} \left( \| \rho \|_{\infty} + \| \theta^{-1} \|_{\infty} + A_1 + A_2 + Z + H \right)
+ \int_0^T \bigl(B + D_w + \|\nabla\dot{u}\|_2^2\bigr)\,dt \le C.
\end{equation}
\end{proposition}

\begin{proof}
We establish uniform flux bounds on density sublevels, exclude arbitrarily large density levels along particle trajectories, and then close the remaining estimates.

\smallskip\noindent\textit{Step 1. Estimates for $G$ on a density sublevel.}
Set $T_2 = 1 + A_2$ and $T_3 = 1 + A_2 + \|\nabla\dot{u}\|_2$. On the sublevel $\{\rho \le r_j\}$, the elliptic bounds in Lemma~\ref{vh:diff} give
\begin{equation}\label{vh:cap-elliptic}
\begin{aligned}
 \|G\|_p &\le C_p T_2 \quad (p > 2),
 &\qquad \|G\|_2 &\le C R_j^{b_h/2} A_1, \\
 \|\nabla G\|_2 &\le C R_j^{1/2} T_2,
 &\qquad \|\nabla G\|_n &\le C n^{3/2} R_j T_3 \quad (n \ge 8).
\end{aligned}
\end{equation}
For any $p, n \ge 8$, balancing the average of $G$ on a ball of radius $l$ against Morrey's oscillation bound yields, with a uniform constant,
\begin{equation}\label{vh:cap-morrey}
 \|G\|_\infty \le C \|G\|_p^{1 - \beta} \|\nabla G\|_n^\beta,
 \qquad
 \beta = \frac{2}{p(1 - 2/n) + 2}.
\end{equation}

Let $\mathcal{I} = \{j \ge 1 : D(r_j) \le (\log r_j)^{\zeta_*}\}$. For every attained $j \in \mathcal{I}$, \eqref{vh:jointest} gives $Q_{R_j} = R_j^{1 + o(1)}$. Taking $\zeta = 1/Q_{R_j}$ in \eqref{vh:tails} and using $\log S \le Q_{R_j}$ yield $\int_0^{t_j} S\,dt \le C Q_{R_j}^2$. Together with $\|\nabla\dot{u}\|_2^2 \le C(D_w + L_R g X)$, $X/Y_R \le C D_R^C$, and \eqref{vh:jointest}, we obtain
\begin{equation}\nonumber
 \int_0^{t_j} \frac{T_3^2}{Y_{R_j}}\,dt
 + \int_0^{t_j} \frac{Y_{R_j}}{S}\,dt
 + \int_0^{t_j} S\,dt
 \le C R_j^C.
\end{equation}
Applying H\"older's inequality to $(T_3^2/Y_{R_j})^{1/3}(Y_{R_j}/S)^{1/3} S^{1/3}$ gives, for some fixed $c_4 > 0$,
\begin{equation}\label{vh:cap-T3}
 \int_0^{t_j} T_3^{2/3}\,dt \le C R_j^{c_4}.
\end{equation}

On $\mathcal{G}_j = [0, t_j] \cap \{S_0 \le R_j^\eta\}$ ($\eta = 1/8$), \eqref{vh:compare} yields $S \le C R_j^{\eta/\sigma}$, and \eqref{vh:jointest} implies
\begin{equation}\nonumber
 \int_{\mathcal{G}_j} T_2^2\,dt \le C R_j^{\eta/\sigma} Q_{R_j}.
\end{equation}
Taking $n = 8$ and a sufficiently large fixed $p$ in \eqref{vh:cap-morrey}, and applying \eqref{vh:cap-elliptic}, \eqref{vh:cap-T3}, and H\"older's inequality, we find a fixed $\xi \in (0, 1)$ such that
\begin{equation}\label{vh:good}
 \int_{\mathcal{G}_j} \|G\|_\infty\,dt
 \le C R_j^\beta
     \left(\int_{\mathcal{G}_j} T_2^2\,dt\right)^{(1-\beta)/2}
     \left(\int_0^{t_j} T_3^{2/3}\,dt\right)^{3\beta/2}
 \le C R_j^\xi.
\end{equation}

On $E_j = [t_{j^-}, t_j] \cap \{S_0 > R_j^\eta\}$, set $p_j = 16 + 10\log(e + R_j)$, $n_j = e^{p_j}$, and
\begin{equation}\nonumber
 \beta_j = \frac{2}{p_j(1 - 2/n_j) + 2},
 \qquad a_j = \left(1 - \frac{2}{n_j}\right)\beta_j,
 \qquad w_j = 1 - \left(2 - \frac{1}{n_j}\right)\beta_j.
\end{equation}
Substituting $\|G\|_{p_j} \le C\sqrt{p_j}\,\|G\|_2^{2/p_j}\|\nabla G\|_2^{1 - 2/p_j}$ and \eqref{vh:cap-elliptic} into \eqref{vh:cap-morrey} and noting that
$$(n_j^{3/2} R_j^{1/2})^{\beta_j} R_j^{(b_h - 1)a_j/2} \le C$$
and $T_2^{\beta_j/n_j} \le C$, we deduce
\begin{equation}\label{vh:cap-linfty}
 \|G\|_\infty \le C\sqrt{p_j R_j}\, A_1^{a_j} T_2^{w_j} T_3^{\beta_j},
\end{equation}
where $\frac{a_j}{2} + w_j + \frac{3\beta_j}{2} = 1$ and $1 - w_j = O(j^{-1})$.
Define
\begin{equation}\nonumber
 \widehat{d}_j = \frac{1}{L_j}\int_{E_j} \frac{T_2^2}{S_0}\,dt
 \le 2\left(d_j + \frac{|E_j|}{L_j}\right),
 \qquad x_j = (e_j \widehat{d}_j)^{1/2}.
\end{equation}
By Lemma~\ref{vh:sum} and the disjointness of the interiors of $[t_{j^-}, t_j]$, we have $\sum_{j \in \mathcal{I}} (\widehat{d}_j + e_j) \le C$, whence $\sum_{j \in \mathcal{I}} x_j \le C$ and $\int_{E_j} T_2\,dt \le \sqrt{L_j}\, x_j$. Applying H\"older's inequality to \eqref{vh:cap-linfty} together with $\int_0^T A_1^2\,dt \le C$, \eqref{vh:cap-T3}, $p_j \asymp \log R_j$, and $\beta_j = O(j^{-1})$ yields
\begin{equation}\label{vh:bad}
 \int_{E_j} \|G\|_\infty\,dt
 \le C\sqrt{p_j R_j}
     \left(\int_{E_j} A_1^2\,dt\right)^{a_j/2}
     \left(\int_{E_j} T_2\,dt\right)^{w_j}
     \left(\int_{E_j} T_3^{2/3}\,dt\right)^{3\beta_j/2}
 \le C L_j x_j^{w_j}.
\end{equation}
Splitting the sum at $x_j = e^{-j}$ and using $1 - w_j = O(j^{-1})$ give
\begin{equation}\nonumber
 \sum_{j \in \mathcal{I}} x_j^{w_j}
 \le C\sum_{j \in \mathcal{I}} x_j + C\sum_{j \ge 1} e^{-j/2}
 \le C.
\end{equation}

\smallskip\noindent\textit{Step 2. Exclusion of arbitrarily large density levels.}
For any attained $j \in \mathcal{I}$, choose a particle trajectory reaching the maximum $\|\rho(t_j)\|_\infty = r_j$ and follow it backwards to its last visit to the level $r_j/2$. Since $r_{j^-} \le r_j/2$, this crossing interval lies in $[t_{j^-}, t_j]$. Integrating $D_t\Lambda(\rho) = -G - P \le -G$ along the trajectory and applying \eqref{vh:good}--\eqref{vh:bad} yield
\begin{equation}\label{vh:cap-increment}
 \frac{\Lambda(r_j) - \Lambda(r_j/2)}{L_j}
 \le \frac{1}{L_j}\int_{t_{j^-}}^{t_j} \|G(t)\|_\infty\,dt
 \le C\frac{R_j^\xi}{L_j} + C x_j^{w_j}.
\end{equation}
Summing \eqref{vh:cap-increment} over all attained $j \le N$ in $\mathcal{I}$, we obtain
\begin{equation}\label{vh:cap-sum}
 \sum_{\substack{j \le N \\ j \in \mathcal{I}}}
 \frac{\Lambda(r_j) - \Lambda(r_j/2)}{L_j} \le C,
\end{equation}
uniformly in $N$.

In cases \textup{(A)}, \textup{(C)}, and \textup{(D)}, \eqref{vac:struct} and monotonicity of $D$ give
\begin{equation}\nonumber
 \frac{\Lambda(r_j) - \Lambda(r_j/2)}{L_j} \ge \frac{c}{j D(r_j)},
 \qquad
 \sum_{j \ge 1} \frac{1}{j D(r_j)} = \infty,
\end{equation}
while the excluded indices satisfy $\sum_{j \notin \mathcal{I}} [j D(r_j)]^{-1} \le C\sum_{j \ge 1} j^{-1-\zeta_*} < \infty$. Thus the subseries over $j \in \mathcal{I}$ diverges. In case \textup{(B)}, $D_h(r) = L(r)^{o(1)}$, so every sufficiently large index belongs to $\mathcal{I}$ and
\begin{equation}\nonumber
 \frac{\Lambda(r_j) - \Lambda(r_j/2)}{L_j}
 \asymp \int_{r_j/2}^{r_j} \frac{h(z)}{z^2 L(z)}\,dz.
\end{equation}
By Tonelli's theorem and \textup{(B)},
\begin{equation}\nonumber
 \infty = \int_{r_0}^\infty \frac{\Lambda(r)}{r^2 L(r)}\,dr
 = \Lambda(r_0)\int_{r_0}^\infty \frac{dr}{r^2 L(r)}
   + \int_{r_0}^\infty \frac{h(z)}{z}\left(\int_z^\infty \frac{dr}{r^2 L(r)}\right)dz,
\end{equation}
and since $\int_z^\infty \frac{dr}{r^2 L(r)} \asymp \frac{1}{z L(z)}$, it follows that $\int_{r_0}^\infty \frac{h(z)}{z^2 L(z)}\,dz = \infty$. In all four cases, the left-hand side of \eqref{vh:cap-sum} diverges as $N \to \infty$, so there exists a fixed $N$ depending only on the initial data and $T$ such that
\begin{equation}\nonumber
 \sup_{0 \le t \le T} \|\rho(t)\|_\infty \le r_N.
\end{equation}

\smallskip\noindent\textit{Step 3. Completion of the estimates.}
With $R = 2 + r_N$ fixed, \eqref{vh:floor} gives $\|1/\theta\|_{L^\infty(\mathbb{R}^2 \times [0,T])} \le C$, and Lemma~\ref{vh:joint} yields
\begin{equation}\nonumber
 \sup_{0 \le t \le T} V_R(t) + \int_0^T \mathcal{D}_R\,dt \le Q_R \le C.
\end{equation}
Since $Y_R \ge e$ and $S \ge 1$, this directly bounds $A_1, A_2, Z, H$ in $L^\infty(0,T)$ and $B, D_w$ in $L^1(0,T)$. Finally, by Lemma~\ref{vh:diff} and \eqref{vh:tails},
\begin{equation}\nonumber
 \int_0^T \|\nabla\dot{u}\|_2^2\,dt
 \le C\int_0^T (D_w + L_R g X)\,dt
 \le C.
\end{equation}
This completes the proof of Proposition~\ref{vh:cap}.
\end{proof}

\subsection{Higher-order estimates}\label{sec:high}

The goal of this subsection is to establish some necessary higher-order estimates, which ensure that the local strong solution can be extended globally in time.

\begin{lemma}\label{vx:ell}
Under the assumptions of Proposition~\ref{vh:cap}, for each $2\le p<\infty$ and almost every $t\in[0,T]$, we have
\begin{align}
 \|u\|_\infty+\|\nabla u\|_p+\|P\|_p \le C_{p},\label{vx:vel}
\end{align}
and
\begin{align}
 \|\theta\|_\infty^{1+\alpha}+\|\nabla^2K(\theta)\|_2
\le C \bigl(1+\|\sqrt{\rho\kappa(\theta)}\,\dot\theta\|_2\bigr).\label{vx:heat}
\end{align}
\end{lemma}

\begin{proof}
By \eqref{vh:anchor} and Proposition~\ref{vh:cap}, we have
\begin{equation}\nonumber
 \|K(\theta)\|_{s,w} \le C\sqrt{s}\,(1 + Z) \le C\sqrt{s},
 \qquad s \ge 2.
\end{equation}
For any fixed $p \ge 2$, choose $b > 1$ large enough and then $\delta > 0$ small enough so that, with $b' = b/(b - 1)$,
\begin{equation}\nonumber
 \frac{pb}{m} \ge 2,
 \qquad
 \delta p b' + \frac{c}{b - 1} \le a_0,
 \qquad
 4\delta \le a_0.
\end{equation}
Applying H\"older's inequality together with $w = W^{-c}$ and \eqref{vac:entropy} yields
\begin{equation}\nonumber
\begin{aligned}
 \|W^\delta P\|_p^p
 &= \int (\theta^{pb} w)^{1/b} \rho^p W^{p\delta + c/b}\,dx \\
 &\le \|\theta\|_{pb,w}^p
      \left(\int \rho^{pb'} W^{\delta p b' + c/(b - 1)}\,dx\right)^{1/b'} \\
 &\le \left(\frac{m}{\kappa_0}\|K(\theta)\|_{pb/m,w}\right)^{p/m}
      \|\rho\|_\infty^{p - 1/b'} \left(\int \rho W^{a_0}\,dx\right)^{1/b'}
  \le C_p,
\end{aligned}
\end{equation}
and hence
\begin{equation}\label{vx:weighted-pressure}
 \|W^\delta P\|_p \le C_p.
\end{equation}
Next, using \eqref{vac:entropy} and H\"older's inequality, we arrive at
\begin{equation}\label{vx:force-weight}
 \|W^{2\delta}\rho\dot{u}\|_1 + \|\rho\dot{u}\|_2
 \le \bigl(\|W^{2\delta}\sqrt{\rho}\|_2 + \|\sqrt{\rho}\|_\infty\bigr)
     \|\sqrt{\rho}\,\dot{u}\|_2
 \le C.
\end{equation}
By the Gagliardo--Nirenberg inequality, \eqref{v2:cap-input}, \eqref{vx:force-weight} and \eqref{vx:flux-system}, we have
\[
 \|(G,\omega)\|_p
 \le C_p\|(G,\omega)\|_2^{2/p}
       \bigl(\|\nabla(G,\omega)\|_2+\|(G,\omega)\|_2\bigr)^{1-2/p}\le C_p .
\]
Substituting $d_u=(G+P)/h(\rho)$ into the div--curl estimate and using \eqref{vx:weighted-pressure}, we obtain
\[
 \|\nabla u\|_p\le C_p\bigl(\|(G,\omega)\|_p+\|P\|_p\bigr)\le C_p .
\]
Moreover, \eqref{vh:M-high} in Lemma~\ref{vh:products} yields $\|W^\delta(G, \omega)\|_2 \le C A_2 \le C$. Since $W^{2\delta} \in A_2$, the weighted div--curl estimate and \eqref{vx:weighted-pressure} imply
\begin{equation}\label{vx:weighted-gradient}
 \|W^\delta\nabla u\|_2
 \le C\bigl(\|W^\delta\operatorname{div} u\|_2 + \|W^\delta\omega\|_2\bigr)
 \le C\bigl(\|W^\delta(G, \omega)\|_2 + \|W^\delta P\|_2\bigr)
 \le C.
\end{equation}

To estimate $\|u\|_\infty$, rewrite the momentum equation as
\begin{equation}\nonumber
 \mu\Delta u + (\mu + \lambda(0))\nabla\operatorname{div} u
 = \rho\dot{u} + \nabla h_1,
 \qquad
 h_1 = P - \bigl(\lambda(\rho) - \lambda(0)\bigr)\operatorname{div} u.
\end{equation}
Using $|\lambda(\rho) - \lambda(0)| \le C \rho$ and the preceding bounds, we obtain
\begin{equation}\nonumber
 \|h_1\|_1 + \|h_1\|_4
 \le \|P\|_1 + \|P\|_4
     + C\bigl(\|\rho\|_2\|\operatorname{div} u\|_2
              + \|\rho\|_\infty\|\operatorname{div} u\|_4\bigr)
 \le C.
\end{equation}
Let $\mathcal{N}$ denote the inverse of the constant-coefficient Lam\'e operator $\mu\Delta + (\mu + \lambda(0))\nabla\operatorname{div}$. Its fundamental solution takes the form
\begin{equation}\nonumber
 N(x, y) = A\log|x - y| + H(x, y),
 \qquad
 |H(x, y)| \le C,
 \qquad
 |\nabla_y N(x, y)| \le \frac{C}{|x - y|},
\end{equation}
and the same bounds hold uniformly for the rescaled periodic Green function. Whenever $\int f\,dx = 0$, subtracting $A\log\langle x\rangle\int f\,dy = 0$ gives
\begin{equation}\nonumber
 |\mathcal{N} f(x)|
 \le C\int \bigl(1 + \log\langle y\rangle
           + \mathbf{1}_{\{|x - y| < 1\}}|\log|x - y||\bigr)|f(y)|\,dy
 \le C_\delta\bigl(\|f\|_2 + \|W^\delta f\|_1\bigr).
\end{equation}
Similarly, integrating by parts in $\mathcal{N}(\nabla h_1)$ and splitting at $|x - y| = 1$ yields, for $p > 2$,
\begin{equation}\nonumber
 |\mathcal{N}(\nabla h_1)(x)|
 \le C\int \frac{|h_1(y)|}{|x - y|}\,dy
 \le C_p\bigl(\|h_1\|_1 + \|h_1\|_p\bigr).
\end{equation}
Thus,
\begin{equation}\label{vx:kernel-bounds}
 \|\mathcal{N} f\|_\infty
 \le C_\delta\bigl(\|f\|_2 + \|W^\delta f\|_1\bigr) \quad \Bigl(\int f\,dx = 0\Bigr),
 \qquad
 \|\mathcal{N}(\nabla h_1)\|_\infty
 \le C_p\bigl(\|h_1\|_1 + \|h_1\|_p\bigr).
\end{equation}
Furthermore, if $\int f\,dx = 0$, then for $|\xi| \le 1$,
\begin{equation}\nonumber
 |\widehat{f}(\xi)|
 = \left|\int (e^{-ix\cdot\xi} - 1)f(x)\,dx\right|
 \le C|\xi|^\delta \|W^\delta f\|_1,
\end{equation}
so splitting the Fourier integral at $|\xi| = 1$ gives
\begin{equation}\label{vx:fourier-energy}
 \int_{\mathbb{R}^2} \frac{|\widehat{f}(\xi)|^2}{|\xi|^2}\,d\xi
 \le C_\delta\bigl(\|W^\delta f\|_1^2 + \|f\|_2^2\bigr).
\end{equation}

Testing $\rho\dot{u} = \nabla G + \mu\nabla^\perp\omega$ against logarithmic cutoffs (or integrating over a periodic cell) gives $\int \rho\dot{u}\,dx = 0$. Applying \eqref{vx:fourier-energy}, \eqref{vx:force-weight}, and the $L^2$ boundedness of second-order Riesz transforms yields
\begin{equation}\nonumber
 \|\nabla\mathcal{N}(\rho\dot{u})\|_2 + \|\nabla\mathcal{N}(\nabla h_1)\|_2
 \le C\bigl(\|W^\delta\rho\dot{u}\|_1 + \|\rho\dot{u}\|_2 + \|h_1\|_2\bigr)
 \le C.
\end{equation}
Setting $U = u - \mathcal{N}(\rho\dot{u}) - \mathcal{N}(\nabla h_1)$, we see that
\begin{equation}\nonumber
 \mu\Delta U + (\mu + \lambda(0))\nabla\operatorname{div} U = 0,
 \qquad
 \nabla U \in L^2,
\end{equation}
so $U \equiv c_0$ for some constant vector $c_0$. By \eqref{vx:kernel-bounds} and energy conservation,
\begin{equation}\nonumber
 \|u - c_0\|_\infty \le C,
 \qquad
 |c_0| \le m_0^{-1/2}\|\sqrt{\rho}\,u\|_2 + \|u - c_0\|_\infty \le C,
\end{equation}
which gives $\|u\|_\infty \le C$ and proves \eqref{vx:vel}.

Finally, set
\begin{equation}\nonumber
 f = \Delta K(\theta) = c_v\rho\dot{\theta} + P\operatorname{div} u - Q.
\end{equation}
Since $\rho$ is bounded and $\theta$ is bounded away from zero, we have
\begin{equation}\nonumber
 \|\sqrt{\rho/\kappa(\theta)}\|_\infty + \|W^\delta\sqrt{\rho/\kappa(\theta)}\|_2 \le C.
\end{equation}
Applying Cauchy--Schwarz together with \eqref{vx:weighted-pressure}, \eqref{vx:weighted-gradient}, and \eqref{vx:vel}, we deduce
\begin{equation}\label{vx:heat-source}
\begin{aligned}
 \|f\|_2 + \|W^\delta f\|_1
 &\le C\bigl(\|\sqrt{\rho/\kappa(\theta)}\|_\infty + \|W^\delta\sqrt{\rho/\kappa(\theta)}\|_2\bigr)\sqrt{B} \\
 &\quad + C\bigl(\|P\|_4 + \|\nabla u\|_4\bigr)\|\nabla u\|_4
        + C\bigl(\|W^\delta P\|_2 + \|W^\delta\nabla u\|_2\bigr)\|\nabla u\|_2 \\
& \le C(1 + \sqrt{B}).
\end{aligned}
\end{equation}
Testing $\Delta K(\theta) = f$ against logarithmic cutoffs $\chi_L$ (with $\chi_L = 1$ on $B_L$ and $\chi_L = 0$ outside $B_{L^2}$) gives
\begin{equation}\nonumber
 \left|\int \chi_L f\,dx\right|
 = \left|\int \nabla\chi_L \cdot \nabla K(\theta)\,dx\right|
 \le C(\log L)^{-1/2}\|\nabla K(\theta)\|_2 \longrightarrow 0
 \qquad (L \to \infty),
\end{equation}
whence $\int f\,dx = 0$. Let $\mathcal{N}_0 = \Delta^{-1}$ be the scalar Newton potential, for which \eqref{vx:kernel-bounds} and \eqref{vx:fourier-energy} hold equally. Then $\Delta(K(\theta) - \mathcal{N}_0 f) = 0$ with $\nabla(K(\theta) - \mathcal{N}_0 f) \in L^2$, so $K(\theta) - \mathcal{N}_0 f$ is constant. Averaging against $\rho\,dx/m_0$ gives
\begin{equation}\nonumber
 K(\theta) - \frac{1}{m_0}\int \rho K(\theta)\,dx
 = \mathcal{N}_0 f - \frac{1}{m_0}\int \rho\,\mathcal{N}_0 f\,dx,
\end{equation}
where H\"older's inequality and the bound on $H$ imply
\begin{equation}\nonumber
 \frac{1}{m_0}\int \rho K(\theta)\,dx
 \le C\left(\frac{1}{m_0}\int \rho\theta^{2+\alpha}\,dx\right)^{m/(m+1)}
 \le C.
\end{equation}
Combining this estimate with \eqref{vx:kernel-bounds} and \eqref{vx:heat-source}, we arrive at
\begin{equation}\nonumber
 \frac{\kappa_0}{m}\|\theta\|_\infty^m
 = \|K(\theta)\|_\infty
 \le 2\|\mathcal{N}_0 f\|_\infty + \frac{1}{m_0}\int \rho K(\theta)\,dx
 \le C(1 + \sqrt{B}),
\end{equation}
together with
\begin{equation}\nonumber
 \|\nabla^2 K(\theta)\|_2 = \|\Delta K(\theta)\|_2 = \|f\|_2 \le C(1 + \sqrt{B}).
\end{equation}
This establishes \eqref{vx:heat} and completes the proof of Lemma~\ref{vx:ell}.
\end{proof}

\begin{lemma}\label{vx:density}
There exists a positive constant $C$ such that
\begin{equation}\label{vx:weighted-density}
 \sup_{0 \le t \le T} \|\bar{x}^a\rho(t)\|_{L^1 \cap H^1 \cap W^{1,q}}
 + \int_0^T \bigl(\|\nabla u\|_\infty^2 + \|\nabla^2 u\|_2^2 + \|\nabla^2 u\|_q^2\bigr)\,dt
 \le C.
\end{equation}
\end{lemma}

\begin{proof}
Set $Y = e + \|\nabla\rho\|_2 + \|\nabla\rho\|_q$ and $A = 1 + \|\nabla\dot{u}\|_2$. By \eqref{vac:sob}, \eqref{vac:density-weight}, and the elliptic equations for $G$ and $\omega$, we have for $p = 2, 4, q$,
\begin{equation}\nonumber
 \|\nabla(G, \omega)\|_p \le C_{p}\|\rho\dot{u}\|_p \le C_{p} A.
\end{equation}
Combining $\|(G, \omega)\|_2 \le C$ with the Gagliardo--Nirenberg inequality yields $\|(G, \omega)\|_\infty \le C(1 + A^{2/3})$, while \eqref{vx:heat} gives $\|\theta\|_\infty \le C(1 + B^{1/(2m)})$ with $m = 1 + \alpha$. From $h(\rho)\operatorname{div} u = G + P$ and the density upper bound, it follows that
\begin{equation}\nonumber
 \|\operatorname{div} u\|_\infty + \|\omega\|_\infty + \|\theta\|_\infty \le F,
 \qquad
 F = C\bigl(1 + A^{2/3} + B^{1/(2m)}\bigr),
\end{equation}
where Proposition~\ref{vh:cap} ensures $\|F\|_{L^{2+\varepsilon}(0,T)} \le C$ for any $0 < \varepsilon < \min\{1, 2\alpha\}$.

Next, since $\theta$ is bounded away from zero and $K'(\theta) = \kappa(\theta)$, we have $\|\nabla\theta\|_p \le C\|\nabla K(\theta)\|_p$. By Gagliardo--Nirenberg interpolation,
\begin{equation}\nonumber
 \|\nabla K(\theta)\|_q
 \le C\|\nabla K(\theta)\|_2^{2/q} \|\nabla^2 K(\theta)\|_2^{1 - 2/q},
\end{equation}
which together with \eqref{vx:heat} and Proposition~\ref{vh:cap} shows that
\begin{equation}\nonumber
 D_0 = C\bigl(1 + A + \|\nabla\theta\|_2 + \|\nabla\theta\|_q\bigr)
 \quad\text{satisfies}\quad
 \|D_0\|_{L^2(0,T)} \le C.
\end{equation}
Differentiating $h(\rho)\operatorname{div} u = G + \rho\theta$ gives
\begin{equation}\nonumber
 h(\rho)\nabla\operatorname{div} u
 = \nabla G + \theta\nabla\rho
   + \rho\nabla\theta - h'(\rho)\operatorname{div} u\,\nabla\rho.
\end{equation}
Since $h \ge 2\mu$ and $\rho$ is bounded, the div--curl estimate yields
\begin{equation}\label{vx:density-hessian}
 \|\nabla^2 u\|_2 + \|\nabla^2 u\|_q \le C(D_0 + F Y).
\end{equation}
Applying the Beale--Kato--Majda logarithmic inequality
\begin{equation}\nonumber
 \|\nabla u\|_\infty
 \le C_q(1 + \|\nabla u\|_2)
     + C_q\bigl(1 + \|\operatorname{div} u\|_\infty + \|\omega\|_\infty\bigr)\log(e + \|\nabla^2 u\|_q)
\end{equation}
together with $\|\nabla u\|_2 \le C$ and \eqref{vx:density-hessian}, we obtain
\begin{equation}\label{vx:density-log}
 \|\nabla u\|_\infty
 \le C(1 + F)\bigl[\log(e + D_0 + F) + \log Y\bigr].
\end{equation}

Differentiating the continuity equation yields, for $p = 2, q$,
\begin{equation}\nonumber
 \frac{d}{dt}\|\nabla\rho\|_p
 \le C\|\nabla u\|_\infty\|\nabla\rho\|_p + C\|\rho\|_\infty\|\nabla\operatorname{div} u\|_p,
\end{equation}
so $Y' \le C\|\nabla u\|_\infty Y + C_{T}(D_0 + F Y)$ by \eqref{vx:density-hessian}. Dividing by $Y$ and substituting \eqref{vx:density-log} gives
\begin{equation}\nonumber
 (\log Y)'
 \le C_{T}(1 + F)\log Y
     + C_{T}(1 + F)\log(e + D_0 + F) + C_{T}(D_0 + F).
\end{equation}
Since $F \in L^{2+\varepsilon}(0,T)$ and $D_0 \in L^2(0,T)$, the right-hand coefficients are integrable on $(0,T)$, and Gr\"onwall's inequality yields
$$\sup_{0 \le t \le T} Y(t) \le C_{T}.$$
Returning to \eqref{vx:density-hessian} and \eqref{vx:density-log}, we deduce
\begin{equation}\nonumber
 \int_0^T \bigl(\|\nabla u\|_\infty^2 + \|\nabla^2 u\|_2^2 + \|\nabla^2 u\|_q^2\bigr)\,dt \le C_{T}.
\end{equation}

It remains to propagate the spatial weight. Set $b = \bar{x}^a\rho$ and $\gamma = a\log\bar{x}$, noting that $\|\nabla\gamma\|_\infty + \|\nabla^2\gamma\|_\infty \le C$. The continuity equation becomes
\begin{equation}\nonumber
 D_t b = b(u \cdot \nabla\gamma - \operatorname{div} u),
\end{equation}
which gives, for $p = 1, 2, q, \infty$,
\begin{equation}\nonumber
 \|b(t)\|_p
 \le \|b(0)\|_p \exp\left(C\int_0^t \bigl(\|u\|_\infty + \|\operatorname{div} u\|_\infty\bigr)\,ds\right).
\end{equation}
By \eqref{vx:vel} and $\|\operatorname{div} u\|_\infty \le F \in L^1(0,T)$, we obtain
\begin{equation}\label{vx:weight-Lp}
 \sup_{0 \le t \le T}
 \bigl(\|b(t)\|_1 + \|b(t)\|_2 + \|b(t)\|_q + \|b(t)\|_\infty\bigr)
 \le C_{T}.
\end{equation}
Differentiating the transport equation for $b$ yields
\begin{equation}\nonumber
 D_t\nabla b
 = -(\nabla u)^\top\nabla b
   + (u \cdot \nabla\gamma - \operatorname{div} u)\nabla b
   + b\nabla(u \cdot \nabla\gamma - \operatorname{div} u).
\end{equation}
For $p = 2, q$, combining this identity with \eqref{vx:weight-Lp} and $\|u\|_\infty + \|\nabla u\|_p \le C_{T,p}$ from \eqref{vx:vel} gives
\begin{equation}\nonumber
 \frac{d}{dt}\|\nabla b\|_p
 \le C_{T}(1 + \|\nabla u\|_\infty)\|\nabla b\|_p
     + C_{T}\bigl(1 + \|\nabla^2 u\|_p\bigr).
\end{equation}
Since the initial norms are finite by \eqref{main:data} and both time-dependent coefficients are integrable on $(0,T)$, Gr\"onwall's inequality implies
$$\sup_{0 \le t \le T}(\|\nabla b(t)\|_2 + \|\nabla b(t)\|_q) \le C_{T}.$$
Together with \eqref{vx:weight-Lp} and the preceding velocity estimates, this establishes \eqref{vx:weighted-density} and completes the proof of Lemma~\ref{vx:density}.
\end{proof}

\begin{lemma}\label{vt:pos}
For each $0 < \tau < T$, one has
\begin{equation}\nonumber
 \sup_{\tau \le t \le T} \bigl(B(t) + \|\nabla^2 K(\theta(t))\|_2^2\bigr)
 + \int_\tau^T \|\nabla D_t K(\theta)\|_2^2\,dt
 \le C_{\tau},
 \qquad
 \int_0^T \|\nabla^2 K(\theta)\|_2^2\,dt \le C.
\end{equation}
\end{lemma}

\begin{proof}
Set
\begin{equation}\nonumber
 k = K(\theta),
 \qquad V = D_t k,
 \qquad b = \frac{\rho}{\kappa(\theta)},
 \qquad \mathcal{S} = 2\mu D(u) + \lambda(\rho)(\operatorname{div} u) I,
\end{equation}
so that the temperature equation reads
\begin{equation}\label{vt:heat-V}
 c_v b V = \Delta k + Q - P\operatorname{div} u,
 \qquad
 B = \int b V^2\,dx.
\end{equation}
By Proposition~\ref{vh:cap} and the bounds on $\rho$ and $\theta^{-1}$, we have
\begin{equation}\nonumber
 0 \le b \le C\rho,
 \qquad
 \left|\frac{\kappa'(\theta)}{\kappa(\theta)^2}\right|
 + |\lambda(\rho)| + |\rho\lambda'(\rho)| \le C,
\end{equation}
while Lemma~\ref{vx:ell} gives, for some sufficiently small $\delta > 0$,
\begin{equation}\label{vt:mult-input}
 \|\nabla u\|_4 + \|P\|_4 + \|W^\delta\nabla u\|_2 + \|W^\delta P\|_2
 \le C_{T}.
\end{equation}
From \eqref{vt:heat-V}, \eqref{vt:mult-input}, and the lower bound on $\theta$, we deduce
\begin{equation}\label{vt:lap}
 \|\nabla^2 k\|_2
 = \|\Delta k\|_2
 \le C\|\sqrt{\rho/\kappa(\theta)}\|_\infty\sqrt{B}
     + \|Q\|_2 + \|P\operatorname{div} u\|_2
 \le C(1 + \sqrt{B}),
\end{equation}
which together with $\|\nabla k\|_2 \le C$ and $\int_0^T B\,dt \le C_{T}$ yields
\begin{equation}\nonumber
 \|\nabla k\|_4^2 \le C\|\nabla k\|_2\|\nabla^2 k\|_2 \le C(1 + \sqrt{B}),
 \qquad
 \int_0^T \|\nabla^2 k\|_2^2\,dt \le C_{T}.
\end{equation}

Next, put $H_V = \sqrt{B} + \|\nabla V\|_2$. On $\Gamma_t$, the bounds on $\rho$ and $\theta$ give $b \ge c_{T} > 0$ and $\|V\|_{L^2(\Gamma_t)} \le C\sqrt{B}$. Choosing $p_0 > 1$ close enough to $1$ so that $c(p_0 - 1) < a_0$ (with $p_0' = p_0/(p_0 - 1)$ and $w = W^{-c}$), we obtain from \eqref{vac:entropy} that
\begin{equation}\nonumber
 \|b W^{c/p_0'}\|_{p_0}^{p_0}
 \le C \int \rho^{p_0} W^{c(p_0 - 1)}\,dx
 \le C \int \rho W^{a_0}\,dx
 \le C.
\end{equation}
Hence H\"older's inequality and \eqref{vac:sob} imply, for every $j \ge 2$,
\begin{equation}\label{vt:V-weighted}
 \|V\|_{L^j(b\,dx)}
 \le C\|V\|_{L^{jp_0'}(w\,dx)}
 \le C\sqrt{j}\,\bigl(\|\nabla V\|_2 + \|V\|_{L^2(\Gamma_t)}\bigr)
 \le C\sqrt{j}\,H_V.
\end{equation}
Set $Y = |\nabla u| + P$. Choose $s \in (2, 4)$ sufficiently close to $2$, put $t = 2s/(s - 2)$, and pick $\eta > 0$ such that
\begin{equation}\nonumber
 0 < \eta < \delta\left(\frac{4}{s} - 1\right),
 \qquad \eta t > 2.
\end{equation}
Interpolating \eqref{vt:mult-input} gives $\|W^\eta Y\|_s \le C$, while \eqref{vac:sob} with the integrable weight $W^{-\eta t}$ gives $\|W^{-\eta} V\|_t \le C H_V$. Consequently,
\begin{equation}\label{vt:YV}
 \|Y V\|_2 \le \|W^\eta Y\|_s \|W^{-\eta} V\|_t \le C H_V.
\end{equation}

We now differentiate \eqref{vt:heat-V} along the flow. With $\mathcal{A} = \nabla u + (\nabla u)^\top - (\operatorname{div} u) I$, direct computation yields
\begin{equation}\label{vt:identities}
\begin{aligned}
 (D_t + \operatorname{div} u) b
 &= -b\frac{\kappa'(\theta)}{\kappa(\theta)^2} V, \\
 (D_t + \operatorname{div} u)\Delta k
 &= \Delta V - \operatorname{div}(\mathcal{A}\nabla k), \\
 (D_t + \operatorname{div} u)(Q - P\operatorname{div} u)
 &= (2\mathcal{S} - P I) : \nabla\dot{u}
    - b(\operatorname{div} u) V + E,
\end{aligned}
\end{equation}
where
\begin{equation}\nonumber
 E = (\operatorname{div} u) Q
     - 2\mathcal{S} : (\nabla u)^2
     - \rho\lambda'(\rho)(\operatorname{div} u)^3
     + P\operatorname{tr}\bigl((\nabla u)^2\bigr),
\end{equation}
which satisfies $|E| \le C Y|\nabla u|^2$ and $|2\mathcal{S} - P I| \le C Y$. Applying $D_t + \operatorname{div} u$ to \eqref{vt:heat-V}, testing by $V$, and integrating by parts, we obtain
\begin{equation}\label{vt:identity}
\begin{aligned}
 \frac{c_v}{2} B' + \|\nabla V\|_2^2
 &= \frac{c_v}{2}\int b\frac{\kappa'(\theta)}{\kappa(\theta)^2} V^3\,dx
    + \int (2\mathcal{S} - P I) : \nabla\dot{u}\,V\,dx \\
 &\quad - \int b(\operatorname{div} u) V^2\,dx
    + \int E V\,dx
    + \int \mathcal{A}\nabla k \cdot \nabla V\,dx.
\end{aligned}
\end{equation}
By \eqref{vt:YV}, \eqref{vt:mult-input}, \eqref{vt:lap}, and Young's inequality, the last four integrals are estimated by
\begin{equation}\nonumber
\begin{aligned}
 & C\bigl(\|\nabla\dot{u}\|_2 + \sqrt{B} + \|\nabla u\|_4^2\bigr) H_V
   + C\|\nabla u\|_4\|\nabla k\|_4\|\nabla V\|_2 \\
 &\quad \le \varepsilon\|\nabla V\|_2^2
          + C_\varepsilon\bigl(1 + B + \|\nabla\dot{u}\|_2^2\bigr).
\end{aligned}
\end{equation}
For the cubic term, interpolating between $L^2(b\,dx)$ and $L^j(b\,dx)$ ($j \ge 8$) via \eqref{vt:V-weighted} and applying Young's inequality gives
\begin{equation}\nonumber
 \int b|V|^3\,dx
 \le B^{\frac{j-3}{j-2}}\|V\|_{L^j(b\,dx)}^{\frac{j}{j-2}}
 \le C B^{\frac{j-3}{j-2}}\bigl(\sqrt{j}\,H_V\bigr)^{\frac{j}{j-2}}
 \le \varepsilon\|\nabla V\|_2^2
     + C_\varepsilon\left(j^{\frac{j}{2(j-2)}} B^{3/2}
                          + j^{\frac{j}{j-4}} B^{\frac{2(j-3)}{j-4}}\right).
\end{equation}
Choosing $j = 8 + \log(e + B)$, so that $B^{2/(j-4)} + j^{4/(j-4)} \le C$, we deduce
\begin{equation}\nonumber
 \int b|V|^3\,dx
 \le \varepsilon\|\nabla V\|_2^2 + C_\varepsilon(e + B)^2\log(e + B).
\end{equation}
Substituting these bounds into \eqref{vt:identity} and taking $\varepsilon > 0$ sufficiently small yields
\begin{equation}\label{vt:B-differential}
 B' + c_1\|\nabla V\|_2^2
 \le C\bigl(1 + \|\nabla\dot{u}\|_2^2\bigr) + C(e + B)^2\log(e + B).
\end{equation}

Finally, since $\int_0^T(B + \|\nabla\dot{u}\|_2^2)\,dt \le C_{T}$ by Proposition~\ref{vh:cap}, there exists $t_* \in (\tau/2, \tau)$ such that
\begin{equation}\nonumber
 B(t_*) \le \frac{2}{\tau}\int_0^\tau B(t)\,dt \le C_{\tau}.
\end{equation}
Dividing \eqref{vt:B-differential} by $e + B$ gives
\begin{equation}\nonumber
 \bigl(\log(e + B)\bigr)'
 \le C\bigl(1 + \|\nabla\dot{u}\|_2^2\bigr) + C(e + B)\log(e + B),
\end{equation}
where both $e + B$ and $\|\nabla\dot{u}\|_2^2$ are integrable on $[t_*, T]$. Gr\"onwall's inequality therefore yields
$$\sup_{\tau \le t \le T} B(t) \le C_{\tau}.$$
Integrating \eqref{vt:B-differential} over $[\tau, T]$ and applying \eqref{vt:lap} give
\begin{equation}\nonumber
 \sup_{\tau \le t \le T} \|\nabla^2 k\|_2^2
 + \int_\tau^T \|\nabla V\|_2^2\,dt
 \le C_{\tau}.
\end{equation}
This completes the proof of Lemma~\ref{vt:pos}.
\end{proof}

\begin{lemma}\label{vt:trace}
Let $z$ be the initial temperature of a smooth approximation. Then
\begin{equation}\nonumber
 \int \rho(t)|\theta(t) - z|^2\,dx \le C t,
 \qquad 0 \le t \le T.
\end{equation}
\end{lemma}

\begin{proof}
Set $Y(t) = \left(\int \rho|\theta - z|^2\,dx\right)^{1/2}$. By the continuity equation and H\"older's inequality,
\begin{equation}\nonumber
 \frac{1}{2}\frac{d}{dt} Y^2
 = \int \rho(\theta - z)(\dot{\theta} - u \cdot \nabla z)\,dx
 \le Y\bigl(\|\sqrt{\rho}\,\dot{\theta}\|_2 + \|\sqrt{\rho}\,u \cdot \nabla z\|_2\bigr).
\end{equation}
Integrating over $(0,t)$ and applying $Y(0) = 0$, we arrive at
\begin{equation}\nonumber
 Y(t) \le \int_0^t \bigl(\|\sqrt{\rho}\,\dot{\theta}\|_2 + \|\sqrt{\rho}\,u \cdot \nabla z\|_2\bigr)\,ds.
\end{equation}
By Proposition~\ref{vh:cap} and \eqref{vx:vel}, the positive lower bound on $\theta$ gives
\begin{equation}\nonumber
 \|\sqrt{\rho}\,\dot{\theta}\|_2 \le C\sqrt{B},
 \qquad
 \|\sqrt{\rho}\,u \cdot \nabla z\|_2
 \le \|\rho\|_\infty^{1/2}\|u\|_\infty\|\nabla z\|_2
 \le C.
\end{equation}
Applying H\"older's inequality together with $\int_0^T B(s)\,ds \le C$, we deduce
\begin{equation}\nonumber
 Y(t)
 \le C\sqrt{t}\left(\int_0^t B(s)\,ds\right)^{1/2} + C t
 \le C(\sqrt{t} + t).
\end{equation}
Squaring this bound and using $t \le T$ yields
$$\int \rho(t)|\theta(t) - z|^2\,dx \le C t.$$
This completes the proof of Lemma~\ref{vt:trace}.
\end{proof}

\section{Proof of Theorem~\ref{main:global}}\label{vx:construct}

With the a priori estimates established in Section 3 at hand, we are now in a position to prove the main result.

\begin{proof}[Proof of Theorem~\ref{main:global}]
Let $T^* \in (0,\infty]$ be the maximal existence time of the local strong solution given by Lemma~\ref{pre:local}. Fix $T_0 > 0$ and set $T_1 = \min\{T^*, T_0\}$. For every $0 < T < T_1$, the a priori estimates of Section~\ref{sec:apriori} give
\begin{equation}\label{vx:global-uniform}
\begin{aligned}
 &\sup_{0 \le t \le T}\Bigl(
   \|\rho\|_\infty + \|\theta^{-1}\|_\infty
   + \| \nabla u \|_2 + \| \sqrt{\rho} \dot{u} \|_2 + \| \nabla K(\theta) \|_2
   + \|\bar{x}^a\rho\|_{L^1 \cap H^1 \cap W^{1,q}}
   + \|u\|_\infty\Bigr) \\
 &\quad + \int_0^T \Bigl(
   \| \sqrt{\rho \kappa(\theta)} \dot{\theta} \|^2_2 + \|\nabla\dot{u}\|_2^2
   + \|\nabla u\|_\infty^2
   + \|\nabla^2 u\|_2^2 + \|\nabla^2 u\|_q^2\Bigr)\,dt
  \le C_{T_0}.
\end{aligned}
\end{equation}
For every $0 < \tau < T_1$, Lemmas~\ref{vx:ell} and~\ref{vt:pos} also yield
\begin{equation}\label{vx:positive-time-uniform}
 \sup_{\tau \le t < T_1}\Bigl(
   \|\theta\|_\infty + \| \sqrt{\rho \kappa(\theta)} \dot{\theta} \|^2_2 + \|\nabla^2 K(\theta)\|_2^2\Bigr)
 + \int_\tau^{T_1} \|\nabla D_t K(\theta)\|_2^2\,dt
 \le C_{\tau,T_0}.
\end{equation}

We next bound $\nabla^2 u$ for positive times. By \eqref{vx:global-uniform}, \eqref{vx:positive-time-uniform}, and $\nabla\theta = \nabla K(\theta)/\kappa(\theta)$, we have for a.e.\ $t \in [\tau, T_1)$,
\begin{equation}\nonumber
 \|\nabla(G, \omega)\|_2 \le C\|\rho\dot{u}\|_2 \le C A_2,
 \qquad
 \|\nabla P\|_2
 \le C\bigl(\|\theta\|_\infty\|\nabla\rho\|_2
            + \|\rho\|_\infty\|\nabla\theta\|_2\bigr)
 \le C_{\tau,T_0}.
\end{equation}
Differentiating $h(\rho)\operatorname{div} u = G + P$ and applying the div--curl estimate together with \eqref{vx:vel}, we obtain
\begin{equation}\nonumber
 \|\nabla^2 u\|_2
 \le C\left(\|\nabla(G, \omega)\|_2 + \|\nabla P\|_2
            + \|\nabla\rho\|_q \|\nabla u\|_{\frac{2q}{q-2}}\right)
 \le C_{\tau,T_0}.
\end{equation}
Hence,
\begin{equation}\label{vx:restart-bounds}
 \sup_{\tau \le t < T_1}\Bigl(
   \|\bar{x}^a\rho\|_{L^1 \cap H^1 \cap W^{1,q}}
   + \|\nabla u\|_{H^1}
   + \|\nabla K(\theta)\|_{H^1}
   + \|\theta\|_\infty + \|\theta^{-1}\|_\infty
   + \|\sqrt{\rho}\,\dot{u}\|_2\Bigr)
 \le C_{\tau,T_0}.
\end{equation}

We claim that $T^* = \infty$. Suppose otherwise that $T^* < \infty$, and take $T_0 = T^* + 1$ (so $T_1 = T^*$) and $\tau = T^*/2$. For a.e.\ $t_0 \in [\tau, T^*)$, the momentum equation gives
\begin{equation}\nonumber
 -\mu\Delta u(t_0)
 - \nabla\bigl((\mu + \lambda(\rho(t_0)))\operatorname{div} u(t_0)\bigr)
 + \nabla P(t_0)
 = \sqrt{\rho(t_0)}\,g_1(t_0),
\end{equation}
where $g_1(t_0) = -\sqrt{\rho(t_0)}\,\dot{u}(t_0)$ obeys $\|g_1(t_0)\|_2 \le C_{T_0}$. Together with \eqref{vx:restart-bounds}, $(\rho, u, \theta)(t_0)$ satisfies the hypotheses of Lemma~\ref{pre:local} with norms bounded uniformly in $t_0$. Lemma~\ref{pre:local} therefore yields a local strong solution starting from $t_0$ on $[t_0, t_0 + \delta]$, where $\delta > 0$ depends only on $T_0$ and $\tau$. Picking $t_0 \in [\tau, T^*)$ with $t_0 + \delta > T^*$ extends the solution past $T^*$ by local uniqueness, a contradiction. Thus $T^* = \infty$.

For any $T > 0$, \eqref{vx:global-uniform}--\eqref{vx:restart-bounds}, Lemma~\ref{vx:density}, and Lemma~\ref{vt:pos} imply \eqref{main:finite}; in particular, the continuity equation and \eqref{vx:vel} give
\begin{equation}\nonumber
 \|\rho_t\|_q
 \le \|u\|_\infty\|\nabla\rho\|_q + \|\rho\|_\infty\|\operatorname{div} u\|_q
 \le C_T.
\end{equation}
By the embedding $W^{1,q}(\mathbb{R}^2) \hookrightarrow L^\infty(\mathbb{R}^2)$ ($q > 2$) and \eqref{vx:global-uniform}, we have
\begin{equation}\nonumber
 0 \le \rho(t,x) \le \|\bar{x}^a\rho(t)\|_\infty \bar{x}^{-a} \le C_T \bar{x}^{-a},
\end{equation}
so $\rho(t,x) \to 0$ as $|x| \to \infty$. The initial traces of $u$ and $\theta$ are given directly by Lemma~\ref{pre:local} and Lemma~\ref{vt:trace}, respectively. Next, fix $\varphi \in C_c^\infty(\mathbb{R}^2)$ with $\varphi = 1$ on $B_1$ and $\varphi = 0$ outside $B_2$, and put $\varphi_L(x) = \varphi(x/L)$. Multiplying the conservative system \eqref{intro:conservative} by $\varphi_L$, integrating over $\mathbb{R}^2 \times (0,t)$, and passing to the limit $L \to \infty$, all integrals involving $\nabla\varphi_L$ vanish in view of \eqref{vx:global-uniform}. Consequently, for every $t \ge 0$,
\begin{equation}\nonumber
 \int \rho(t)\,dx = m_0,
 \qquad
 \int (\rho u)(t)\,dx = \int \rho_0 u_0\,dx,
 \qquad
 \int \left(\frac{1}{2}\rho(t)|u(t)|^2 + c_v\rho(t)\theta(t)\right)dx = E_0.
\end{equation}
The uniqueness of $(\rho,u,\theta)$ follows from arguments similar to those in the proof of \cite[Theorem 1]{ChoKim2006}.
This completes the proof of Theorem~\ref{main:global}.
\end{proof}

\bigskip

\noindent\textbf{Data availability.} No data was used for the research described in the article.

\bigskip

\noindent\textbf{Conflict of interest.} The authors declare that they have no conflict of interest.

\section*{Acknowledgments}
The authors acknowledge the use of generative AI tools in developing and checking mathematical arguments and in drafting and revising this manuscript.

\end{document}